\documentclass[letterpaper,11pt]{article}

\usepackage[margin=1in]{geometry} 
\usepackage{xspace,exscale,relsize}
\usepackage{fancybox,shadow}
\usepackage{graphicx}
\usepackage{color}

\usepackage{amsthm,amsfonts}
\usepackage{amssymb}
\usepackage{authblk}
\usepackage[title]{appendix}

\usepackage{extarrows}

\usepackage[noadjust]{cite}

\usepackage{amsmath}
	\makeatletter
	\let\over=\@@over \let\overwithdelims=\@@overwithdelims
	\let\atop=\@@atop \let\atopwithdelims=\@@atopwithdelims
  	\let\above=\@@above \let\abovewithdelims=\@@abovewithdelims
  	\makeatother
\allowdisplaybreaks

\usepackage{rotating}

\usepackage{booktabs}
\usepackage{array}
\usepackage{ifpdf}
\usepackage{float}
\usepackage{subfigure}
\usepackage{psfrag}

\usepackage{prettyref}
\usepackage{enumitem}

\usepackage{graphicx}
\usepackage{subcaption}

\usepackage{tikz}
\usetikzlibrary{arrows}
\tikzstyle{int}=[draw, fill=blue!20, minimum size=2em]
\tikzstyle{dot}=[circle, draw, fill=blue!20, minimum size=2em]
\tikzstyle{init} = [pin edge={to-,thin,black}]

\usepackage[
            CJKbookmarks=true,
            bookmarksnumbered=true,
            bookmarksopen=true,
            colorlinks=true,
            citecolor=red,
            linkcolor=blue,
            anchorcolor=red,
            urlcolor=blue,
            pdfauthor={Polyanskiy and Wu}
            ]{hyperref}

\usepackage[all]{xy}

\usepackage{mathtools}

\usepackage{ifthen}
\newboolean{aos}
\setboolean{aos}{FALSE}

\ifx\eqref\undefined
	\newcommand{\eqref}[1]{~(\ref{#1})}
\fi
\ifx\mod\undefined
	\def\mod{\mathop{\rm mod}}
\fi

\usepackage{bm}

\newcommand{\norm}[1]{{\left\Vert #1 \right\Vert}}

\def\argmax{\mathop{\rm argmax}}

\def\exp{\mathop{\rm exp}}

\def\Var{\mathrm{Var}}
\def\Cov{\mathrm{Cov}}

\def\simiid{\stackrel{\mathrm{iid}}{\sim}}

\newcommand{\id}{\text{id}}

\newcommand{\Unif}{\mathrm{Uniform}}

\newcommand{\reals}{\mathbb{R}}

\newcommand{\Expect}{\mathbb{E}}
\newcommand{\expect}[1]{\mathbb{E}\left[#1\right]}

\newcommand{\prob}[1]{\mathbb{P}\left[#1\right]}

\newcommand{\TV}{{\rm TV}}

\newcommand{\diff}{{\rm d}}

\newcommand{\pth}[1]{\left( #1 \right)}
\newcommand{\qth}[1]{\left[ #1 \right]}

\newcommand{\iiddistr}{\simiid}

\newcommand{\var}{\Var}
\newcommand\indep{\protect\mathpalette{\protect\independenT}{\perp}}
\def\independenT#1#2{\mathrel{\rlap{$#1#2$}\mkern2mu{#1#2}}}

\newcommand{\Binom}{\text{Binom}}

\newcommand{\Iprod}[2]{\langle #1, #2 \rangle}
\newcommand{\iprod}[2]{\langle #1, #2 \rangle}
\newcommand{\indc}[1]{{1\left\{{#1}\right\}}}

\definecolor{myblue}{rgb}{.8, .8, 1}
\definecolor{mathblue}{rgb}{0.2472, 0.24, 0.6} 
\definecolor{mathred}{rgb}{0.6, 0.24, 0.442893}
\definecolor{mathyellow}{rgb}{0.6, 0.547014, 0.24}

\newcommand{\calG}{{\mathcal{G}}}
\newcommand{\calH}{{\mathcal{H}}}

\newcommand{\calR}{{\mathcal{R}}}

\newrefformat{eq}{(\ref{#1})}
\newrefformat{thm}{Theorem~\ref{#1}}
\newrefformat{th}{Theorem~\ref{#1}}
\newrefformat{chap}{Chapter~\ref{#1}}
\newrefformat{sec}{Section~\ref{#1}}
\newrefformat{seca}{Section~\ref{#1}}
\newrefformat{algo}{Algorithm~\ref{#1}}
\newrefformat{fig}{Fig.~\ref{#1}}
\newrefformat{tab}{Table~\ref{#1}}
\newrefformat{rmk}{Remark~\ref{#1}}
\newrefformat{clm}{Claim~\ref{#1}}
\newrefformat{def}{Definition~\ref{#1}}
\newrefformat{cor}{Corollary~\ref{#1}}
\newrefformat{lmm}{Lemma~\ref{#1}}
\newrefformat{lem}{Lemma~\ref{#1}}
\newrefformat{prop}{Proposition~\ref{#1}}
\newrefformat{pr}{Proposition~\ref{#1}}
\newrefformat{app}{Appendix~\ref{#1}}
\newrefformat{apx}{Appendix~\ref{#1}}
\newrefformat{ex}{Example~\ref{#1}}
\newrefformat{exer}{Exercise~\ref{#1}}
\newrefformat{soln}{Solution~\ref{#1}}

\def\unifto{\mathop{{\mskip 3mu plus 2mu minus 1mu%
	\setbox0=\hbox{$\mathchar"3221$}%
	\raise.6ex\copy0\kern-\wd0%
	\lower0.5ex\hbox{$\mathchar"3221$}}\mskip 3mu plus 2mu minus 1mu}}

\ifx\lesssim\undefined
\def\simleq{{{\mskip 3mu plus 2mu minus 1mu%
	\setbox0=\hbox{$\mathchar"013C$}%
	\raise.2ex\copy0\kern-\wd0%
	\lower0.9ex\hbox{$\mathchar"0218$}}\mskip 3mu plus 2mu minus 1mu}}
\else
\def\simleq{\lesssim}
\fi

\ifx\gtrsim\undefined
\def\simgeq{{{\mskip 3mu plus 2mu minus 1mu%
	\setbox0=\hbox{$\mathchar"013E$}%
	\raise.2ex\copy0\kern-\wd0%
	\lower0.9ex\hbox{$\mathchar"0218$}}\mskip 3mu plus 2mu minus 1mu}}
\else
\def\simgeq{\gtrsim}
\fi

\newtheorem{theorem}{Theorem}
\newtheorem{lemma}[theorem]{Lemma}

\newtheorem{proposition}[theorem]{Proposition}

\theoremstyle{definition}
\newtheorem{definition}{Definition}

\newtheorem{remark}{Remark}

\newif\ifmapx
{\catcode`/=0 \catcode`\\=12/gdef/mkillslash\#1{#1}}
\edef\jobnametmp{\expandafter\string\csname embayes2_apx\endcsname}
\edef\jobnameapx{\expandafter\mkillslash\jobnametmp}
\edef\jobnameexpand{\jobname}
\ifx\jobnameexpand\jobnameapx
\mapxtrue
\else
\mapxfalse
\fi

\newcommand{\poly}{\mathsf{poly}}
\newcommand{\polylog}{\mathsf{polylog}}

\newcommand{\Law}{\text{Law}}

\renewcommand{\hat}{\widehat}
\renewcommand{\tilde}{\widetilde}
\newcommand{\CMT}{\mathsf{CMT}}
\newcommand{\ASWT}{\mathsf{ASWT}}
\newcommand{\SMT}{\mathsf{SMT}}
\newcommand{\SSMT}{\mathsf{SSMT}}

\title{ 
The broken-sample problem revisited II: \\Detecting hidden linear dependencies}

\author{Shuyang Gong
\footnote{Department of Statistics and Data Science and Institute for Foundations of Data Science, Yale University. Email: \texttt{shuyang.gong@yale.edu}.} \hspace{1cm} Yihong Wu 
\footnote{Department of Statistics and Data Science, Yale University. Email: \texttt{yihong.wu@yale.edu}.} \hspace{1cm} Jiaming Xu  \footnote{The Fuqua School of Business, Duke University. Email: \texttt{jx77@duke.edu}.}}
\begin{document}
\maketitle
\ifpdf
\DeclareGraphicsExtensions{.pgf}
\graphicspath{{figures/}{plots/}}
\fi
\begin{abstract}

We study the detection problem for shuffled linear regression, in which the covariates and the responses are either independent or linearly related through a latent correspondence. Consider independent standard Gaussian covariates \(X_1,\ldots,X_n\in\mathbb R^d\) and scalar responses \(Y_1,\ldots,Y_n\). In a shuffled linear model, \(Y_i=\langle X_{\pi_*(i)},\beta_*\rangle+\sigma Z_i\), where 
both the permutation \(\pi_*\) and the regression coefficient $\beta_*$ are latent and $Z_i$'s are white Gaussian noise.
The goal is to distinguish this model from one in which the covariates and responses are independent with the same marginal distribution.

In the low-dimensional regime of $d \ll \frac{\log n}{\log\log n}$, we show that the detection threshold 
is given by $\sigma=o(1)$ for fixed $d$ and $\sigma^2 = \Theta(\frac{1}{d\log d})$ for growing $d$; both thresholds are achieved in nearly linear time using moment-based test statistics inspired by matched filters in signal processing.
For higher dimensions, tests based on sliced Wasserstein distance succeed when $\sigma^2 = o(\frac{1}{(d\log d)^2})$, while the optimal threshold remains open.

As a by-product, we show that the shuffled linear model exhibits an ``all-or-nothing'' phase transition for estimating $\beta_*$ at $\sigma^2\asymp \frac{1}{d\log d}$ in low dimensions. In particular, $\beta_*$ can be consistently estimated in nearly linear time, even when it is information-theoretically impossible to recover any positive fraction of the latent matching $\pi^*$.
\end{abstract}

\tableofcontents

\section{Introduction}

Suppose two surveys are carried out, each measuring a different attribute of the same subjects, but the records contain no common identifier. We observe both samples without knowing which measurements belong to the same individual. Can we still tell whether the two attributes are correlated? This is the detection question in the classical \emph{broken sample problem},
which grew out of record linkage in survey and census \cite{fellegi1969theory,degroot1971matchmaking,degroot1980estimation,bai2005broken}. More recently, it has garnered renewed interest due to applications in  privacy-preserving data publishing~\cite{muralidhar2006data,xiao2006anatomy,pananjady2017linear}, database alignment \cite{cullina2018fundamental,paslev2023testing}, 
multi-sensor data integration~\cite{kaplan2008assignment,wang2021multisensor}, and multimodal data analysis~\cite{zhu2023robust}.

The mathematical setup of the detection problem is as follows. Suppose we observe two samples of iid data points $X_1,\dots,X_n$ and $Y_1,\dots,Y_n$, but not their correspondence. Under the null hypothesis, the samples are independent and respectively sampled from distributions $\mathcal{P}_X$ and $\mathcal{P}_Y$. Under the alternative, there is an unobserved permutation $\pi_*$ such that, conditional on $\pi_*$, the pairs $(X_{\pi_*(i)},Y_i),i\in [n]$ are iid from a joint distribution $\mathcal{P}_{X,Y}$ with the same marginals. The goal is to detect the dependency between the two samples, and if so, recover their hidden correspondence.

Most existing works assume that the joint distribution $\mathcal{P}_{X,Y}$ is known. Notably, in a companion paper \cite{jiao2025DetectionBroken} of this work, the authors 
establish a sharp detection threshold characterized by both spectral and $L_2$ conditions of the likelihood ratio operator $\diff \mathcal{P}_{XY}/\diff(\mathcal{P}_X \otimes \mathcal{P}_Y)$, resolving a conjecture of Bai and Hsing~\cite{bai2005broken} in the positive. 

In many practical scenarios, however, the joint law $\mathcal{P}_{X,Y}$ is unknown. As a first step toward this more challenging setting, this paper studies the case where  $X$ and $Y$ are related through a linear model with unknown regression coefficients.
We formalize this setting below.

\begin{definition}[Shuffled linear regression]\label{def:slr}
    Suppose we observe two datasets $\{X_1,\ldots,X_n\}$ and $\{Y_1,\ldots,Y_n\}$ distributed according to one of the following hypotheses:
\begin{align}
\mathcal H_{0}: \quad &
  X_{i}\stackrel{\mathrm{iid}}{\sim}N(0,I_d), 
\quad Y_i \stackrel{\mathrm{iid}}{\sim}N(0,1+\sigma^{2}), 
 \nonumber \\
\mathcal H_{1}: \quad &
  X_{i}\stackrel{\mathrm{iid}}{\sim}N(0,I_d),
\quad  Y_i= \iprod{X_{\pi_*(i)}}{\beta_*} + \sigma Z_i 
\label{eq-def-of-the-problem}
\end{align}
where under $\mathcal H_{1}$, $\beta_*$ is uniformly distributed on the unit sphere in $\mathbb R^d$, $\pi_*$ is a uniform random permutation on $[n]$, $\{Z_i\}$ are iid $N(0,1)$, and $(\beta_*,\pi_*,\{Z_i\},\{X_{i}\})$ are mutually independent.
\end{definition}

If $\pi^*$ were known, we could pair each response with its covariate, reducing the model \eqref{eq-def-of-the-problem} to standard linear regression.
If $\beta_*$ were known, we could project each $X_i$ to $\Iprod{X_i}{\beta_*}$ without loss of information, reducing the model to the broken sample problem with a bivariate Gaussian distribution. 
When both $\pi^*$ and $\beta_*$ are unknown, this is the model of \emph{shuffled linear regression} \cite{pananjady2017linear,abid2017linear}, also known as unlabeled sensing \cite{unnikrishnan2018unlabeled,zhang2019permutation,zhang2020optimal}, or linear regression with permuted or mismatched data \cite{slawski2019linear,mazumder2023linear,slawski2020two}. 

Note that both the covariates and the responses have identical marginal distribution under both hypotheses. As such, neither sample alone can distinguish the hypotheses: detection must rely on their dependency under \(\mathcal H_1\) by 
accounting for the hidden permutation $\pi_*$ and regression coefficients $\beta_*$. This makes the problem significantly more challenging than standard broken sample problem, both statistically and computationally.

In this paper, we focus on the task of \textit{strong detection}, namely distinguishing $\mathcal{H}_0$ and $\mathcal{H}_1$ with vanishing type-I and type-II errors as $n\to\infty$. Our main contributions are summarized below:

\begin{itemize}
    \item \textbf{Low-dimensional regime
    ($d=o(\log n/\log\log n)$).}
    For fixed $d$, strong detection is possible if and only if
    $\sigma^2\to0$. For $d\to\infty$, we determine the threshold
    up to constant factors: strong detection is achievable when
    $\sigma^2\le c/(d\log d)$, whereas even weak detection is
    impossible when $\sigma^2\ge C/(d\log d)$, for sufficiently
    small $c>0$ and large $C>0$. Our tests attain the optimal
    noise scale in $n^{1+o(1)}$ time.
    \item \textbf{High-dimensional regime.}
    Strong detection is achievable when
    $\sigma^2=o((d\log d)^{-2})$ and $d=o(n/\log n)$.
    Conversely, for $d\to\infty$ and $d\le n/(\log n)^3$,
    even weak detection is impossible when
    $\sigma^2\ge C/(d\log d)$.
		Determining the sharp detection threshold remains an open question.\footnote{  
        The recent work \cite{Zhangsong25+} provided evidence on a statistical-computational gap in  high dimensions under the low-degree polynomial framework \cite{hopkins2018thesis}. 
        When $d=\omega((\log n)^4)$, it is shown in \cite{Zhangsong25+}  that all polynomials of degree $o(d^{1/4})$ fail to distinguish $\mathcal H_0$ and $\mathcal H_1$, even when $\sigma=0$.}
\end{itemize}

\begin{table}[htbp]
  \centering
  \begin{tabular*}{\linewidth}{@{\extracolsep{\fill}}llll@{}}
    \toprule
     & \textbf{Noise variance} & \textbf{Dimension} & \textbf{Running time} \\
    \midrule
    Cross-moment test
      & $\sigma^2=o\left(\dfrac{1}{d^2}\right)$
      & $d \le (\log n)^{0.3}$ & $n^{1+o(1)}$ \\
    \addlinespace
    Sliced-moment test
      & $\sigma^2=O\left(\dfrac{1}{d\log d}\right)$
      & $d \le (\log n)^{1-o(1)}$ & $n^{1+o(1)}$ \\
    \addlinespace
    Approximate sliced-Wasserstein test
      & $\sigma^2=o\left(\dfrac{1}{(d\log d)^2}\right)$
      & $d \le n^{1-o(1)}$ & $n^{1+o(1)}d^{\Theta(d)}$ \\
    \midrule
    Information-theoretic lower bound
      & $\sigma^2=\Omega\left(\dfrac{1}{d\log d}\right)$
      & $d \leq n^{1-o(1)}$ & \textemdash \\
    \bottomrule
  \end{tabular*}
  \caption{Performance guarantees and time complexity for proposed test statistics, along with information-theoretic lower bound for the impossibility of detection.}
  \label{tab:summary_results}
\end{table}

 Specifically, as summarized in \prettyref{tab:summary_results}, we introduce three complementary test statistics, based on \emph{cross moments}, \emph{sliced moments}, and \emph{approximate sliced-Wasserstein distance}. 
Inspired by \emph{matched filters} in signal processing, the cross-moment test is simple and practical, while the sliced-moment test attains the optimal detection threshold in the low-dimensional regime. The approximate sliced-Wasserstein test delivers the strongest performance in high dimensions, at a higher computational cost.
For low dimensions, all tests run in nearly linear time $n^{1+o(1)}$.

It is instructive to examine how the noise tolerance varies with dimension $d$. When $d$ is fixed, all three tests succeed as $\sigma\to 0$, matching the broken sample setting where $\beta_*$ is known \cite[Corollary 2]{jiao2025DetectionBroken}. In this regime, the unknown $\beta_*$ has negligible effect on the statistical limit. However, as the latent dimension $d$ grows, inferring $\beta_*$ becomes intrinsically harder, and the maximum tolerable noise level decreases. This behavior stands in stark contrast to the $d$-dimensional broken-sample problem, where the entire vector $Y_i=X_{\pi_*(i)}+\sigma Z_i\in\mathbb R^d$ is observed instead of a one-dimensional projection: there, strong detection remains possible even when $\sigma^2$ grows with $d$, as long as $\sigma^2=o(d)$ \cite{elimelech2023phase,jiao2025DetectionBroken}.
Our results are further extended to the sparse setting, where the detection threshold improves to $\sigma^2\asymp (s\log (d/s))^{-1}$, with $s=\|\beta_*\|_0$ denoting the sparsity level of $\beta_*$. 

Closely related to detection is the \emph{estimation} problem in shuffled linear regression, where the goal is to recover $\pi_*$ and $\beta_*$. The existing literature largely focuses on reconstructing the matching $\pi_*$, which amounts to solving an NP-hard quadratic assignment problem based on projection matrices constructed from $\{X_i\}$ and $\{Y_i\}$. 
For regression, however, the parameter of primary interest is the regression coefficient \(\beta_*\). 
A surprising implication of our detection results is that $\beta_*$ can be estimated accurately in nearly linear time, even when it is information-theoretically impossible to match any positive fraction of the data. In particular, we establish that in the low-dimensional regime $d=o(\log n /\log \log n)$,  the estimation of $\beta_*$ undergoes an \emph{all-or-nothing phase transition}:  Consistent estimation is possible when
$\sigma^2 \lesssim (d\log d)^{-1}$, while 
non-trivial estimation of $\beta_*$ becomes impossible when $\sigma^2 \gtrsim (d\log d)^{-1}$.
In comparison, estimating the permutation $\pi_*$ requires far lower noise levels, $\sigma^2=n^{-2+o(1)}$ for almost exact recovery and $\sigma^2=n^{-4+o(1)}$ for exact recovery \cite{wang2022random}.

\paragraph{Notation}
Let $H_k(x)$ be the monic degree-$k$ Hermite polynomial and $h_k(x)\triangleq H_k(x)/\sqrt{k!}$ be the corresponding normalized version, so that $\Expect_{Z\sim N(0,1)}[h_k(Z)h_j(Z)]=\indc{k=j}$. Moreover, we will use $\mathbb P_0,\mathbb P_1$ and $\mathbb E_0,\mathbb E_1$ to denote the probability and expectation under $\mathcal{H}_0,\mathcal H_1$, and for a random variable $W$, we denote by $\mathbb{E}_W$ the expectation taken with respect to $W$ while holding all other variables fixed. We will use both notations interchangeably when there is no ambiguity. Let $\mathfrak{S}_n$ denote the set of all permutations of $[n]$. For $\pi_*\in\mathfrak{S}_n$, let $\Pi_*$ denote the corresponding
permutation matrix. 
We denote the identity permutation by $\mathsf{Id}$ and the identity matrix by $\mathbf{I}_n$.
Let $\|A\|_{\operatorname{op}}$ denote the operator norm of a matrix $A$.

We use standard asymptotic notations: for two sequences $a_n$ and $b_n$ of positive numbers, we write $a_n=O(b_n)$ or $b_n=\Omega(a_n)$ if $a_n\leq Cb_n$ for an absolute constant $C$; we write $a_n=\Theta(b_n)$ or $a_n\asymp b_n$ if both $a_n=O(b_n)$ and $a_n=\Omega(b_n)$ hold; we use $a_n=o(b_n)$ or $a_n\ll b_n$ to denote $a_n/b_n\to 0$, and $a_n=\omega(b_n)$ or $a_n\gg b_n$ to denote $a_n/b_n\to\infty$. 
We write $a_n=\poly(n)$ if $a_n=n^{\Theta(1)}$ and $a_n=\polylog(n)$ if $a_n=\poly(\log n)$. For $a,b\in\mathbb R$, we write $a\wedge b\triangleq\min\{a,b\}$ and $a\vee b\triangleq\max\{a,b\}$.

\paragraph{Organization of this paper}
Section~\ref{sec:main-results} presents the main results of the paper, including the proposed test statistics and their performance guarantees, together with an information-theoretic lower bound.
Section~\ref{sec:ext} determines the sharp threshold for estimating the regression coefficients $\beta_*$ in the low-dimensional regime, and extends the detection results to sparse $\beta_*$.
Section~\ref{sec:numerical-results} contains numerical experiments.
Sections~\ref{sec:integral-over-beta-testing-statistic} and~\ref{sec:remaining-proofs} prove the main results.
Appendix~\ref{sec:computation-CMT} shows how to compute the cross-moment test and bounds its running time. Appendices~\ref{sec:auxiliary-results} and~\ref{sec:normal-approx} collect auxiliary results on Hermite polynomials and normal approximation.

\section{Test Statistics and Lower bounds}

\label{sec:main-results}
In this section, we introduce the test statistics and provide statistical guarantees, complemented by an information-theoretic lower bound. 

\subsection{Hermite expansion of the empirical measure}

We first describe a sequence of reductions of the testing problem and a heuristic that guides the design of our test statistics proposed in Sections \ref{sec:cross-moment}--\ref{sec:ASWT}. Along the way, we compare with existing results and algorithms for the broken sample problem and highlight the added difficulty in shuffled regression due to the latent regression coefficients. 
For simplicity, we will focus on fixed dimensionality and repeatedly appeal to large-sample limiting theorems to gain insights.

Recall that the observations consist of the set of covariates $X=(X_1,\ldots,X_n)$
and the set of responses 
$Y=(Y_1,\ldots,Y_n)$, which are independent under the null $\calH_0$, but are related through a linear model
$Y_i = \Iprod{\beta_*}{X_{\pi_*(i)}} + \sigma Z_i$ under the alternative $\calH_1,$ where $\beta_* \sim \Unif(\mathbb S^{d-1})$ is the latent regression coefficient, $\pi_* \sim \Unif(\mathfrak S_n)$ is the latent permutation, and $Z_i\sim N(0,1)$.  Therefore, the strong detection problem boils down to testing whether $Y$ is independent of $X$, i.e., determining whether $\TV(\Law(Y|X), \Law(Y))\to 1$ in probability with respect to $X$.

As a first step of the reduction, in light of the permutation invariance of the problem, it is sufficient to summarize $Y$ into its empirical distribution:
\[
\pi_Y \triangleq \frac{1}{n} \sum_{i=1}^n \delta_{Y_i} \, .
\]
Furthermore, we may encode this distribution by its action on test functions. In other words, a further sufficient statistic is the collection of  empirical averages 
$\{\pi_Y(g): g \in \calG \}$ for some sufficiently rich class of test functions, where
\[
\pi_Y(g) \triangleq \frac{1}{n} \sum_{i=1}^n g(Y_i).
\]

A commonly used class of test functions consists of polynomials. For convenience, we will choose the basis of Hermite polynomials, normalized to be orthonormal under the marginal distribution $Y_i \sim N(0,1+\sigma^2)$ and further scaled by $\frac{1}{\sqrt{n}}$ to be $O_P(1)$, to obtain the following Hermite coefficients:
\begin{equation}\label{eq-def-Tk}
T_k  \triangleq \frac{1}{\sqrt{n}} \sum_{i=1}^n h_k\pth{\frac{Y_i}{\sqrt{1+\sigma^2}}}.
\end{equation}
Thus, the goal is to distinguish the joint law of the sequence $(T_k: k\geq 1)$ with or without conditioning on $X$ (corresponding to $\calH_1$ or $\calH_0$ respectively).

Without conditioning on $X$, by CLT and thanks to the orthonormality of the Hermite polynomial
\[
\Expect\qth{h_k\pth{\frac{Y_i}{\sqrt{1+\sigma^2}}} h_m\pth{\frac{Y_i}{\sqrt{1+\sigma^2}}}} = \indc{k=m},
\]
as $n\to\infty$, 
the limiting joint distribution of the $T_k$'s is independent standard normal.

However, conditioned on $X$, this limiting distribution no longer holds. To describe the conditional law, let us 
introduce the correlation coefficient $\rho \equiv \frac{1}{\sqrt{1+\sigma^2}}$ and rewrite
\[
T_k  = \frac{1}{\sqrt{n}} \sum_{i=1}^n h_k \left(\rho \Iprod{X_i}{\beta_*} + 
\sqrt{1-\rho^2} Z_i \right).
\]
Note that, \textit{crucially}, the latent permutation $\pi_*$ no longer appears in the empirical average.
Next, we apply the decomposition rule of Hermite polynomials (see \prettyref{lem:scaling-law-hermite} in Appendix \ref{sec:auxiliary-results})
\[
    H_k(\rho x+\sqrt{1-\rho^2}y) = \sum_{\ell=0}^k \binom{k}{\ell}\rho^{\ell}(1-\rho^2)^{(k-\ell)/2}H_\ell(x)H_{k-\ell}(y)
\]
to obtain
 \begin{equation}
     T_k = \underbrace{\rho^k 
     \frac{1}{\sqrt{n}}\sum_{i=1}^n h_k(\langle X_i,\beta_*\rangle)}_{\equiv S_k(\beta_*)} + \underbrace{\frac{1}{\sqrt{n}}\sum_{i=1}^n \frac{1}{\sqrt{k!}}\sum_{\ell=0}^{k-1} \binom{k}{\ell} \rho^\ell (1-\rho^2)^{(k-\ell)/2}H_\ell(\langle X_i,\beta_*\rangle)H_{k-\ell}(Z_i)}_{\equiv U_k}\,. \label{eq:T_k_decomp}
 \end{equation}
  The first term $S_k(\beta_*)$ is the signal,  and the second term $U_k$ turns out to be approximately independent Gaussian noise and independent of $\beta_*$.
To see this, conditioned on $X_i$'s and $\beta_*$, applying CLT to $Z_i\iiddistr  N(0,1)$ and using the orthogonality of Hermite polynomials, we have that
$U_k$ is approximately centered normal with variance
\[
\Var(U_k) = 
\frac{1}{n} \sum_{\ell=0}^{k-1} \frac{(k-\ell)!}{k!} \binom{k}{\ell}^2 \rho^{2\ell} (1-\rho^2)^{k-\ell}
\sum_{i=1}^n H_\ell(\langle X_i,\beta_*\rangle)^2
\approx
 \sum_{\ell=0}^{k-1} \binom{k}{\ell} \rho^{2\ell} (1-\rho^2)^{k-\ell} = 1-\rho^{2k},
\]
where the second step applies LLN to $\Iprod{X_i}{\beta_*}\iiddistr N(0,1)$.
Furthermore, for $k\neq m$, by CLT we have
$\Cov(U_k,U_m) = O(\frac{1}{\sqrt{n}})$ in probability with respect to $X$. 
Thus, we arrive at an effective model in the large-$n$ limit:
\begin{equation}
T_k = S_k(\beta_*) + U_k
    \label{eq:effective-model}
\end{equation}
where $U_k$'s are approximately distributed as independent 
$N(0,1-\rho^{2k})$.

If $\beta_*$ is known (broken sample problem), by CLT $S_k(\beta_*)$ is approximately $N(0,\rho^{2k})$. 
Thus the effective model \prettyref{eq:effective-model} is approximately Gaussian, for which the optimal test is the likelihood ratio test for testing  $\prod_{k\geq 1} N\Big(\Big[\begin{smallmatrix}
    0\\0
\end{smallmatrix}\Big], \Big[\begin{smallmatrix}
    \rho^{2k} & \rho^{2k}\\\rho^{2k} & 1
\end{smallmatrix}\Big]\Big)$
against 
$\prod_{k\geq 1} N\Big(\Big[\begin{smallmatrix}
    0\\0
\end{smallmatrix}\Big], \Big[\begin{smallmatrix}
    \rho^{2k} & 0 \\ 0 & 1
\end{smallmatrix}\Big]\Big)$,
known as the \textit{quadratic discriminant analysis}. 
Alternatively, one may consider the bilinear test statistic
\begin{equation}
\sum_{k \geq 1} T_k S_k(\beta_*)
\label{eq:LRT}
\end{equation}
which is a special case of the  eigenfunction-based test in the Gaussian broken sample problem \cite[Eq.~(15)]{jiao2025DetectionBroken}. 

For the present shuffled regression problem, $\beta_*$ is \emph{unknown}, and the effective model \prettyref{eq:effective-model} 
is no longer normal. 
This \emph{non-normality} is manifested in two aspects:
a) each marginal of $T_k$ is a Gaussian mixture;
b) different $T_k$'s are dependent through a common latent $\beta_*$.
It is instructive to consider the first two Hermite coefficients $T_1$ and $T_2$ as an example. 
As defined in \prettyref{eq:T_k_decomp}, their conditional means can be expressed as a linear and quadratic form of the same latent $\beta_*$:
\[
S_1(\beta_*) = \rho \beta_*^{\top} \bar X, \quad
S_2(\beta_*) = \frac{\rho^2}{\sqrt{2}} \beta_*^{\top} \bar \Sigma\beta_*,
\]
where $\bar X \equiv \frac{1}{\sqrt{n}} \sum_{i=1}^n X_i$ and 
$\bar \Sigma \equiv \frac{1}{\sqrt{n}} \sum_{i=1}^n (X_iX_i^\top-I)$ are scaled versions of the
sample mean and sample covariance matrix of $X$.
Since $\beta_* \sim \Unif(\mathbb S^{d-1})$, conditioned on $X$, 
$S_1$ is equal in law to 
$\rho \|\bar X\| \cdot \beta_{*,1}$, where $(\beta_{*,1}+1)/2 \sim \text{Beta}(\frac{d-1}{2},\frac{d-1}{2})$. Convolving this with $N(0,1-\rho^2)$ yields the distribution of $T_1$, which is a Gaussian mixture. 
Furthermore, \emph{coupled by the common $\beta_*$}, $S_1$ and $S_2$ are dependent and so are $T_1$ and $T_2$ even in the large-$n$ limit.

Using only $T_1$, testing its Gaussian mixture distribution against the marginal $N(0,1)$ has already achieved non-trivial detection. To attain strong detection, one needs to pool information from all Hermite coefficients and exploit their dependency.
Since it is difficult to describe the joint law of $T_k$'s (which are polynomials of the uniform $\beta_*$), we modify the likelihood ratio test \prettyref{eq:LRT} for known $\beta_*$  by either \textit{averaging} or \textit{optimizing} the latent $\beta_*$.
This motivates two different detection strategies, which 
are analyzed in Sections \ref{sec:cross-moment} and \ref{sec:sliced-moment}.
In \prettyref{sec:ASWT}, we turn to another test function class (step functions), which yields a test that remedies the underperformance of polynomials in high dimensions.

\subsection{Cross-moment test}
\label{sec:cross-moment}
Motivated by the decomposition \prettyref{eq:T_k_decomp}, we begin by constructing our first test statistic by averaging the latent variable $\beta_*$ in \eqref{eq:LRT}. Specifically, we define the statistic
\begin{align}
\sum_{k \geq 1} T_k \underbrace{\mathbb{E}_{\beta_*}\left[S_k(\beta_*)\right]}_{\equiv S_k}, \label{eq:test_vallina}
\end{align}
which follows the \emph{matched-filter} principle in signal processing \cite{oppenheim1997signals}: we correlate the observed Hermite coefficients $T_k$ with their conditional means $S_k=\mathbb E_1[T_k\mid X]$ under the alternative.
Because under $\calH_0$, $T_k$'s are centered with identity covariance conditioned on $X$, \prettyref{eq:test_vallina} is the best linear statistic in the sense of maximizing the signal-to-noise ratio.

Specifically, under the null hypothesis, for all $k \ge 1$, $T_k$ and $S_k$ are independent, zero-mean, so that $\Expect_0[T_kS_k]=0$. Under the alternative, 
each $S_k$ ``matches'' the signal part of $T_k$ in \prettyref{eq:T_k_decomp}, so that their product has a positive mean. Indeed,
\begin{equation}
 \Expect_1 \left[ T_k S_k\right]
= \mathbb{E}_{X} \left[  S_k \expect{T_k \mid X} \right]
= \mathbb{E}_X \left[  S_k^2 \right]
=\mathbb{E}_{\beta \indep \beta_*} \mathbb{E}_X [S_k(\beta)S_k(\beta_*)]
=\rho^{2k} \mathbb{E}
\left[\beta_{*,1}^k\right],
\label{eq:mf-intuition}
\end{equation}
where the second equality follows from
$\expect{T_k\mid X,\beta_*}=S_k(\beta_*)$ in view of \prettyref{eq:T_k_decomp}; the third equality applies the definition of $S_k$ and writes the second moment as averaging over two independent replicas $\beta$ and $\beta_*$; and the last equality holds because $\mathbb{E}_X [S_k(\beta)S_k(\beta_*)]=\rho^{2k} \iprod{\beta}{\beta_*}^k$, which has the same distribution as 
$\rho^{2k}\beta_{*,1}^k$.

Overall, the mean separation of the test statistic \prettyref{eq:test_vallina} is
\begin{align}
\sum_{k\ge 1} \Expect_1 \left[ T_k S_k\right]
= \sum_{k\ge 1} \rho^{2k} \mathbb{E}
\left[\beta_{*,1}^k\right]
= \mathbb{E}\left[  \frac{\rho^2 \beta_{*,1}}{1-\rho^2 \beta_{*,1}}\right] . \label{eq-first-moment-d=123}
\end{align}
Note that $\beta_{*,1}$ is uniformly distributed over $\{\pm 1\}$ for $d=1$ and 
has density proportional to $(1-x^2)^{(d-3)/2}$ for $d\ge 2$. 
This implies that the RHS of \eqref{eq-first-moment-d=123} diverges to infinity when $\rho \to 1$ for $d=1,2,3$  and remains bounded otherwise. Thus, we expect the statistic in~\prettyref{eq:test_vallina} to achieve strong detection when $\rho^2 \to 1$ for $d=1,2,3$, but fail otherwise. 

To achieve strong detection for $d\ge 4,$ we further leverage the dependency among $T_k$'s via their higher cross-moments and compute its conditional mean as the matched filter. Specifically, given a set $A$ of positive integers, let
\begin{equation}\label{eq-def-TASA}
T_A = \prod_{k\in A} T_k \quad \text{ and } \quad S_A = \mathbb E_1[T_A \mid X]\,.
\end{equation}
Extending~\prettyref{eq:test_vallina}, define the \emph{cross-moment test} statistic as
\begin{equation}\label{eq-testing-statistic(I)}
    \mathsf{CMT} \triangleq \sum_{A\in\binom{[K]}{D}}T_A S_A\,, 
\end{equation}
where the sum ranges over all $D$-element subsets of $[K]$.
Similarly as before, under the null hypothesis, $T_A$ and $S_A$ are independent with means nearly zero, so that 
$\mathbb{E}_0[\mathsf{CMT}]=o(1)$. Under the alternative,
$$
\Expect_1 \left[ T_A S_A \right]
= \mathbb{E}_{X} \left[  S_A \expect{T_A \mid X} \right]
= \mathbb{E}_X \left[  S_A S_A \right]
\approx \mathbb{E}_{\beta \indep \beta_*} \mathbb{E}_X \left[\prod_{k \in A} S_k(\beta)S_k(\beta_*) \right]
\approx  \mathbb{E}
\left[ \prod_{k \in A} \left( \rho^2 \beta_{*,1}\right)^{k} \right],
$$
where the approximations can be justified in the large-$n$ limit in view of~\prettyref{eq:effective-model} and the orthogonality of the Hermite polynomials. It follows that 
$$
 \Expect_1 \left[ \mathsf{CMT} \right]=
 \sum_{A\in\binom{[K]}{D}}
 \Expect_1 \left[ T_A S_A \right]
 \approx \sum_{A\in\binom{[K]}{D}}
 \Expect_1 \left[ \prod_{k \in A} \left( \rho^2 \beta_{*,1}\right)^{k} \right]
 \approx \frac{1}{D!} \mathbb{E}\left[  
 \left(\sum_{k=1}^K \left( \rho^2 \beta_{*,1}\right)^{k}   \right)^{D} \right]
$$
where the factor of $D!$ is due to the elements of $A$ being distinct. Recall that $\rho^2=\frac{1}{1+\sigma^2}$. By choosing $K\approx 1/\sigma^2$, and focusing on the probability density of $\beta_{*,1}$ around $1-\sigma^2$, we can get that 
$$
 \Expect_1 \left[ \mathsf{CMT} \right] \approx \frac{1}{D!} \mathbb{E}\left[  
 \left(\sum_{k=1}^K \left( \rho^2 \beta_{*,1}\right)^{k}   \right)^{D} \right]
 \approx \frac{1}{D!} \left(\frac{1}{\sigma^2}\right)^D \left( 
 \sigma^2\right)^{(d-3)/2} \sigma^2,
$$
which diverges to infinity, when $1/\sigma^2 \gtrsim D\gtrsim d$ and $D \to \infty.$ It turns out that to bound the variance of $\mathsf{CMT}$, we further need $K\gg D^2$, so that, by birthday paradox,  two subsets chosen from $\binom{[K]}{D}$ are disjoint with high probability. Together these yield the  success condition $1/\sigma \gg d.$

Building on the above intuitions, we establish the following theoretical guarantee for the cross-moment test. 
\begin{theorem}\label{thm:CMT}
    When 
    $d\leq (\log n)^{0.3}$, we choose
   \begin{equation}\label{eq-choice-par-CMT}
    D=30d,\quad K= \left\lfloor D^2\Big(\frac{1}{\sigma d}\wedge (\log n)^{0.01}\Big) \right\rfloor\,.
    \end{equation}
    If $\sigma^2=o(1/d^2)$, then the cross-moment test~\eqref{eq-testing-statistic(I)} achieves strong detection. Moreover, it can be computed in $n^{1+o(1)}$ time (see Appendix~\ref{sec:computation-CMT}).
\end{theorem}

\subsection{Sliced-moment test}
\label{sec:sliced-moment}

In the previous construction of the cross-moment test, we averaged over the unknown regression direction \(\beta_*\). Next, we seek a new test statistic by \emph{optimizing} over the \(\beta_*\) direction. Specifically, we truncate the bilinear statistic in \eqref{eq:LRT} at degree \(M\) and maximize it over candidate directions \(\beta\), obtaining the \emph{sliced-moment test}:
 \begin{align}
     \mathsf{SMT} &\triangleq \max_{\beta\in N_\delta}  
     \sum_{k=1}^M T_k  S_k(\beta)
     \label{eq-def-SMT}\,,
\end{align} 
 where $M$ is an integer and $N_\delta$ is a $\delta$-net of the unit sphere. 
In other words, we slice each $X_i$ along different candidate directions $\beta$ to capture dependencies between $\langle X_i, \beta \rangle$'s and $Y_i$'s under the alternative hypothesis $\calH_1$.
Choosing $\delta=1/\sqrt{M}$, there exists $\hat{\beta}$ such that $\iprod{\beta_*}{\hat{\beta}} \ge 1-\frac{1}{M}$. Using this approximate $\hat{\beta}$ is almost as good as $\beta_*$, because, similar to \prettyref{eq:mf-intuition} we have for all $k\leq M$,
 $$
  \Expect_{1}\left[T_k  S_k(\hat{\beta}) \right]
  = \Expect_{\beta_*}
  \Expect_{X} \left[ S_k(\beta_*) S_k(\hat{\beta}) \right]
  =  \Expect_{\beta_*} \left[ \rho^{2k} \iprod{\beta_*}{\hat{\beta}}^k\right]
  \ge  \rho^{2k} \left( 1- \frac{1}{M} \right)^k \gtrsim \rho^{2k}.
 $$
Therefore, recalling $\rho^2 = \frac{1}{1+\sigma^2}$, for $M\lesssim 1/\sigma^2 $ we obtain 
$\Expect_1\left[\mathsf{SMT}\right] \gtrsim M.$ Under the null, for any $\beta \in N_\delta,$ by CLT $\{S_k(\beta)\}_{k=1}^M$ and $\{T_k\}_{k=1}^M$ converge to independent standard normals. Thus,  for any constant $c>0,$
$$
\mathbb{P}_0 \left\{ \mathsf{SMT}  \ge c M \right \}
\le \sum_{\beta \in N_\delta} \mathbb{P}_0 \left\{ 
 \sum_{k=1}^M T_k  S_k(\beta) \ge c M \right\}
 \le \exp\left( - \Theta(M) + d \log \frac{3}{\delta} \right),
$$
where the last inequality follows from the standard concentration inequality for sub-exponential random variables and $|N_\delta| \leq (3/\delta)^d$ \cite[Corollary 27.4]{polyanskiy2025information}. 
In conclusion, $\mathsf{SMT}$ achieves strong detection when $d \log M \lesssim M \lesssim 1/\sigma^2$.

By formalizing the above derivations, we obtain the following theoretical guarantee for the sliced-moment test. 
\begin{theorem}\label{thm:SMT}
Suppose $d=o(\frac{\log n}{\log \log n})$. Let $\delta=\frac{1}{\sqrt{M}}$ with 
\begin{equation}\label{eq-def-M-d-large}
    M\triangleq\left\lfloor\frac{1}{\sigma^2}\wedge 
    \sqrt{\frac{d \log n}{\log\log n}}
    {\log(d+1)}\right\rfloor.
\end{equation}
Then the sliced-moment test achieves strong detection if $\sigma^2\to 0$ and $\sigma^2\leq\frac{c}{d\log(d+1)}$ for some sufficiently small constant $c>0$. In addition, the sliced-moment test can be computed in  $n^{1+o(1)}$ time.
\end{theorem}

\subsection{Approximate sliced-Wasserstein test}
\label{sec:ASWT}

Both the cross-moment test and the sliced-moment test rely on Hermite polynomials of degree $k = \poly(d)$. However, these tests face substantial challenges in the very high-dimensional regime. The main obstacle is that the $p$-th moments of normalized Hermite polynomials grow exponentially in $k$ for $p \geq 3$:
\begin{equation}\label{eq-higher-moments-hermite-explode}
\mathbb E_{z\sim N(0,1)}[|h_k(z)|^p] = \exp(\Theta(k)), \quad p\geq 3\,.
\end{equation}
This exponential growth condition \eqref{eq-higher-moments-hermite-explode} affects our analysis in two significant ways:
\begin{itemize}
\item \textbf{Second moment of the cross-moment test.} Expanding $\mathbb E_1[T_A T_B S_A S_B]$ reveals dominant contributions from index pairings, while additional ``catastrophic'' terms arise when certain indices appear with multiplicity $\geq 3$. These terms invoke \eqref{eq-higher-moments-hermite-explode} yet decay only as $n^{-1/2}$. Controlling them requires $k = O(\log n)$. 
    \item \textbf{Gaussian approximation in the sliced-moment test.} 
    To establish CLT for each $T_k$, 
    the error term
$(1/\sqrt{n})\mathbb E_{z\sim N(0,1)}[|h_k(z)|^p]$
must remain small, as required by Berry–Esseen bounds or Stein’s method. This condition holds only when $k = O(\log n)$. 
\end{itemize}
Since both tests require $k=\poly(d)$ to extract enough signal, we are restricted to $d \leq \polylog(n)$.
To address this restriction, we  resort to comparing the empirical distributions via one-dimensional projections. 
Specifically, let 
\begin{align*}
     \mu_{X\beta}\triangleq\frac{1}{n}\sum_{i=1}^n \delta_{\langle X_i,\beta\rangle}\,, \quad  \mu_Y \triangleq \frac{1}{n}\sum_{i=1}^n \delta_{\frac{Y_i}{\sqrt{1+\sigma^2}}}\,,\notag
\end{align*}
and let $F_{n,\beta}$ and $G_{n}$ be the cumulative distribution function (CDF) of $\mu_{X\beta}$ and $\mu_Y$. Define the \emph{sliced-Wasserstein test} as 
\begin{equation}\label{eq-def-SWT}
    \mathsf{SWT} \triangleq \min_{\beta\in N_\delta}\mathsf{d}_p(F_{n,\beta},G_n),
\end{equation}
where  $N_\delta$ is a $\delta$-net of the unit sphere, and 
\begin{equation*}
\mathsf{d}_p(F_{n,\beta},G_n)=\|F_{n,\beta}-G_n\|_p\triangleq\left( \int_{\mathbb R} |F_{n,\beta}(t)-G_n(t)|^p\mathrm{d}t \right)^{\frac{1}{p}}\,.
\end{equation*}
 Under $\mathcal{H}_1$, by choosing $\delta=\Theta(\sigma^2)$, there exists $\widehat\beta\in N_{\delta}$ such that $\|\widehat{\beta}-\beta_*\| = O(\sigma^2)$. Therefore, the correlation between $\langle X_{\pi^*(i)},\widehat\beta\rangle$ and $Y_i/\sqrt{1+\sigma^2}$ is $1-O(\sigma^2)$. Following the discussion in \cite[Equation (9) to (14)]{ding2021efficient}, we can obtain
\begin{equation}
    \mathsf d_p(F_{n,\widehat{\beta}},G_n) = O\left( \sqrt{\frac{\sigma}{n}} \right) \text{ with probability } 1-o(1)\,.\label{eq-H1-SWT}
\end{equation}
Under $\mathcal{H}_0$, for any $\beta\in N_\delta$,
$
    \mathbb E_0[\mathsf{d}_p(F_{n,\beta},G_n)]=\Theta\left(\frac{1}{\sqrt{n}}\right)\,.
$
Moreover, for each $\beta$, we have the small ball probability estimate:
\begin{equation}\label{eq:small-ball-est-SWT}
\mathbb P_0\left[ \mathsf{d}_p(F_{n,\beta},G_n)\leq \sqrt{\frac{\sigma}{n}} \right] = \exp\left( -\Theta\left(\frac{1}{\sigma}\right) \right).
\end{equation}
By applying a union bound on $N_\delta$, we have
\begin{equation*}
    \mathbb P_0\left[ \min_{\beta\in N_\delta}\mathsf{d}_p(F_{n,\beta},G_n)\leq \sqrt{\frac{\sigma}{n}} \right] \leq \exp\left(- \Theta\left(\frac{1}{\sigma}\right) + d \log \left(\frac{1}{\sigma^2}\right)\right)\,,
\end{equation*}
which is $o(1)$, provided that  $\sigma^2 = o(1/(d\log d)^2)$. Combining this with \eqref{eq-H1-SWT}, we see that the sliced-Wasserstein test achieves strong detection. 
However, obtaining the required small ball probability in \eqref{eq:small-ball-est-SWT} is delicate. To bypass this analytical difficulty, we introduce a pseudo-metric approximating $\mathsf{d}_1(F,G)$ that is easier to control. Concretely, define $I_1,\dots,I_L$ to be the uniform partition of the interval $[-\frac{1}{2},\frac{1}{2}]$ (i.e., $I_k=[-\frac{1}{2}+\frac{k-1}{L},-\frac{1}{2}+\frac{k}{L}]$) and let 
\begin{equation}
    \mathsf{d}(\mu_{X\beta},\mu_{Y}) \triangleq \sum_{\ell \in [L]} |\mu_{X\beta}(I_\ell)-\mu_Y(I_\ell)|\,.
    \notag
\end{equation}
 Then the approximate sliced-Wasserstein test is defined as :
\begin{equation}
    \label{eq-testing-statistic-high-dimension}
    \mathsf{ASWT}=\min_{\beta\in N_\delta} \mathsf{d}(\mu_{X\beta},\mu_{Y})\,.
\end{equation}
The following theorem shows that the approximate sliced-Wasserstein test achieves strong detection when $\sigma^2 = o(1/(d\log d)^2)$ 
and $d=o(n/\log(n)).$
\begin{theorem}\label{thm:ASWT}
    Suppose $d=o\big(\frac{n}{\log n}\big)$ and $\sigma^2=o\big(\frac{1}{(d\log(d+1))^2}\big)$. Let $L \triangleq \lfloor \omega_n^{1/3}d\log(d+1) \rfloor$ with $\omega_n\triangleq\frac{n}{d\log n}\wedge\frac{1}{\sigma^2(d\log(d+1))^2}\wedge\log n$, and $\delta=1/L^2$.
    Then the approximate sliced-Wasserstein test achieves strong detection.
  The computational complexity is bounded by $d(3L^2)^dn^{1+o(1)}$.
\end{theorem}

\subsection{Information-theoretic lower bound}

In the following, we present the impossibility result. Let $L(X,Y)=
\frac{\diff \mathbb P_1[X,Y]}{\diff \mathbb P_0[X,Y]}$ denote the likelihood ratio. It is well-known that 
\begin{align*}
\mathbb E_0 \left[ L^2(X,Y) \right]= 1+o(1)  & \Longrightarrow  \TV(\mathbb P_0, \mathbb P_1) =o(1), \\
 \mathbb E_0 \left[ L^2(X,Y) \right]= O(1) & \Longrightarrow   \TV(\mathbb P_0, \mathbb P_1) =1-\Omega(1).
\end{align*}
Therefore, it suffices to bound the second moment of the likelihood ratio.
\begin{theorem}\label{thm:lower-bound}
    We have the following: 
    \begin{enumerate}[label=(\alph*)]
        \item $d\to \infty$ and $d\leq \frac{n}{(\log n)^3}$:  If $\sigma^2\geq \frac{(1+\varepsilon)\pi^2}{3d\log d}$ for any constant $\varepsilon>0$, then
    \begin{equation}
        \mathbb E_0 \left[ L^2(X,Y) \right]= 1+o(1)\,.\notag
    \end{equation}
        \item $d=O(1)$: If $\sigma^2>\varepsilon$ for any constant $\varepsilon>0$, 
        then
        \begin{equation}
        \mathbb E_0 \left[ L^2(X,Y) \right] = O(1)\,.\notag
    \end{equation}
     Moreover, if $\sigma^2\to\infty$, then $\mathbb E_0[L^2(X,Y)]=1+o(1)$.
    \end{enumerate}
\end{theorem}

\begin{remark}
\prettyref{thm:lower-bound}(a) implies that in the growing-dimension regime, \emph{weak detection} (i.e., the sum of Type-I and Type-II errors is bounded away from $1$) is information-theoretically impossible provided
$\sigma^2 \gtrsim \frac{1}{d \log d}$.  
In contrast, if $\sigma^2 \lesssim \frac{1}{d \log d}$, \prettyref{thm:SMT} shows that strong detection is attainable by the sliced-moment test when $d = o(\log n / \log\log n)$. Therefore, in this regime, the detection problem exhibits a phase transition--from the impossibility of weak detection to the possibility of strong detection-- at the boundary $
\sigma^2 \asymp \frac{1}{d \log d}.$
In higher dimensions $d=\omega(\log n/\log \log n)$, the approximate sliced-Wasserstein test achieves strong detection when $\sigma^2 \lesssim \frac{1}{(d \log d)^2}$. Whether detection is possible in the intermediate region
$
\frac{1}{(d \log d)^2} \lesssim \sigma^2 \lesssim \frac{1}{d \log d}
$
remains an open question.

\prettyref{thm:lower-bound}(b) implies that for fixed dimension, \emph{strong detection} is impossible when $\sigma^2 = \Omega(1)$, and even \emph{weak detection} becomes impossible once $\sigma^2 = \omega(1)$. On the other hand, all three of our proposed tests achieve strong detection in polynomial time when $\sigma^2 \to 0$, implying that the boundary for strong detection lies at $\sigma^2 \asymp 1$. 
Moreover, weak detection is achievable in polynomial time whenever $\sigma^2=O(1)$. For instance, one can show that the distributions of $T_2S_2$
under the two hypotheses have means separated by a positive constant and bounded variances, whenever $\sigma^2=O(1)$ and $d=O(1)$. 
Consequently, their total variation distance is bounded away from zero,
implying that the optimal likelihood ratio test based on $T_2S_2$ achieves weak detection. 
\end{remark}

\section{Extensions}
\label{sec:ext}
In the following two subsections, we provide two important extensions of our results: (i) The all-or-nothing phase transition in the recovery of $\beta_*$; (ii) The detection threshold when $\beta_*$ is sparse.

\subsection{Recovery of regression coefficients}
\label{subsec-recovery-of-beta^*}

The preceding subsections focus on the detection problem. Interestingly, the detection results have direct implications for recovery of the regression coefficient \texorpdfstring{$\beta_*$}{beta}. Specifically, when $\sigma^2 \lesssim (d \log d)^{-1}$, Theorem~\ref{thm:SMT} establishes that the sliced-moment test achieves strong detection. Moreover, the proof shows that the maximizer in \eqref{eq-def-SMT} concentrates around $\beta_*$ under $\mathcal{H}_1$. 
Hence, \prettyref{thm:beta-recovery} below shows that this maximizer constitutes a consistent estimator of $\beta_*$. In contrast, when $\sigma^2 \gtrsim (d \log d)^{-1}$, Theorem~\ref{thm:lower-bound} shows that $X$ and $Y$ become nearly independent. Building upon this, \prettyref{thm:impossibility-beta-recovery} proves that nontrivial estimation of $\beta_*$ is information-theoretically impossible in this regime.

The next theorem shows that the maximizer of the sliced-moment test achieves vanishing recovery error when $\sigma^2 \lesssim (d \log d)^{-1}$.

\begin{theorem}\label{thm:beta-recovery}
Let $\tilde{\beta}$ be a maximizer to~\eqref{eq-def-SMT}. 
    Under the same conditions as~\prettyref{thm:SMT}, 
    we have $\|\widetilde \beta-\beta_*\|_2 = o\left(1\right) $ with high probability
    and $\expect{\|\widetilde \beta-\beta_*\|_2^2}=o(1)$.
\end{theorem}
The next theorem shows that non-trivial estimation for $\beta_*$ is fundamentally impossible when $\sigma^2 \gtrsim (d\log d)^{-1}$. 
\begin{theorem}
    \label{thm:impossibility-beta-recovery}
    If $d\to\infty$ and $\sigma^2\geq \frac{(1+\varepsilon)\pi^2}{3d\log d}$ for any constant $\varepsilon>0$, then the minimum mean-squared error satisfies
    \begin{equation}
        \mathsf{MMSE}(\beta_*) \triangleq  \min_{\hat{\beta}}\mathbb E_1\big[\| \beta_* -  \hat{\beta}(X,Y)\|_2^2\big] = 
1-o(1)\,.\notag
    \end{equation}
\end{theorem}

\prettyref{thm:beta-recovery} and \prettyref{thm:impossibility-beta-recovery} together imply that in the regime 
$\omega(1) \le d \le o(\log n / \log\log n)$, the recovery of $\beta_*$ exhibits a sharp ``all-or-nothing'' transition from vanishing estimation error to trivial estimation error at the boundary $
\sigma^2 \asymp (d \log d)^{-1}$.

In stark contrast, the recovery of the hidden pairing $\pi_*$ requires substantially smaller noise levels. In particular, perfect recovery of $\pi_*$ requires 
$\sigma^2 = n^{-4+o(1)}$, and \emph{almost} perfect recovery requires 
$\sigma^2 = n^{-2+o(1)}$. Moreover, as shown in \cite[Proposition~1]{wang2022random}, when $\sigma^2 = \omega(n^{-\epsilon})$ for any $\epsilon > 0$, even partial recovery of $\pi_*$ becomes impossible. 
Thus, there exists a wide regime in which recovering $\pi_*$ is information-theoretically impossible while $\beta_*$ can still be estimated  with vanishing error. This has important implications for practical applications such as flow cytometry, where the primary objective is to estimate the regression coefficient $\beta_*$ without access to the true ordering of the labels~\cite{abid2017linear}.

\subsection{Sparse regime}
\label{sec:sparse-beta}

In this section, we extend the detection problem to the case where the regression coefficient $\beta_*$ is sparse. 
Let $s = \|\beta_*\|_0$ denote the sparsity level, i.e., the number of nonzero coordinates.
We consider a simple setting in which $\beta_*\in \{0,\frac{1}{\sqrt{s}}\}^d$. This binary setting has been studied before without shuffling in~\cite{reeves2019all,gamarnik2022sparse}.

\begin{definition}[Shuffled sparse linear regression]
\label{def:slr_sparse}
 Consider the hypothesis testing problem given in~\prettyref{def:slr}, where $\beta_*$ is uniformly distributed over the set 
 $
 \mathbb A_{\mathsf{sparse}}
 =\frac{1}{\sqrt{s}} \cdot \{v \in \{0,1\}^d:\|v\|_0=s \}.
 $
\end{definition}

Inspired by \eqref{eq-def-SMT}, we propose the following \emph{sparse sliced-moment test} ($\mathsf{SSMT}$):
\begin{align}
     \mathsf{SSMT} &\triangleq\max_{\beta\in \mathbb A_{\mathsf{sparse}}}  \; \sum_{k=1}^{M'} S_k(\beta)T_k \, ,
     \label{eq-testing-statistic-sparse}
\end{align}
where $M'$ is chosen as
\begin{equation}\label{eq-def-M'}
    M'=\Big\lfloor \frac{1}{\sigma^2}\wedge \left(C_0 \log \binom{d}{s}\right) \Big\rfloor\,,
\end{equation}
where $C_0>0$ is a sufficiently large constant.
The following theorem establishes sufficient and necessary conditions for  achieving strong detection. 
\begin{theorem}
    \label{thm:sparse-beta}
        Suppose that $d\to\infty$, $s=o(d)$ and $s \log d \leq c_0 \log n $ for a sufficiently small constant $c_0>0$.  
    \begin{itemize}
        \item If $\sigma^2\leq \frac{c}{\log\binom{d}{s}}$ for a sufficiently small constant $c>0$,
        then  the test statistic \eqref{eq-testing-statistic-sparse} achieves strong detection;
        \item If $\sigma^2\geq \frac{(1+\epsilon)\pi^2}{6\log\binom{d}{s}}$ for any constant $\epsilon>0$, then $\mathbb E_0[L^2_{\mathsf{sparse}}]=1+o(1)$ and hence weak detection is impossible. 
    \end{itemize}
\end{theorem}

\prettyref{thm:sparse-beta} shows that the detection threshold can be improved from 
$\sigma^2 \asymp (d \log d )^{-1}$ in the non-sparse case to $\sigma^2 \asymp (s \log (d/s))^{-1}$ in the sparse case. 

\section{Numerical experiments}
\label{sec:numerical-results}

In this section, we investigate the finite-sample performance of the proposed tests on synthetic data. All experiments use dimension $d=3$ and sample size $n=1000$. For each parameter configuration, we generate $R=200$ independent datasets under $\mathcal H_0$ and another $R=200$ independent datasets under $\mathcal H_1$, resampling the complete model in each replication. Applying a test statistic to these datasets produces null scores $S_1^{(0)},\ldots,S_R^{(0)}$ and alternative scores $S_1^{(1)},\ldots,S_R^{(1)}$. The scores are chosen so that larger values provide stronger evidence against $\mathcal H_0$; in particular, the score used for the sliced-Wasserstein test is the negative of its minimized empirical Wasserstein distance. The empirical false-positive and true-positive rates at threshold $t$ are
\[
\widehat{\operatorname{FPR}}(t)
=
\frac{1}{R}\sum_{r=1}^R
\mathbf 1\{S_r^{(0)}\geq t\},
\qquad
\widehat{\operatorname{TPR}}(t)
=
\frac{1}{R}\sum_{r=1}^R
\mathbf 1\{S_r^{(1)}\geq t\}.
\]
Varying $t$ over the observed scores yields the empirical receiver operating characteristic (ROC) curve.

\begin{figure}[H]
    \centering
    \begin{subfigure} 
        \centering
        \includegraphics[width=0.45\textwidth]{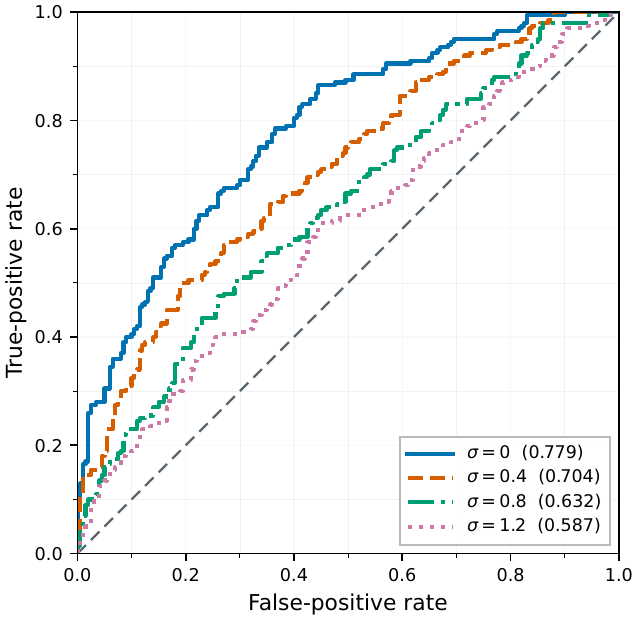}
        \label{fig:left-noise-CMT}
    \end{subfigure}
    \hspace{0.02\textwidth} 
    \begin{subfigure}
        \centering
        \includegraphics[width=0.45\textwidth]{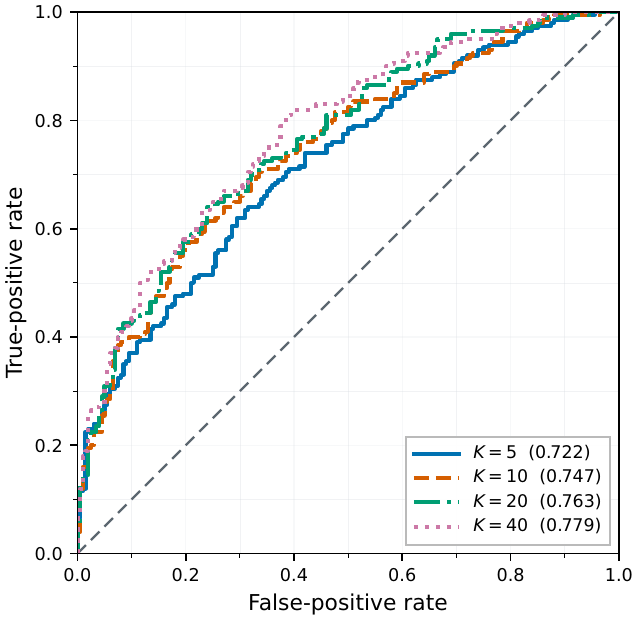}
        \label{fig:right-different-tests-CMT}
    \end{subfigure}
    \caption{Empirical ROC curves for the cross-moment test with $d=3$, $n=1000$, and $D=1$. Left: $K=50$ and $\sigma\in\{0,0.4,0.8,1.2\}$. Right: $\sigma=0$ and $K\in\{5,10,20,40\}$. The numbers in parentheses are the empirical AUC values defined in~\eqref{eq:empirical-AUC}.}
    \label{fig:two_images_CMT}
\end{figure}

\paragraph{Cross-moment test}
We consider the test statistic \eqref{eq-testing-statistic(I)} with $D=1$, namely $
    f=\sum_{k=1}^K T_kS_k$. 
In Figure~\ref{fig:two_images_CMT}, we sample $\beta_*$ uniformly from $\mathbb S^2$. The left panel fixes the degree $K=50$ and varies $\sigma\in\{0,0.4,0.8,1.2\}$, whereas the right panel fixes $\sigma=0$ and varies $K\in\{5,10,20,40\}$. The test becomes less powerful as the noise increases, while including more moments improves its performance in the noiseless setting.

\begin{figure}[H]
    \centering
    \begin{subfigure} 
        \centering
        \includegraphics[width=0.45\textwidth]{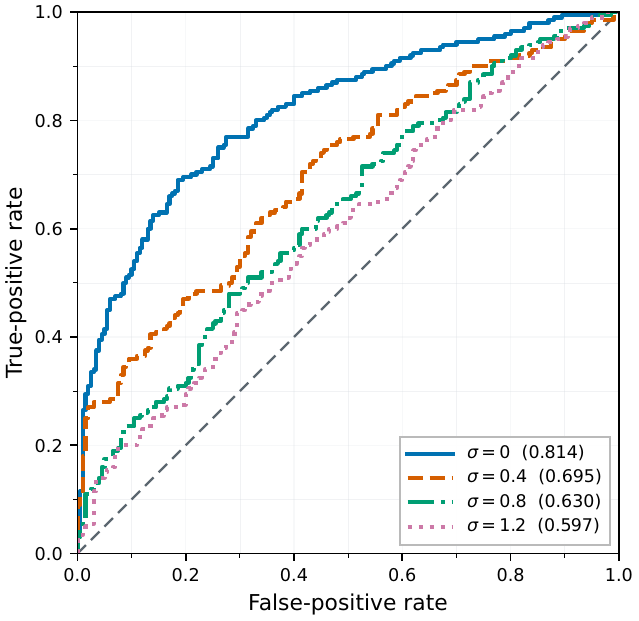}
        \label{fig:left-noise-SMT}
    \end{subfigure}
    \hspace{0.02\textwidth} 
    \begin{subfigure}
        \centering
       \includegraphics[width=0.45\textwidth]{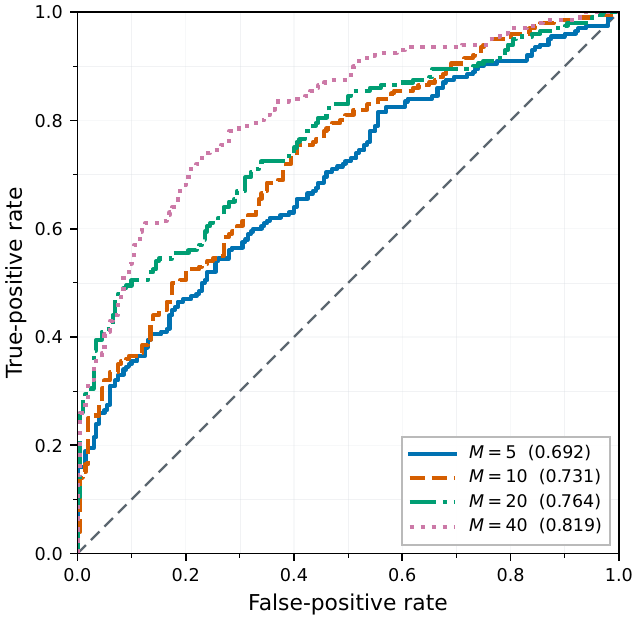}
        \label{fig:right-different-tests-SMT}
    \end{subfigure}
    \caption{Empirical ROC curves for the sliced-moment test with $d=3$, $n=1000$, and $|N_\delta|=100$. Left: $M=50$ and $\sigma\in\{0,0.4,0.8,1.2\}$. Right: $\sigma=0$ and $M\in\{5,10,20,40\}$.}
    \label{fig:two_images_SMT}
\end{figure}

\paragraph{Sliced-moment test} 
The simulation results for the test in \eqref{eq-def-SMT} are summarized in Figure~\ref{fig:two_images_SMT}. We sample $\beta_*$ uniformly from $\mathbb S^2$ and maximize over a fixed grid of size $100$.\footnote{For the sliced tests below, we construct a fixed nested collection of directions on $\mathbb S^2$ by greedy farthest-point selection from a deterministic $10^5$-point spherical Fibonacci candidate set~\cite{gonzalez2010measurement,marques2013spherical}.} The left panel fixes $M=50$ and varies $\sigma\in\{0,0.4,0.8,1.2\}$; the right panel fixes $\sigma=0$ and varies $M\in\{5,10,20,40\}$. Similar to the cross-moment test, performance degrades with increasing noise level and improves with the number of moments in the noiseless setting.

\begin{figure}[H]
    \centering
    \begin{subfigure} 
        \centering
        \includegraphics[width=0.45\textwidth]{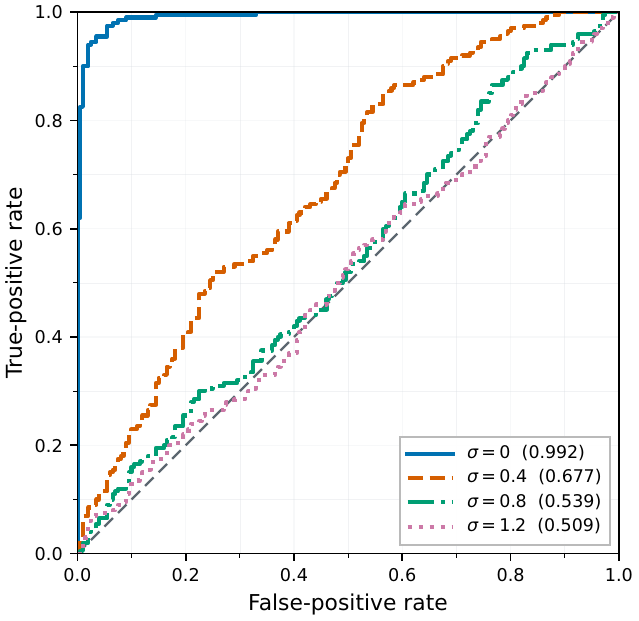}
        \label{fig:left-noise-SWT}
    \end{subfigure}
    \hspace{0.02\textwidth} 
    \begin{subfigure}
        \centering
       \includegraphics[width=0.45\textwidth]{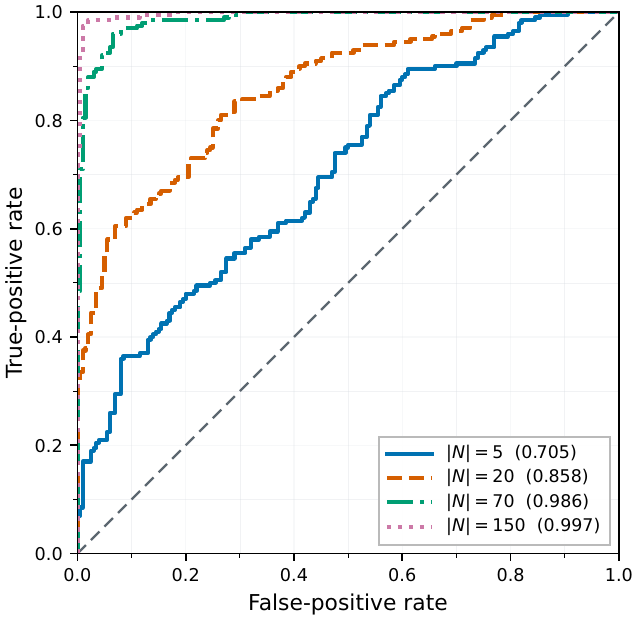}
        \label{fig:right-different-tests-SWT}
    \end{subfigure}
    \caption{Empirical ROC curves for the sliced-Wasserstein test with $d=3$, $n=1000$, and $p=1$. Left: $|N_\delta|=100$ and $\sigma\in\{0,0.4,0.8,1.2\}$. Right: $\sigma=0$ and $|N_\delta|\in\{5,20,70,150\}$.}
    \label{fig:two_images_SWT}
\end{figure}

\paragraph{Sliced-Wasserstein test} 
We evaluate the sliced-Wasserstein test in \eqref{eq-def-SWT} with $p=1$ and sample $\beta_*$ uniformly from $\mathbb S^2$. In the left panel of Figure~\ref{fig:two_images_SWT}, we minimize over the same size-$100$ net of $\mathbb S^2$ used above and vary $\sigma\in\{0,0.4,0.8,1.2\}$. In the right panel, we fix $\sigma=0$ and vary the net size over $\{5,20,70,150\}$. The test is particularly powerful at low noise, and its performance improves substantially as the net is refined.

\begin{figure}[H]
    \centering
    \includegraphics[width=0.72\textwidth]{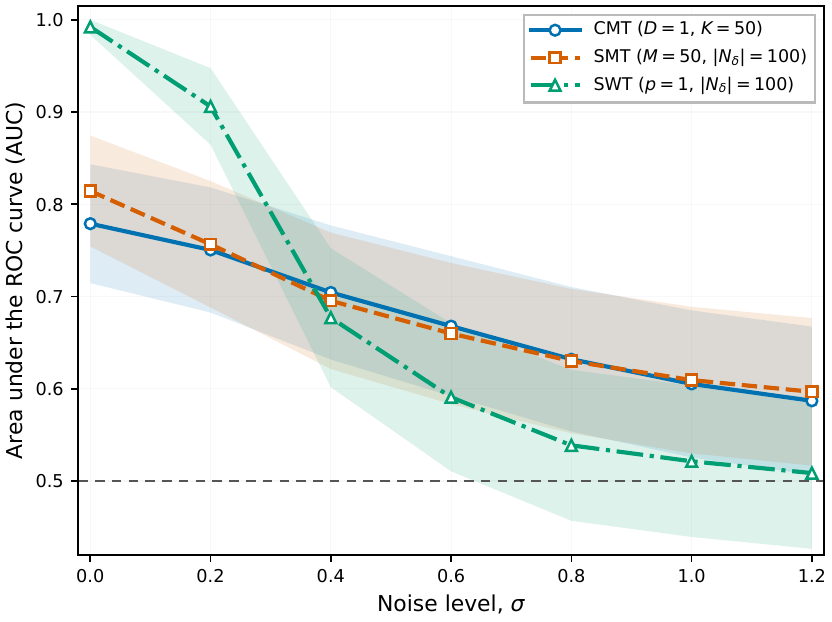}
    \caption{Estimated AUC versus noise level for the three tests with $d=3$, $n=1000$, and $200$ datasets under each hypothesis.}
    \label{fig:auc-sigma-comparison}
\end{figure}

\paragraph{Area under the ROC curve}
To compare the three methods directly, let $S^{(0)}$ and $S^{(1)}$ denote independent population scores under $\mathcal H_0$ and $\mathcal H_1$, respectively. The area under the population ROC curve is
\[
\operatorname{AUC}
=
\mathbb P\bigl(S^{(1)}>S^{(0)}\bigr)
+\frac12\mathbb P\bigl(S^{(1)}=S^{(0)}\bigr).
\]
Thus, $\text{AUC}=1/2$ corresponds to random guessing and $\text{AUC}=1$ means perfect separation.
The empirical estimator of AUC is the two-sample Mann--Whitney statistic
\begin{equation}\label{eq:empirical-AUC}
\widehat{\operatorname{AUC}}
=
\frac{1}{R^2}\sum_{r=1}^R\sum_{s=1}^R
\left[
\mathbf 1\bigl\{S_r^{(1)}>S_s^{(0)}\bigr\}
+\frac12\mathbf 1\bigl\{S_r^{(1)}=S_s^{(0)}\bigr\}
\right].
\end{equation}
The values in parentheses in the legends of Figures~\ref{fig:two_images_CMT}--\ref{fig:two_images_SWT} are the corresponding empirical AUC values. Figure~\ref{fig:auc-sigma-comparison} plots the estimated AUC versus $\sigma$ for all three tests. The shaded regions are $95\%$ simultaneous confidence bands constructed using the nonparametric estimator in \cite{delong1988comparing}.
\prettyref{fig:auc-sigma-comparison} shows that SWT has the largest AUC at low noise, but its advantage largely disappears as $\sigma$ increases, and for moderate to high noise the two moment-based tests achieve higher AUC. 

\section{Proof of Theorem~\ref{thm:CMT}}
\label{sec:integral-over-beta-testing-statistic}

The computational guarantee in Theorem~\ref{thm:CMT} is proved in Appendix~\ref{sec:computation-CMT}. In this section, we prove the statistical guarantee.

\subsection{Statistical analysis of cross-moment test}

We prove Theorem~\ref{thm:CMT} via the second moment method. We assume 
$d \le (\log n)^{0.3}$ and $K, D$ are chosen according to~\eqref{eq-choice-par-CMT} throughout the analysis. 
Propositions~\ref{prop-H_0-moments-CMT} and \ref{prop-H_1-moments-CMT} below bound the first and second moments of $\CMT$~\eqref{eq-testing-statistic(I)} under $\mathcal{H}_0$ and $\mathcal{H}_1$, respectively.

\begin{proposition}
    \label{prop-H_0-moments-CMT}
     We have
    \begin{align}
        & \left| \mathbb E_0\left[ \CMT \right] \right| \leq n^{-D/3+o(1)}\,,\label{eq-H0-first-moment-CMT}\\
        & \mathbb E_0\left[\CMT^2\right]=\sum_{A\in \binom{[K]}{D}}\mathbb E_{\beta_1,\beta_2}\left[ \prod_{k\in A}\left(\frac{\langle\beta_1,\beta_2\rangle }{1+\sigma^2}\right)^k \right]+n^{-1+o(1)}\,.\label{eq-H0-second-moment-CMT}
    \end{align}
\end{proposition}

\begin{proposition}
    \label{prop-H_1-moments-CMT}
   We have
    \begin{align}
       & \mathbb E_1\left[ \CMT \right] = \sum_{A\in\binom{[K]}{D}}\mathbb E_{\beta_1,\beta_2}\left[ \prod_{k\in A}\left(\frac{\langle\beta_1,\beta_2\rangle }{1+\sigma^2}\right)^k \right]+n^{-1+o(1)}\,.\label{eq-H1-first-moment-CMT}\\
      & \var_1[\CMT]
      \le \sum_{1 \le t \le D} \frac{(3K)^t}{ [(D-t)!]^2t!} 
     \left( \mathbb E_{\beta_1,\beta_2} \left[ \left(\sum_{k=1}^K \left(\frac{|\langle\beta_1,\beta_2\rangle|}{1+\sigma^2}\right)^k \right)^{D-t}\right]\right)^2 + n^{-1+o(1)}.
        \label{eq-H1-second-moment-CMT}
    \end{align}
\end{proposition}
We also need the following two auxiliary lemmas to obtain more explicit estimates of the first and second moments. 
Let 
$$ \mu_{d,K,D} \triangleq \mathbb E_{\beta_1,\beta_2} \left[ \left(\sum_{k=1}^K \left(\frac{\langle\beta_1,\beta_2\rangle}{1+\sigma^2}\right)^k \right)^{D}\right] \quad \text{and} \quad
\nu_{d,K,D}(t) \triangleq \mathbb E_{\beta_1,\beta_2} \left[ \left(\sum_{k=1}^K \left(\frac{|\langle\beta_1,\beta_2\rangle|}{1+\sigma^2}\right)^k \right)^{D-t}\right].
$$
\begin{lemma}
    \label{lem:est-sum-first-second-moment-growing-d}
   It holds that 
    \begin{align}
        \label{eq-est-first-moment-growing-d}
        & \sum_{A\in\binom{[K]}{D}}\mathbb E_{\beta_1,\beta_2}\left[ \prod_{k\in A}\left(\frac{\langle\beta_1,\beta_2\rangle }{1+\sigma^2}\right)^k \right] = \frac{1+o(1)}{D!} \mu_{d,K,D}\,.
    \end{align}
\end{lemma}

\begin{lemma}
    \label{lem:est-ratio-Dminust-versus-D}
    We have the following estimate:
    \begin{equation}
\begin{split}
    & \mu_{d,K,D} \geq K^{D-\frac{d-1}{2}}\left(\frac{1-e^{-2}}{3}\right)^D\frac{1}{2\pi}\,.
\end{split}\label{eq:mean_mu_lower_bound}
\end{equation}

Moreover, for $1\leq t \leq D-5d,$
    \begin{align}
        &\label{eq-est-first-moment-growing-d-goal1-case1}
    \nu_{d,K,D}(t)  \leq 4 \left( \frac{1-e^{-10}}{11}K \right)^{-t} \mu_{d,K,D}\,, 
    \end{align}
    and for $D-5d\leq t\leq D$, $\nu_{d,K,D}(t)  \leq K^{D-t} \le K^{5d}.$
\end{lemma}

The proofs of Propositions~\ref{prop-H_0-moments-CMT} and \ref{prop-H_1-moments-CMT}
are deferred to \prettyref{sec:propositions_cross_moment},
while those of Lemmas~\ref{lem:est-sum-first-second-moment-growing-d} and \ref{lem:est-ratio-Dminust-versus-D} are provided in Subsection~\ref{subsec-two-lemmas}. With these ingredients in hand, we are ready to prove Theorem~\ref{thm:CMT}.
\begin{proof}[Proof of Theorem~\ref{thm:CMT}]

     By Proposition~\ref{prop-H_1-moments-CMT}, Lemma~\ref{lem:est-sum-first-second-moment-growing-d} and~\prettyref{eq:mean_mu_lower_bound} in Lemma~\ref{lem:est-ratio-Dminust-versus-D}, we have
\begin{equation}
    \mathbb E_1[\CMT] 
    =\frac{1+o(1)}{D!} \mu_{d,K,D} + n^{-1+o(1)} 
    \geq \frac{K^{D-\frac{d-1}{2}}}{D!}\left(\frac{1-e^{-2}}{3}\right)^D\frac{1+o(1)}{2\pi}
    \to \infty\,, \label{eq:mean_CMT_lower_bound}
\end{equation}
where the divergence holds because $D=30d$ and $K=\lfloor D^2\widehat\omega_n\rfloor$ with $\widehat\omega_n\triangleq\frac{1}{\sigma d}\wedge(\log n)^{0.01}\to\infty$, as $\sigma^2=o(1/d^2)$.
Moreover, by Proposition~\ref{prop-H_0-moments-CMT} and \ref{prop-H_1-moments-CMT}, $\mathbb{E}_0[\CMT^2] = \mathbb E_1[\CMT] + n^{-1+o(1)}. $  Therefore, applying Markov's inequality yields that for any constant $0<c<1$, 
\begin{equation}
    \begin{split}
        \mathbb P_0\left[ \CMT\geq c\mathbb E_1[\CMT] \right] 
        \le \mathbb P_0\left[ \CMT^2 \geq  \left(c \mathbb E_1[\CMT]\right)^2 \right]
        \le \frac{\mathbb{E}_0\left[\CMT^2\right]}{\left(c\mathbb E_1[\CMT]\right)^2}=o(1)\,.
    \end{split}\notag
\end{equation}
Moreover, applying Chebyshev's inequality yields that
\begin{equation}
\begin{split}
 \mathbb P_1\left[ \CMT\leq c\mathbb E_1[\CMT] \right] = \mathbb P_1\left[ \CMT- \mathbb E_1[\CMT] \leq -(1-c)\mathbb E_1[\CMT]\right]
    \leq \frac{\operatorname{Var}_1[\CMT]}{(1-c)^2\left(\mathbb E_1[\CMT]\right)^2}.
\end{split}\notag
\end{equation}
To complete the proof, it remains to show $\var_1[\CMT]=o(\left(\mathbb E_1[\CMT]\right)^2).$
   By Lemma~\ref{lem:est-ratio-Dminust-versus-D}, 
   \begin{align*}
   \var_1[\CMT]
       &\le \sum_{1 \le t \le D-5d} \frac{( O(K))^t}{ [(D-t)!]^2t!} 
       K^{-2t} \mu^2_{d,K,D}
       + 
       O\left( \frac{(3K)^{D+10d}}{(D-5d)!}\right) + n^{-1+o(1)} \\
       & \le O\left( \frac{D^2}{K} \right)\cdot 
       \frac{\mu^2_{d,K,D}}{[D!]^2}
       +O\left( \frac{(3K)^{D+10d}}{(D-5d)!}\right) + n^{-1+o(1)} ,
   \end{align*}
   where the last inequality holds because $D!/(D-t)! \le D^t$  and $K \gg D^2$. It then follows from~\prettyref{eq:mean_CMT_lower_bound} that
   $$
   \frac{\operatorname{Var}_1[\CMT]}{\left(\mathbb E_1[\CMT]\right)^2}
   \le O\left( \frac{D^2}{K} \right)
   + \frac{K^{D+10d}}{(D-5d)!} \frac{[D!]^2[O(1)]^D}{K^{2D-d+1}}  
   =o(1),
   $$
   where the last inequality holds because $D=30d$ and $K\ge D^2\widehat\omega_n/2$,
   $$
   \frac{K^{D+10d}}{(D-5d)!}  \frac{[D!]^2[O(1)]^D}{K^{2D-d+1}}
\le \frac{D^{D+5d} [O(1)]^D}{ K^{D-11d+1}}
\le \frac{D^{7D/6} [O(1)]^D}{ (D^2 \hat{\omega}_n)^{19D/30}}
= o(1).
   $$

\end{proof}

\subsection{Proof of Propositions~\ref{prop-H_0-moments-CMT} and \ref{prop-H_1-moments-CMT}} \label{sec:propositions_cross_moment}

\subsubsection{Proof of Proposition~\ref{prop-H_0-moments-CMT}}
\label{subsec-proof-prop-H_0-moments-CMT}
Under $\mathcal H_0$, we have
\begin{equation}\label{eq-H0-first-moment-d-consant}
    \mathbb E_0[\CMT]
    =\sum_{A\in\binom{[K]}D}\mathbb E_0[T_AS_A]
    =\sum_{A\in\binom{[K]}D}\mathbb E_0[T_A]\mathbb E_0[S_A]
    =\sum_{A\in\binom{[K]}D}\bigl(\mathbb E_0[T_A]\bigr)^2.
\end{equation}
For fixed subset $A$ of degrees, 
denote the index tuple $\mathbf i_A=(i_k)_{k\in A}$ where each index $i_k\in[n]$. Let
\begin{equation*}
    h(\mathbf i_A)\equiv
    \mathbb E_0\left[\prod_{k\in A}h_k\left(\frac{Y_{i_k}}{\sqrt{1+\sigma^2}}\right)\right].
\end{equation*}
Expanding $T_A$, we obtain
\begin{equation}\label{eq-H0-first-moment-expansion-d-constant}
    \mathbb E_0[T_A]
    =\frac1{n^{D/2}}\sum_{\mathbf i_A\in[n]^A}h(\mathbf i_A).
\end{equation}
By independence and orthogonality of Hermite polynomials, we conclude
\begin{equation}\label{eq-H0-cancelation-first-moment-d-constant}
    h(\mathbf i_A)=0
    \quad\text{if an index occurs once or twice in }\mathbf i_A.
\end{equation}
Thus if $h(\mathbf i_A)\neq 0$, every appearing index occurs at least three times, and there are at most $\lfloor D/3\rfloor$ distinct indices. By H\"older's inequality and Lemma~\ref{lem:fourth-moment-hermite-estimate},
\begin{equation}\label{eq-crude-est-D-tuple}
    |h(\mathbf i_A)|
    \le\prod_{k\in A}\mathbb E_0\left[\left|h_k\left(\frac{Y_{i_k}}{\sqrt{1+\sigma^2}}\right)\right|^D\right]^{1/D}
    \le\exp(KD\log(KD)).
\end{equation}
There are at most $\binom n\ell\ell^D$ index tuples using exactly $\ell$ distinct indices. Hence
\begin{equation}
    |\mathbb E_0[T_A]|
    \le\frac1{n^{D/2}}\sum_{\ell=1}^{\lfloor D/3\rfloor}
    \binom n\ell\ell^D\exp(KD\log(KD))
    \le n^{-D/6+o(1)},
\end{equation}
where we used $D=30d$, $K\le D^2(\log n)^{0.01}$ and $d\le(\log n)^{0.3}$. Substituting this bound into \eqref{eq-H0-first-moment-d-consant}, we obtain
\begin{equation*}
    0\le\mathbb E_0[\CMT]
    \le\sum_{A\in\binom{[K]}D}n^{-D/3+o(1)}
    \le K^D n^{-D/3+o(1)}
    =n^{-D/3+o(1)}.
\end{equation*}
This proves \eqref{eq-H0-first-moment-CMT}.

For the second-moment calculations here and later under $\mathcal H_1$ (Sec.~\ref{subsec-proof-prop-H_1-moments-CMT}), let us introduce the replicas $\beta_1,\beta_2,\beta_3 \iiddistr \Unif(\mathbb S^{d-1})$ and $\{Z_i^{(a)}: i\in[n],a=1,2,3\}$ be iid standard Gaussian, which, together with $X$, are mutually independent. 
Here $\beta_1$ represents the true $\beta_*$, while $\beta_2,\beta_3$ account for the averaging in the conditional means.
We use the shorthand
\begin{equation}
    U_i^{(a)}\equiv
    \frac{\langle X_i,\beta_a\rangle+\sigma Z_i^{(a)}}{\sqrt{1+\sigma^2}},
    \qquad
    c_{ab}\equiv\frac{\langle\beta_a,\beta_b\rangle}{1+\sigma^2}
    \quad(a\ne b),
		\label{eq:replicaU}
\end{equation}
and $c_{aa}\equiv1$. Conditional on the directions $\beta_a$'s, the $U$-variables are standard Gaussian, independent across indices, with $\mathbb E[U_i^{(a)}U_i^{(b)}\mid\beta_1,\beta_2,\beta_3]=c_{ab}$.

We will establish
\begin{align}
    \mathbb E_0[S_A^2]
    &=\mathbb E_{\beta_1,\beta_2}\left[\prod_{k\in A}c_{12}^k\right]+O(n^{-1+o(1)}),
    \label{eq-second-moment-H0-d-constant-each-term-3}\\
    |\mathbb E_0[S_AS_B]|
    &\le n^{-1+o(1)}\quad\text{if }A\ne B,
    \label{eq-second-moment-H0-d-constant-each-term-4}\\
    \mathbb E_0[T_A^2]
    &=1+O(n^{-1+o(1)}),
    \label{eq-second-moment-H0-d-constant-each-term-1}\\
    |\mathbb E_0[T_AT_B]|
    &\le n^{-1+o(1)}\quad\text{if }A\ne B.
    \label{eq-second-moment-H0-d-constant-each-term-2}
\end{align}
We prove the two estimates involving $S_A,S_B$; the estimates involving $T_A,T_B$ follow by the same argument.

For $\mathbf i_A,\mathbf i'_A\in[n]^A$, write
\begin{equation*}
    h(\mathbf i_A,\mathbf i'_A)\equiv
    \mathbb E\left[
    \prod_{k\in A}h_k(U_{i_k}^{(1)})
    \prod_{k\in A}h_k(U_{i'_k}^{(2)})
    \,\middle|\,\beta_1,\beta_2\right],
\end{equation*}
suppressing the dependence on $\beta_1,\beta_2$. Since $S_A=\mathbb E_1[T_A\mid X]$, we have
\begin{equation}\label{eq-second-moment-H0-d-constant-each-term-3-10}
    \mathbb E_0[S_A^2]
    =\mathbb E_{\beta_1,\beta_2}\left[
    \frac1{n^D}\sum_{\mathbf i_A,\mathbf i'_A\in[n]^A}
    h(\mathbf i_A,\mathbf i'_A)\right].
\end{equation}
We partition the pairs $(\mathbf i_A,\mathbf i'_A)\in[n]^A\times[n]^A$ into three classes (occurrences counted jointly in $\mathbf i_A$ and $\mathbf i'_A$):
\begin{itemize}
    \item $\mathcal R_1$ (singletons): at least one index occurs exactly once;
    \item $\mathcal R_2$ (exact pairings): every appearing index occurs exactly twice;
    \item $\mathcal R_3$ (over-represented): every appearing index occurs at least twice, and at least one index occurs three or more times.
\end{itemize}
Thus
\begin{equation*}
    \sum_{\mathbf i_A,\mathbf i'_A\in[n]^A}h(\mathbf i_A,\mathbf i'_A)
    =\sum_{r=1}^3
    \sum_{(\mathbf i_A,\mathbf i'_A)\in\mathcal R_r}h(\mathbf i_A,\mathbf i'_A).
\end{equation*}
We bound the three contributions separately and uniformly in $\beta_1,\beta_2$.

For $\mathcal R_1$, independence across indices and centering imply that each summand vanishes:
\begin{equation}\label{eq-second-moment-H0-d-constant-each-term-3-2}
    h(\mathbf i_A,\mathbf i'_A)=0, \quad (\mathbf i_A,\mathbf i'_A)\in\mathcal R_1.
\end{equation}

For $\mathcal R_2$, by Hermite orthogonality, nonzero contributions require that all indices must be paired exactly: $i_k=i'_k$ for every $k\in A$. The total number of ways to choose the distinct indices in these $D$ pairs is $(n)_D\equiv n(n-1)\cdots(n-D+1)$. Applying Lemma~\ref{lem:hermite_cov}, we obtain
\begin{equation}\label{eq-second-moment-H0-d-constant-each-term-3-3}
\begin{aligned}
    \sum_{(\mathbf i_A,\mathbf i'_A)\in\mathcal R_2}
    h(\mathbf i_A,\mathbf i'_A)
    &=(n)_D\prod_{k\in A}
    \mathbb E[h_k(U_1^{(1)})h_k(U_1^{(2)})\mid\beta_1,\beta_2]\\
    &=(n)_D\prod_{k\in A}c_{12}^k
    =\bigl(1+O(D^2/n)\bigr)n^D\prod_{k\in A}c_{12}^k.
\end{aligned}
\end{equation}

For $\mathcal R_3$, $\mathbf i_A$ and $\mathbf i'_A$ together contain at most $D-1$ distinct indices. By H\"older's inequality and Lemma~\ref{lem:fourth-moment-hermite-estimate},
\begin{equation*}
    |h(\mathbf i_A,\mathbf i'_A)|
    \le\prod_{k\in A}
    \mathbb E_0\left[\left|h_k\left(\frac{Y_1}{\sqrt{1+\sigma^2}}\right)\right|^{2D}\right]^{1/D}
    \le\exp(2DK\log(2D)).
\end{equation*}
There are at most $\binom n\ell\ell^{2D}$ pairs of tuples using exactly $\ell$ distinct indices. Therefore,
\begin{equation}\label{eq-second-moment-H0-d-constant-each-term-3-4}
    \sum_{(\mathbf i_A,\mathbf i'_A)\in\mathcal R_3}
    |h(\mathbf i_A,\mathbf i'_A)|
    \le\sum_{\ell=1}^{D-1}\binom n\ell\ell^{2D}\exp(2DK\log(2D))
    \le n^{D-1+o(1)}.
\end{equation}
Substituting the three contributions into \eqref{eq-second-moment-H0-d-constant-each-term-3-10} proves \eqref{eq-second-moment-H0-d-constant-each-term-3}. For $\mathbb E_0[T_A^2]$, repeat the argument with $U_i^{(1)}$ and $U_i^{(2)}$ both replaced by $Y_i/\sqrt{1+\sigma^2}$. Each surviving pair contributes $\mathbb E_0[h_k(Y_1/\sqrt{1+\sigma^2})^2]=1$, proving \eqref{eq-second-moment-H0-d-constant-each-term-1}.

For $A\ne B$, use $\mathbf i_A\in[n]^A$ and $\mathbf j_B\in[n]^B$, and write
\begin{equation*}
    h(\mathbf i_A,\mathbf j_B)\equiv
    \mathbb E\left[
    \prod_{k\in A}h_k(U_{i_k}^{(1)})
    \prod_{\ell\in B}h_\ell(U_{j_\ell}^{(2)})
    \,\middle|\,\beta_1,\beta_2\right].
\end{equation*}
The corresponding expansion is
\begin{equation}\label{eq-second-moment-H0-d-constant-each-term-4-1}
    \mathbb E_0[S_AS_B]
    =\mathbb E_{\beta_1,\beta_2}\left[
    \frac1{n^D}
    \sum_{\mathbf i_A\in[n]^A}
    \sum_{\mathbf j_B\in[n]^B}
    h(\mathbf i_A,\mathbf j_B)\right].
\end{equation}
For pairs $(\mathbf i_A,\mathbf j_B)\in[n]^A\times[n]^B$, define $\mathcal R_1,\mathcal R_2,\mathcal R_3$ by the same multiplicity conditions. The contribution from $\mathcal R_1$ again vanishes by independence and centering. For $\mathcal R_2$, a nonzero summand would require every degree in $A$ to be paired with the same degree in $B$. Since $A\ne B$, this is impossible, and
\begin{equation}\label{eq-second-moment-H0-d-constant-each-term-4-2}
    h(\mathbf i_A,\mathbf j_B)=0, \quad (\mathbf i_A,\mathbf j_B) \in \calR_2.
\end{equation}
For $\mathcal R_3$, $\mathbf i_A$ and $\mathbf j_B$ together contain at most $D-1$ distinct indices. Applying the same counting and moment bounds as in \eqref{eq-second-moment-H0-d-constant-each-term-3-4}, we obtain
\begin{equation}\label{eq-second-moment-H0-d-constant-each-term-4-3}
    \sum_{(\mathbf i_A,\mathbf j_B)\in\mathcal R_3}
    |h(\mathbf i_A,\mathbf j_B)|
    \le n^{D-1+o(1)}.
\end{equation}
Substituting the three contributions into \eqref{eq-second-moment-H0-d-constant-each-term-4-1} proves \eqref{eq-second-moment-H0-d-constant-each-term-4}. The same argument gives \eqref{eq-second-moment-H0-d-constant-each-term-2}.

Finally, independence under $\mathcal H_0$ gives
\begin{align*}
    \mathbb E_0[\CMT^2]
    &=\sum_{A,B\in\binom{[K]}D}\mathbb E_0[T_AT_B]\mathbb E_0[S_AS_B]\\
    &=\sum_{A\in\binom{[K]}D}
    \left(\mathbb E_{\beta_1,\beta_2}\left[\prod_{k\in A}c_{12}^k\right]+O(n^{-1+o(1)})\right)
    +\sum_{\substack{A,B\in\binom{[K]}D\\A\ne B}}O(n^{-1+o(1)})\\
    &=\sum_{A\in\binom{[K]}D}
    \mathbb E_{\beta_1,\beta_2}\left[\prod_{k\in A}c_{12}^k\right]+O(n^{-1+o(1)}),
\end{align*}
where the errors are uniform in $A,B$ and $K^{2D}=n^{o(1)}$. This proves \eqref{eq-H0-second-moment-CMT}.

\subsubsection{Proof of Proposition~\ref{prop-H_1-moments-CMT}}
\label{subsec-proof-prop-H_1-moments-CMT}
We start from the first moment under $\mathcal{H}_1$. For each $A\in\binom{[K]}{D}$, we have
$$
\mathbb E_1\left[ T_AS_A \right]
=\mathbb E_1\left[\mathbb E_1\left[ T_A \mid X \right] S_A\right]
=\mathbb E_1\left[S^2_A \right]=
\mathbb E_0\left[S^2_A \right],
$$
where the second equality follows from the definition of $S_A=
\mathbb E_1\left[ T_A \mid X \right]$ and the last equality holds because $X$ has the same distribution under $\calH_1$
and $\calH_0$. It then follows from \eqref{eq-second-moment-H0-d-constant-each-term-3} that 
\begin{equation}\label{eq-second-moment-H1-CMT-each-term}
    \mathbb E_1[T_A S_A] = \mathbb E_{\beta_1,\beta_2}\left[ \prod_{k\in A}\left(\frac{\langle\beta_1,\beta_2\rangle }{1+\sigma^2}\right)^k \right]+O(n^{-1+o(1)})\,.\notag
\end{equation}
Summing over $A\in\binom{[K]}{D}$ yields \eqref{eq-H1-first-moment-CMT}. 

Next, we turn to the second moment under $\mathcal H_1$. We use the three replicas introduced in \prettyref{eq:replicaU}: $U^{(1)}$ represents the observations in both $T_A$ and $T_B$, while $U^{(2)}$ and $U^{(3)}$ account for $S_A$ and $S_B$, respectively. Fix $A,B\in\binom{[K]}{D}$, and let $\mathbf i_A,\mathbf i'_A\in[n]^A$ and $\mathbf j_B,\mathbf j'_B\in[n]^B$ be the index tuples corresponding to $T_A,S_A,T_B,S_B$. Conditional on $\beta_1,\beta_2,\beta_3$, write
\begin{equation*}
h(\mathbf i_A,\mathbf i'_A,\mathbf j_B,\mathbf j'_B)
\triangleq
\mathbb E\left[
\prod_{k\in A}h_k(U_{i_k}^{(1)})h_k(U_{i'_k}^{(2)})
\prod_{\ell\in B}h_\ell(U_{j_\ell}^{(1)})h_\ell(U_{j'_\ell}^{(3)})
\right].
\end{equation*}
Here the expectation is over $X$ and the noise, with $\beta_1,\beta_2,\beta_3$ held fixed. Expanding the four factors and averaging over the directions gives
\begin{equation}\label{eq-second-moment-H1-each-term-d-constant-non-intersect-1}
\mathbb E_1[T_AS_AT_BS_B]
=
\frac{1}{n^{2D}}\mathbb E_{\beta_1,\beta_2,\beta_3}
\left[
\sum_{\mathbf i_A,\mathbf i'_A\in[n]^A}
\sum_{\mathbf j_B,\mathbf j'_B\in[n]^B}
h(\mathbf i_A,\mathbf i'_A,\mathbf j_B,\mathbf j'_B)
\right].
\end{equation}
Similar to \eqref{eq-second-moment-H0-d-constant-each-term-3-10}, we decompose the sum in \eqref{eq-second-moment-H1-each-term-d-constant-non-intersect-1} 
over $(\mathbf i_A,\mathbf i'_A,\mathbf j_B,\mathbf j'_B)$
as before into $\mathcal R_1,\mathcal R_2,\mathcal R_3$, with occurrences counted in the four tuples altogether.
We evaluate the three contributions separately and uniformly in the $\beta$'s.

For $\mathcal R_1$, similar to \eqref{eq-second-moment-H0-d-constant-each-term-3-2}, by 
independence and centering, we get
\begin{equation}\label{eq-second-moment-H1-each-term-d-constant-non-intersect-2}
    h(\mathbf i_A,\mathbf i'_A,\mathbf j_B,\mathbf j'_B)=0, \quad (\mathbf i_A,\mathbf i'_A,\mathbf j_B,\mathbf j'_B)\in\mathcal R_1.
\end{equation}

For $\mathcal R_2$, by orthogonality, we must have exact pairing for nonzero contribution:
For the non-overlapping part, $i_k=i'_k$ for $k\in A\setminus B$, $j_k=j'_k$ for $k\in B\setminus A$. For each $k\in A\cap B$, there are three possible pairings. Pairing $(i_k,i'_k)$ and $(j_k,j'_k)$, or pairing $(i_k,j'_k)$ and $(i'_k,j_k)$, contributes $c_{12}^kc_{13}^k$. Pairing $(i_k,j_k)$ and $(i'_k,j'_k)$ contributes $c_{23}^k$, because the first pair consists of two copies of $h_k(U_{i_k}^{(1)})$, whose second moment equals $1$. Each resulting pairing has $2D$ pairs, to which we assign distinct indices in $(n)_{2D}$ ways. Consequently,
\begin{equation}\label{eq-second-moment-H1-CMT-Mequals2}
\sum_{(\mathbf i_A,\mathbf i'_A,\mathbf j_B,\mathbf j'_B)\in\mathcal R_2}
h(\mathbf i_A,\mathbf i'_A,\mathbf j_B,\mathbf j'_B)
=
(n)_{2D}
\prod_{k\in A\setminus B}c_{12}^k
\prod_{k\in B\setminus A}c_{13}^k
\prod_{k\in A\cap B}\left(c_{23}^k+2c_{12}^kc_{13}^k\right).
\end{equation}

For $\mathcal R_3$, $\mathbf i_A,\mathbf i'_A,\mathbf j_B,\mathbf j'_B$ together contain at most $2D-1$ distinct indices. Applying the same counting argument and Hermite moment bound as in~\eqref{eq-second-moment-H0-d-constant-each-term-3-4}, we obtain
\begin{equation}\label{eq-second-moment-H1-d-const-non-intersect-geq3}
\sum_{(\mathbf i_A,\mathbf i'_A,\mathbf j_B,\mathbf j'_B)\in\mathcal R_3}
\left|h(\mathbf i_A,\mathbf i'_A,\mathbf j_B,\mathbf j'_B)\right|
\le n^{2D-1+o(1)}.
\end{equation}
Substituting the three contributions into \eqref{eq-second-moment-H1-each-term-d-constant-non-intersect-1} and using $(n)_{2D}/n^{2D}=1+O(D^2/n)$, we conclude that
\begin{equation}
\mathbb E_1[T_AS_AT_BS_B]
=
\mathbb E_{\beta_1,\beta_2,\beta_3}\left[
\prod_{k\in A\setminus B}c_{12}^k
\prod_{k\in B\setminus A}c_{13}^k
\prod_{k\in A\cap B}\left(c_{23}^k+2c_{12}^kc_{13}^k\right)
\right]
+O(n^{-1+o(1)})
\label{eq:TASATBSB}
\end{equation}

Let $g(A,B)$ denote the expectation on the right-hand side. 
We sum $g(A,B)$ according to the overlap size $t=|A\cap B|$. For $t=0$, conditioning on $\beta_1$ and using independence of $\beta_2,\beta_3$, we obtain
\begin{align*}
g(A,B)
&=
\mathbb E_{\beta_1,\beta_2,\beta_3}\left[
\prod_{k\in A}c_{12}^k\prod_{k\in B}c_{13}^k
\right]\\
&=
\mathbb E_{\beta_1}\left[
\mathbb E_{\beta_2}\left[\prod_{k\in A}c_{12}^k\mid\beta_1\right]
\mathbb E_{\beta_3}\left[\prod_{k\in B}c_{13}^k\mid\beta_1\right]
\right]\\
&=
\mathbb E_{\beta_1,\beta_2}\left[\prod_{k\in A}c_{12}^k\right]
\mathbb E_{\beta_1,\beta_3}\left[\prod_{k\in B}c_{13}^k\right]
\end{align*}
where the last equality follows from rotational invariance so that each conditional expectation is independent of $\beta_1$. The first factor is zero if $\sum_{k\in A}k$ is odd and nonnegative if it is even; the same holds for the second factor with $B$. We may therefore drop the disjointness restriction and obtain
\begin{align}
\sum_{\substack{A,B\in\binom{[K]}{D}\\A\cap B=\emptyset}}g(A,B)
\le
\left(\sum_{A\in\binom{[K]}{D}}
\mathbb E_{\beta_1,\beta_2}\left[\prod_{k\in A}c_{12}^k\right]\right)^2
=
\left(\mathbb E_1[\CMT]\right)^2+O(n^{-1+o(1)}), \label{eq:2nd-disjoint}
\end{align}
where the last equality follows from~\eqref{eq-H1-first-moment-CMT}.

For $1\le t\le D$, each factor indexed by $A\cap B$ has absolute value at most $3$. We therefore obtain
\begin{align}
\sum_{\substack{A,B\in\binom{[K]}{D}\\|A\cap B|=t}}g(A,B)
&\le
3^t\sum_{\substack{A,B\in\binom{[K]}{D}\\|A\cap B|=t}}
\mathbb E_{\beta_1,\beta_2,\beta_3}\left[
\prod_{k\in A\setminus B}|c_{12}|^k
\prod_{k\in B\setminus A}|c_{13}|^k
\right]\notag\\
&\le
3^t\binom Kt
\left(\sum_{C\in\binom{[K]}{D-t}}
\mathbb E_{\beta_1,\beta_2}\left[\prod_{k\in C}|c_{12}|^k\right]\right)^2.
\label{eq-summ-over-each-t-second-moment-H1-d-grow-slowly-1}
\end{align}
Finally, $\binom Kt\le K^t/t!$ and
\begin{equation}
\sum_{C\in\binom{[K]}{D-t}}
\mathbb E_{\beta_1,\beta_2}\left[\prod_{k\in C}|c_{12}|^k\right]
\le
\frac{1}{(D-t)!}
\mathbb E_{\beta_1,\beta_2}\left[
\left(\sum_{k=1}^K|c_{12}|^k\right)^{D-t}
\right].
\label{eq:2nd-nondisjoint}
\end{equation}
Combining \prettyref{eq:TASATBSB} and \prettyref{eq:2nd-disjoint}--\prettyref{eq:2nd-nondisjoint}, we obtain as desired
\begin{align*}
\var_1[\CMT]
&= \sum_{A,B\in\binom{[K]}{D}} \mathbb E_1[T_AS_AT_BS_B] - \left(\mathbb E_1[\CMT]\right)^2\\
&\le
\sum_{t=1}^D\sum_{\substack{A,B\in\binom{[K]}{D}\\|A\cap B|=t}}g(A,B)
+O(n^{-1+o(1)})\\
&\le
\sum_{t=1}^D\frac{(3K)^t}{[(D-t)!]^2t!}
\left(\mathbb E_{\beta_1,\beta_2}\left[
\left(\sum_{k=1}^K
\left(\frac{|\langle\beta_1,\beta_2\rangle|}{1+\sigma^2}\right)^k
\right)^{D-t}\right]\right)^2
+O(n^{-1+o(1)}).
\end{align*}

\subsection{Proof of Lemmas~\ref{lem:est-ratio-Dminust-versus-D} and \ref{lem:est-sum-first-second-moment-growing-d}}
\label{subsec-two-lemmas}
\begin{proof}[Proof of Lemma~\ref{lem:est-ratio-Dminust-versus-D}]
Let 
\begin{equation}
    \label{eq:def-g_0-g}
    g_0(x) = \sum_{k=1}^K x^k=\frac{x(1-x^K)}{1-x}\text{ and }g(x)=g_0(x/(1+\sigma^2)).
\end{equation}
We start with the case $d=1$. In this case, $\langle\beta_1,\beta_2\rangle$ is uniformly distributed on the set $\{-1,1\}$. Therefore,
\begin{equation}\label{eq:mu-1d}
    \mu_{d,K,D}=\mathbb E_{\beta_1,\beta_2}[g(\langle\beta_1,\beta_2\rangle)^D]
    = \frac{1}{2} g^D(1) + \frac{1}{2} g^D(-1) 
    \geq \frac{1}{2} \left(K e^{-K\sigma^2}\right)^D   -\frac{1}2 \to +\infty \,,
\end{equation}
where the first inequality holds because
$g(-1) \ge -1$ and $g(1) \ge K (1+\sigma^2)^{-K} \ge K e^{-K\sigma^2}$;
the last line follows from $ D^2 \ll K\ll 1/\sigma^2$ under our choice of $K, D$. This completes the proof of \eqref{eq:mean_mu_lower_bound}. Similarly,
\begin{equation}\label{eq:nu-1d}
    \nu_{d,K,D}(t)=\mathbb E_{\beta_1,\beta_2}[g(|\langle\beta_1,\beta_2\rangle|)^{D-t}] = g^{D-t}(1) \leq K^{D-t}.
\end{equation}
Furthermore, combining \eqref{eq:nu-1d} and \eqref{eq:mu-1d} yields that for $1 \le t \le D-5d$
\begin{equation}
    \frac{\nu_{d,K,D}(t)}{\mu_{d,K,D}} \leq \frac{2g^{D-t}(1)}{g^D(1)-1 } 
    \le 4 g^{-t}(1) \le 4 \left( K e^{-K\sigma^2} \right)^{-t}
    \le 4 \left( K e^{-1} \right)^{-t}. \notag
\end{equation}
This shows \eqref{eq-est-first-moment-growing-d-goal1-case1}.

For $d\geq 2$, let $f(\rho)$ denote the probability density function of $\iprod{\beta_1}{\beta_2}$, that is, 
$$
f(\rho) \triangleq c_d (1-\rho^2)^{(d-3)/2}, \forall \rho \in [-1,1],  \text{ where } c_d=\frac{1}{B\left(\frac{d-1}{2},\frac{1}{2}\right)}.
$$ 
where $B(x,y)=\int_{0}^{1} t^{x-1}(1-t)^{y-1} dt$ is  the beta function. Note that 
$B\left(\frac{d-1}{2},\frac{1}{2}\right) \le 
B\left(\frac{1}{2},\frac{1}{2}\right)=\pi$,
so that $c_d \ge 1/\pi.$ 
For $d\ge 3$, $f$ is nonincreasing on $[0,1]$ and $c_d(3x/2)^{(d-3)/2}\le f(1-x)\le c_d(2x)^{(d-3)/2}$ for $0<x\le 1/2$. For $d=2$, $f$ is increasing on $[0,1]$, so below we instead use $f(1-\frac{i+1}{K})\le f(1-\frac iK)\le c_2(K/i)^{1/2}$ and $f(1-\frac1K)\ge f(1-\frac2K)\ge c_2(K/4)^{1/2}$; the resulting bounds differ from the case $d\ge3$ by at most a factor of $2$, which is absorbed in \eqref{eq-est-first-moment-growing-d-goal1-case1-numerator-first-second--ratio} and by the constant $1/(2\pi)$ in \eqref{eq:mean_mu_lower_bound}.

We first bound $\nu_{d,K,D}(t)$.  Consider two cases regarding the values of $t$: (i) $1\leq t\leq D-5d$; (ii) $D-5d\leq t\leq D$. We start from case (i), getting that 
\begin{equation}\label{eq-est-first-moment-growing-d-goal1-case1-numerator}
    \begin{split}
       \nu_{d,K,D}(t) \triangleq & \mathbb E_{\beta_1,\beta_2} \left[ g \left( |\langle\beta_1,\beta_2\rangle| \right)^{D-t}\right] 
        = 2\int_0^{1-\frac{10}{K}}g \left( \rho \right)^{D-t}f(\rho)\mathrm{d}\rho+2\int_{1-\frac{10}{K}}^1 g\left(\rho \right)^{D-t}f(\rho)\mathrm{d}\rho\,.
    \end{split}
\end{equation}
For the first term in \eqref{eq-est-first-moment-growing-d-goal1-case1-numerator},
\begin{equation}\label{eq-est-first-moment-growing-d-goal1-case1-numerator-first-term}
\begin{split}
    \int_0^{1-\frac{10}{K}} g\left( \rho \right)^{D-t}f(\rho)\mathrm{d}\rho 
    = & \sum_{i=10}^{K-1}\int_{1-\frac{i+1}{K}}^{1-\frac{i}{K}} g \left( \rho \right)^{D-t} f(\rho) \mathrm{d}\rho \\
    \leq& \sum_{i=10}^{K-1}\frac{1}{K} g\left( 1-\frac{i}{K} \right)^{D-t}f\left(1-\frac{i+1}{K}\right)\\
    \leq & \sum_{i=10}^{K-1} \frac{1}{K} \left(\frac{K}{i}\right)^{D-t} c_d \left(\frac{2(i+1)}{K}\right)^{\frac{d-3}{2}}\\ 
    \leq & c_d K^{D-t-\frac{d-1}{2}} \sum_{i=10}^{K-1} \frac{3^{\frac{d-3}{2}}}{i^{D-t-\frac{d-3}{2}}} \\
     \leq  & c_d K^{D-t-\frac{d-1}{2}}  10^{-(D-t-\frac{d-1}{2})}\cdot 3^{\frac{d-3}{2}}\,,
\end{split}\notag
\end{equation}
where we use the fact that $f(1-x)\le c_d (2x)^{(d-3)/2}$  , $\sum_{i\ge 10} i^{-a}\le 10^{-a}+\int_{10}^\infty x^{-a}\,dx\le 10^{-(a-1)}$ for $a=D-t-\frac{d-3}{2}\ge 6$, and
\begin{align*}
    g\left(1-\frac{i}{K}\right) \le g_0\left(1-\frac{i}{K}\right)= \frac{ \left(1-\frac{i}{K} \right)\left(1-\left(1-\frac{i}{K}\right)^K \right)}{\frac{i}{K} }\leq \frac{K}{i}.
\end{align*}

For the second term in \eqref{eq-est-first-moment-growing-d-goal1-case1-numerator}, 
\begin{equation}
    \label{eq-est-first-moment-growing-d-goal1-case1-numerator-second-term}
    \begin{split}
     \int_{1-\frac{10}{K}}^1 g\left( \rho \right)^{D-t}f(\rho)\mathrm{d}\rho 
     &\geq \int_{1-\frac{2}{K}}^{1-\frac{1}{K}} g \left(\rho\right)^{D-t}f(\rho)\mathrm{d}\rho
     \geq \frac{1}{K} g\left( 1-\frac{2}{K} \right)^{D-t}f\left(1-\frac{1}{K}\right) \\
    &\geq \frac{1}{K}\left( \frac{1-e^{-2}}{3}K \right)^{D-t} \cdot c_d \left(\frac{3}{2K}\right)^{\frac{d-3}{2}}\\
    &= c_d K^{D-t-\frac{d-1}{2}}  \left( \frac{1-e^{-2}}{3} \right)^{D-t}\left(\frac{3}{2}\right)^{\frac{d-3}{2}}\,,
    \end{split}\notag
\end{equation}
where the last inequality holds because
$f(1-x) \ge c_d (3x/2)^{(d-3)/2}$ for $x \le 1/2$
and
\begin{align}\label{eq:low-bound-of-g}
    g\left(1-\frac{i}{K}\right)
    = \frac{\left(1-\frac{i}{K}\right)\left(1-\left(\frac{1-\frac{i}{K}}{1+\sigma^2}\right)^K \right)}{\sigma^2+\frac{i}{K}}
    \geq \frac{\left(1-\frac{i}{K}\right)\left(1-\left(1-\frac{i}{K}\right)^K \right)}{\sigma^2+\frac{i}{K}} \geq
    \frac{(1-e^{-i})K}{i+1},
\end{align}
for all $i$ such that $i(i+1) \le K(1-K\sigma^2).$
Therefore,  
\begin{equation}
    \label{eq-est-first-moment-growing-d-goal1-case1-numerator-first-second--ratio}
    \begin{split}
     \frac{\int_0^{1-\frac{10}{K}} g\left( \rho \right)^{D-t}f(\rho)\mathrm{d}\rho }{\int_{1-\frac{10}{K}}^1 g \left(\rho \right)^{D-t}f(\rho)\mathrm{d}\rho}
    \leq \left( \frac{3}{10(1-e^{-2})} \right)^{D-t}20^{\frac{d-1}{2}}
    \leq  \left( \left( \frac{3}{10(1-e^{-2})} \right)^5\sqrt{20} \right)^d<1\,.
    \end{split}
\end{equation}
It follows from~\eqref{eq-est-first-moment-growing-d-goal1-case1-numerator} that 
$$
\nu_{d,K,D}(t) \le 4\int_{1-\frac{10}{K}}^1 g\left(\rho \right)^{D-t}f(\rho)\mathrm{d}\rho\,.
$$
Moreover, 
\begin{equation}\label{eq-est-first-moment-growing-d-goal1-case1-denominator}
\begin{split}
   \mu_{d,K,D} &\triangleq \int_{-1}^1 g\left(\rho \right)^{D}f(\rho)\mathrm{d}\rho
    \geq \int_{1-\frac{10}{K}}^1 g\left( \rho \right)^{D}f(\rho)\mathrm{d}\rho \\
    &\geq g\left( 1- \frac{10}{K}\right)^{t} \int_{1-\frac{10}{K}}^1 g \left(\rho\right)^{D-t}f(\rho)\mathrm{d}\rho
    \geq \left( \frac{1-e^{-10}}{11}K \right)^{t}\int_{1-\frac{10}{K}}^1 g \left(\rho \right)^{D-t}f(\rho)\mathrm{d}\rho\,,
\end{split}
\end{equation}
where the last inequality holds because 
$g(1-10/K) \ge (1-e^{-10})K/11$ in view of \eqref{eq:low-bound-of-g}.
Combining the bounds on $\nu_{d,K,D}(t),\mu_{d,K,D}$ yields the desired \eqref{eq-est-first-moment-growing-d-goal1-case1}.

For case (ii) $t \geq D-5d$, we have
\begin{equation}
\begin{split}
    \mathbb E_{\beta_1,\beta_2} \left[ g \left(|\langle\beta_1,\beta_2\rangle| \right)^{D-t}\right]\leq K^{D-t}\leq K^{5d}\,.
\end{split}\notag
\end{equation}

Finally, 
\begin{equation}
\begin{split}
    \mu_{d,K,D}  \geq \int_{1-\frac{2}{K}}^{1-\frac{1}{K}} g \left( \rho \right)^{D}f(\rho)\mathrm{d}\rho 
    \geq \frac{1}{K} g \left( 1- \frac{2}{K}\right)^D f\left(1 - \frac{1}{K}\right) 
    \ge 
    \frac{c_d}{ K}\left(\frac{1-e^{-2}}{3}K\right)^D\frac{1}{K^{\frac{d-3}{2}}} \,.
\end{split}\notag
\end{equation}

This concludes the proof of~\prettyref{eq:mean_mu_lower_bound}.

\end{proof}

\begin{proof}[Proof of Lemma~\ref{lem:est-sum-first-second-moment-growing-d}]
We define $\mathsf A(K,D)$ to be the left-hand side of \eqref{eq-est-first-moment-growing-d}, that is, 
$$
\mathsf A(K,D) \triangleq \sum_{A\in\binom{[K]}{D}}\mathbb E_{\beta_1,\beta_2}\left[ \prod_{k\in A}\left(\frac{\langle\beta_1,\beta_2\rangle }{1+\sigma^2}\right)^k \right] 
$$
and let 
\begin{equation}
\mathsf B(K,D)\triangleq \mathbb E_{\beta_1,\beta_2} \left[\sum_{1\leq k_1\leq k_2\leq \cdots\leq k_D\leq K}\left(\frac{\langle \beta_1,\beta_2\rangle}{1+\sigma^2}\right)^{\sum_{i=1}^Dk_i} \right]\,.\notag
\end{equation}
Then we have 
\begin{align*}
    \mathsf A(K,D) & = 
    \frac{1}{D!} \mathbb E_{\beta_1,\beta_2} \left[\sum_{\substack{1\leq k_1, k_2, \cdots, k_D\leq K\\ k_i\neq k_j \forall i\neq j}}\left(\frac{\langle \beta_1,\beta_2\rangle}{1+\sigma^2}\right)^{\sum_{i=1}^Dk_i} \right]   \\
   &  \leq \frac{1}{D!}
    \mathbb E_{\beta_1,\beta_2} \Bigg[ \left( \sum_{k=1}^K \left(\frac{\langle \beta_1,\beta_2\rangle}{1+\sigma^2} \right)^k\right)^D \Bigg]
    \leq \mathsf B(K,D)\,.\notag
\end{align*}
  Here both inequalities hold because each term is a nonnegative multiple of $\mathbb E[\langle\beta_1,\beta_2\rangle^m]\ge0$, as odd moments vanish by symmetry. Thus, it suffices to show  that 
\begin{equation}\label{eq-est-first-moment-growing-d-ultimate-goal}
    \mathsf B(K,D)-\mathsf A(K,D) = \frac{o(1)}{D!}\mathbb E_{\beta_1,\beta_2} \Bigg[ \left( \sum_{k=1}^K \left(\frac{\langle \beta_1,\beta_2\rangle}{1+\sigma^2} \right)^k\right)^D \Bigg]
    = \frac{o\left(1\right)}{D!} \mu_{d,K,D}  \,.
\end{equation}
Specifically, we have
\begin{equation}
\begin{split}
 \mathsf B(K,D)-\mathsf A(K,D)
 =  \sum_{\substack{1\leq k_1\leq k_2\leq \cdots\leq k_D\leq K\\\exists i, k_i= k_{i+1} }}
 \mathbb E_{\beta_1,\beta_2} \left[\left(\frac{\langle \beta_1,\beta_2\rangle}{1+\sigma^2}\right)^{\sum_{i=1}^Dk_i} \right]
  = & \sum_{t=1}^{D-1} \mathsf C_t \,,
\end{split}\notag
\end{equation}
where $t$ denotes the number of positions where the inequalities hold as equalities, and 
\begin{equation*}
\begin{split}
    \mathsf C_t & \triangleq \sum_{\substack{1\leq i_1<i_2<\dots<i_t\leq D-1}} \sum_{\substack{1\leq k_1\leq k_2\leq \cdots\leq k_D\leq K:\\k_{i_s}=k_{i_s+1}\forall 1\leq s\leq t\\ \forall j \notin \{i_1,\dots,i_t\},\,k_j<k_{j+1} }}
    \mathbb E_{\beta_1,\beta_2} \left[\left(\frac{\langle \beta_1,\beta_2\rangle}{1+\sigma^2}\right)^{\sum_{i=1}^Dk_i} \right] \\
    & \le \binom{D-1}{t}\times \sum_{1\leq k_1<\cdots<k_{D-t}\leq K}  \mathbb E_{\beta_1,\beta_2}\left[\left( \frac{|\langle\beta_1,\beta_2\rangle|}{1+\sigma^2} \right)^{\sum_{i=1}^{D-t}k_i}\right] \\
    & \leq \binom{D-1}{t}\frac{1}{(D-t)!}\mathbb E_{\beta_1,\beta_2} \left[ \left(\sum_{k=1}^K \left(\frac{|\langle\beta_1,\beta_2\rangle|}{1+\sigma^2}\right)^k \right)^{D-t}\right]\,.\notag
\end{split}
\end{equation*}
By Lemma~\ref{lem:est-ratio-Dminust-versus-D}, for $t\leq D-5d$,
\begin{equation}\label{eq-Ct-tsmall}
\begin{split}
    \mathsf C_t\leq   \binom{D-1}{t}\frac{1}{(D-t)!} \times 4 \left( \frac{1-e^{-10}}{11} K \right)^{-t} \mu_{d,K,D}
    \leq  \left( O\left(\frac{D
    ^2}{K}\right) \right)^{t} \frac{\mu_{d,K,D}}{D!}  \,.
\end{split}
\end{equation}
For $D-5d\leq t \leq D-1$,
\begin{align*}
 \mathsf C_t\leq \binom{D-1}{t}\frac{1}{(D-t)!} K^{5d} 
 \leq \frac{D^{D-1}}{D!} K^{5d}
 \leq  \frac{D^{D-1}}{K^{D-6d}} \frac{(O(1))^D}{\pi} \frac{\mu_{d,K,D} }{D!} \,,
 \end{align*}
 where the second inequality holds because
 $\binom{D-1}{t} \le D^{D-1-t}$ and $D!/(D-t)! \le D^t$ and the last inequality holds in view of the lower bound of $\mu_{d,K,D}$ in~\prettyref{eq:mean_mu_lower_bound}. 
Therefore, summing over $t\geq 1$ yields the desired \eqref{eq-est-first-moment-growing-d-ultimate-goal}:
$$
\sum_{t=1}^{D-1} 
\mathsf C_t
\le \sum_{1 \le t \le D-5d} \left( O\left(\frac{D
    ^2}{K}\right) \right)^{t} 
    \frac{\mu_{d,K,D}}{D!}
    + 
    \frac{D^{D}}{K^{D-6d}} \frac{(O(1))^D}{\pi} \frac{\mu_{d,K,D} }{D!}
    =\frac{o\left(1\right)}{D!} \mu_{d,K,D}, 
$$
where the last inequality holds by our choice of $K,D$. 
\end{proof}

\section{Proofs of Theorems~\ref{thm:SMT}--\ref{thm:sparse-beta}}
\label{sec:remaining-proofs}

\subsection{Proof of Theorem~\ref{thm:SMT}}
\label{sec:max-over-beta-testing-statistic}

 We first establish the computational guarantee of $\SMT$. For each $\beta \in N_\delta$, computing $S_k(\beta)$ and $T_k$ requires $n^{1+o(1)}$ time. Moreover, $|N_\delta| \le (3/\delta)^d =n^{o(1)}$, as $\delta=1/\sqrt{M}$ and $d\log M=o(\log n)$ (as $M\le\log n$).
Therefore, the total run time to compute $\SMT$ is $n^{1+o(1)}$. In the following, we separately bound the testing errors of $\SMT$ under $\calH_1$
and $\calH_0.$

\paragraph*{Analysis under $\calH_1$}
Without loss of generality, we assume the true permutation $\pi_*=\id.$
We aim to show that for each possible realization of $\beta_*$,
    \begin{equation}
        \mathbb P_1\left[ \SMT \leq \mathsf c_*M \,\big|\,\beta_*=\beta\right]=o(1)\,,
        \label{eq:SMT_type_1_error}
    \end{equation}
    where $\mathsf c_*=e^{-1}-e^{-\frac{3}{2}}$.
    By the definition of $N_\delta$, there exists a $\widehat{\beta}$ such that $\|\beta-\widehat{\beta}\| \leq \delta$, or equivalently $\langle\beta,\widehat \beta\rangle \geq 1-\frac{\delta^2}{2} = 1-\frac{1}{2M}. $
    Denoting $\mathbb P_{1,\beta}\left[\cdot\right]=\mathbb P_1\left[\cdot\vert \beta_*=\beta\right]$ and writing $\SMT(\beta) = \sum_{k=1}^M T_k  S_k(\beta) \le \SMT$, it suffices to prove
    \begin{equation}\label{eq-final-goal-H1-d-large}
        \mathbb P_{1,\beta}\left[ \SMT(\widehat\beta) \leq \mathsf c_*M \right]\leq o(1)\,.
    \end{equation}
     We apply the second moment method. 
    For the first moment, we have
    \begin{align*}
    \mathbb E_{1,\beta}\left[T_k  S_k(\widehat\beta)\right]
    = & \frac{1}{(1+\sigma^2)^{\frac{k}{2}}} \mathbb E_{1,\beta} \left[ \left( \frac{1}{\sqrt{n}}\sum_{j=1}^nh_k\left(\langle X_j,\widehat{\beta}\rangle\right) \right)\left( \frac{1}{\sqrt{n}}\sum_{j=1}^nh_k \left(\frac{\langle X_j,\beta\rangle+\sigma Z_{j}}{\sqrt{1+\sigma^2}}\right) \right)\right] \\
       = & \frac{ \langle\beta,\widehat\beta\rangle^k}{(1+\sigma^2)^{k}}
      ,
    \end{align*}
    where the last equality follows from~\prettyref{lem:hermite_cov}. Therefore, 
    \begin{equation}
        \begin{split}
             \mathbb E_{1,\beta}\left[ \SMT(\widehat\beta) \right] =& \sum_{k=1}^M \mathbb E_{1,\beta}\left[T_k  S_k(\widehat\beta)\right]=
             \sum_{k=1}^M 
             \frac{ \langle\beta,\widehat\beta\rangle^k}{(1+\sigma^2)^{k}}
            \geq \frac{1}{(1+\sigma^2)^M} 
            \sum_{k=1}^M 
             \left(1- \frac{1}{2M}\right)^k
            \\
             = & \frac{\left( 2M-1 \right)}{(1+\sigma^2)^M}  \left(1-\left( 1-\frac{1}{2M} \right)^M\right)\geq \frac{3\mathsf c_*M}{2}\,,
        \end{split}
    \end{equation}
    where the last inequality follows from $\left( 1-\frac{1}{2M} \right)^M\leq e^{-\frac{1}{2}}$ and $(1+\sigma^2)^M\leq e^{\sigma^2M} \le e$  by our choice of $M$ in \eqref{eq-def-M-d-large}. 
    
    For the second moment, we have
    \begin{align*}
    & \mathbb E_{1,\beta}\left[ T_k  S_k(\widehat\beta)T_\ell  S_\ell(\widehat\beta) \right] \\
    = & \frac{1}{(1+\sigma^2)^{\frac{k+\ell}{2}}} 
            \mathbb E_{1,\beta} \Bigg[ \left( \frac{1}{\sqrt{n}}\sum_{i=1}^nh_k\left(\langle X_i,\widehat{\beta}\rangle\right) \right)\left( \frac{1}{\sqrt{n}}\sum_{j=1}^nh_k \left(\frac{\langle X_j,\beta\rangle+\sigma Z_{j}}{\sqrt{1+\sigma^2}}\right) \right)\\
            &\qquad \qquad\qquad\qquad \times\left( \frac{1}{\sqrt{n}}\sum_{i=1}^nh_\ell\left(\langle X_i,\widehat{\beta}\rangle\right) \right)\left( \frac{1}{\sqrt{n}}\sum_{j=1}^nh_\ell \left(\frac{\langle X_j,\beta\rangle+\sigma Z_{j}}{\sqrt{1+\sigma^2}}\right) \right) \Bigg]\\
            & = 
            \frac{1}{(1+\sigma^2)^{\frac{k+\ell}{2}}} \times \frac{1}{n^2} 
            \sum_{i_k,j_k,r_\ell,s_\ell \in [n]^4}
            \bar{h}\left(i_k,j_k,r_\ell,s_\ell\right)\,,
         \end{align*}   
         where 
         \[
    \bar h_{k,\ell}\left(i_k,j_k,r_\ell,s_\ell\right)\triangleq \mathbb E_{1,\beta} \left[ h_k\left(\langle X_{i_k},\widehat{\beta}\rangle\right) h_k \left(\frac{\langle X_{j_k},\beta\rangle+\sigma Z_{j_k}}{\sqrt{1+\sigma^2}}\right)h_\ell\left(\langle X_{r_\ell},\widehat{\beta}\rangle\right) h_\ell \left(\frac{\langle X_{s_\ell},\beta\rangle+\sigma Z_{s_\ell}}{\sqrt{1+\sigma^2}}\right) \right] \,. 
    \]
        If among $(i_k, j_k, r_\ell, s_\ell)$, there is an index that appears only once, then by independence, 
        $h_{k,\ell}\left(i_k,j_k,r_\ell,s_\ell\right)=0.$
        Therefore, for $h_{k,\ell}\left(i_k,j_k,r_\ell,s_\ell\right) \neq 0$,  each index must appear at least twice. If each index appears exactly twice, then we have the following: when $k\neq \ell$, by the orthogonality of Hermite polynomials, we must have $i_k=j_k$
        and $r_\ell=s_\ell$;
        when $k=\ell$, there are three different ways to pair $(i_k,j_k,r_k,s_k).$ If all four indices are the same, that is, $i_k=j_k=r_\ell=s_\ell$ then 
        by H\"older's inequality,
    \begin{equation}
    \begin{split}
     \left| \bar{h}_{k,\ell} (i_k,j_k, r_\ell, s_\ell)\right|
        \leq  \left(\mathbb E[h^4_k(Z)]\mathbb E[h_\ell^4(Z)]\right)^{1/2} \le \exp (O(k+\ell))\,,\text{ where }Z\sim N(0,1)\,,
    \end{split}\notag
    \end{equation}
    where the last inequality follows from the moment bound of Hermite polynomials in Lemma~\ref{lem:fourth-moment-hermite-estimate}. 
It follows that for $k\neq \ell,$
$$
 \mathbb E_{1,\beta}\left[ T_k  S_k(\widehat\beta)T_\ell  S_\ell(\widehat\beta) \right] 
\le    
            \frac{n(n-1)}{n^2} \left(\frac{\langle\widehat\beta,\beta\rangle}{1+\sigma^2}\right)^{k+\ell}
            +\frac{\exp(O(k+\ell))}{n}\,;
 $$       
for $k = \ell,$
$$
 \mathbb E_{1,\beta}\left[ T_k^2  S^2_k(\widehat\beta)\right] 
\le    
            \frac{n(n-1)}{n^2} \left[ 2\left(\frac{\langle\widehat\beta,\beta\rangle}{1+\sigma^2}\right)^{2k} + 1 \right]
            +\frac{\exp(O(k+\ell))}{n} \,.
 $$       
Combining the last displayed equation yields     
    \begin{equation}\label{eq-expansion-second-moment-H1-d-large}
        \begin{split}
             \mathbb E_{1,\beta}\left[ \SMT^2(\widehat\beta) \right] = \sum_{k,\ell=1}^M \mathbb E_{1,\beta}\left[ T_k  S_k(\widehat\beta)T_\ell  S_\ell(\widehat\beta) \right]
            \le  \left( \sum_{k} \left(\frac{\langle\widehat\beta,\beta\rangle}{1+\sigma^2}\right)^{k} \right)^2 + 2M + \frac{\exp(O(M))}{n},
        \end{split}
    \end{equation}
and hence
$$
\var_{1,\beta}\left[ \SMT(\hat{\beta}) \right] = \mathbb E_{1,\beta}\left[ \SMT^2(\widehat\beta) \right]  - \left(\mathbb E_{1,\beta}\left[ \SMT(\widehat\beta) \right]\right)^2 \le 2M + \frac{\exp(O(M))}{n}.
$$
    Therefore, by Chebyshev's inequality,
    \begin{equation}
        \begin{split}
            \mathbb P_{1,\beta}\left[\SMT(\widehat\beta)< \mathsf c_*M\right]\leq \mathbb P_{1,\beta}\left[ \SMT(\widehat\beta) - \mathbb{E}_{1,\beta} \left[\SMT(\widehat\beta)\right] \le - \frac{c_*}{2} M\right]
            \leq \frac{4\operatorname{Var}_{1,\beta}\left[\SMT(\widehat\beta)\right]}{c_*^2 M^2}=o(1)\,,
        \end{split}
    \end{equation}
where the last inequality holds because $M \to \infty$ and $M=o(\log n)$ since $M\le\sqrt{\frac{d\log n}{\log\log n}}\log(d+1)=o(\log n)$ as $d=o(\log n/\log\log n)$.

 \paragraph*{Analysis under $\calH_0$}
 We next show for any fixed $\beta\in N_\delta$ and any constant $c>0,$
    \begin{equation} \label{lem:eachbeta-analysis-H0-d-large}
         \mathbb P_0\left[ \SMT(\beta) \geq cM \right] \leq \exp\left(-\Theta(M)\right)\,.
    \end{equation}
    For ease of notation, we let $T_{M+k} =(1+\sigma^2)^{\frac{k}{2}}S_k(\beta)$ for $k \in [M].$
Applying Lemma~\ref{lem:gaussian-approximation-hermite-w/correlation} with $\rho=0$, we get that 
    \begin{equation}
        \mathcal{W}_1\left\{\left[T_k\right]_{k=1}^{2M},\left[\mathsf U_k\right]_{k=1}^{2M}\right\} \leq n^{-0.49}\,,
    \end{equation}
    where $\mathcal{W}_1$ denotes the Wasserstein-1 distance, and $[\mathsf U_k]_{k=1}^{2M}$ are iid standard Gaussians.
    Thus, there is a coupling between $\mathsf{T}=[T_k]_{k=1}^{2M}$ and  $\mathsf{U}=[\mathsf U_k]_{k=1}^{2M}$ such that 
    \begin{equation}
        \mathbb E\left\{ \left\Vert \mathsf{T} -   \mathsf{U} \right\Vert_2 \right\} \leq 2 n^{-0.49}\,.
    \end{equation}
     Let $\mathsf{\Delta}=\mathsf{T}-\mathsf{U}$ and $\mathcal{G_{\mathsf{cp}}}$ denote the event that  $\|\mathsf{\Delta}\|_2 \le n^{-0.3}$. Then, by Markov's inequality, $\prob{\mathcal{G_{\mathsf{cp}}}} \ge 1-2 n^{-0.19}$. 
     Therefore, 
    \begin{equation}\label{eq-each-beta-H0-d-large-calculation-1}
        \begin{split}
            &\mathbb P_0\left[ \SMT(\beta) \geq cM \right]\leq \mathbb P\left[ \SMT(\beta) \geq cM,\mathcal{G}_{\mathsf{cp}} \right] + \mathbb P\left[ \mathcal{G}^c_{\mathsf{cp}} \right]\\
            = & \mathbb P\left[ \sum_{k=1}^M (1+\sigma^2)^{-\frac{k}{2}}\left( \mathsf U_k+\mathsf{\Delta}_k \right)\left( \mathsf U_{k+M}+\mathsf{\Delta}_{k+M} \right)\geq cM , \mathcal{G}_{\mathsf{cp}}\right]+2n^{-0.19}\\
            \leq & \mathbb P\left[ \sum_{k=1}^M (1+\sigma^2)^{-\frac{k}{2}}\mathsf U_k\mathsf U_{k+M} \geq cM-n^{-0.6}-n^{-0.3} \| \mathsf{U}\|_2 \right]+2n^{-0.19}\,,
        \end{split}
    \end{equation}
    where in the last inequality we apply the Cauchy-Schwarz inequality to bound $\sum_{k=1}^M \mathsf{\Delta}_k  \mathsf{\Delta}_{k+M} \le \|\mathsf{\Delta}\|_2^2 $ and  $\sum_{k=1}^M \mathsf{\Delta}_k \mathsf{U}_{k+M} +\sum_{k=1}^M \mathsf{\Delta}_{k+M} \mathsf{U}_{k} \le  \|\mathsf{\Delta}\|_2  \|\mathsf{U}\|_2.$
    By Bernstein's inequality for the sum of sub-exponential random variables
    (see e.g.~\cite[Theorem~2.8.1]{vershynin2018high}),
    \begin{align}
        &\mathbb P\left[ \| \mathsf{U}\|_2 \geq n^{0.1} \right]\leq \exp(-n^{0.1})\,.\notag\\
        & \mathbb P\left[  \sum_{k=1}^M (1+\sigma^2)^{-\frac{k}{2}}\mathsf U_k\mathsf U_{k+M} \geq cM/2 \right]\leq \exp(-\Theta(M))\notag\,,
    \end{align}
    where the first line follows from $M=n^{o(1)}$ by \eqref{eq-def-M-d-large}, and the second line uses the fact that $\{\mathsf U_k \mathsf U_{k+M}\}$ are iid sub-exponential random variables. 
    Plugging the above tail probability into \eqref{eq-each-beta-H0-d-large-calculation-1}, we get that 
    \begin{equation}
        \mathbb P_0\left[ \SMT(\beta) \geq cM \right] \leq \exp(-\Theta(M))+n^{-0.19}+\exp(-n^{0.1})=\exp(-\Theta(M))\,,\notag
    \end{equation}
    where the last equality holds due to $M=o(\log n)$. We have proved \eqref{lem:eachbeta-analysis-H0-d-large}. Applying a union bound on $N_\delta$, we get
    \begin{equation*}
            \mathbb P_0\left[ \SMT\geq cM\right]
            \leq  \sum_{\beta\in N_\delta}\mathbb P\left[ \SMT(\beta)\geq cM\right]\leq |N_\delta| \exp(-\Theta(M))
            \leq  \exp(-\Theta(M)+d \log M)=o(1)\,,
    \end{equation*}
    where the last inequality holds due to $|N_\delta| \le (3/\delta)^d$, $\delta=1/\sqrt{M},$ and $d \log M \le c'M$ for some sufficiently small constant $c'>0.$

Combining the above with~\prettyref{eq:SMT_type_1_error}, we deduce that for any constant $0<c<c_*,$
 $$
 \mathbb{P}_0\left[ \SMT\geq cM\right] + \mathbb{P}_1\left[\SMT \leq c M\right] = o(1),
 $$
 and hence $\SMT$ achieves strong detection.

\subsection{Proof of Theorem~\ref{thm:ASWT}}
\label{sec:empirical-testing-statistic}

We bound the computational complexity of $\ASWT$. For each $\beta \in N_\delta$, computing $\langle X_i,\beta\rangle$ for all $i$ takes $O(nd)$ time, and determining which interval $I_\ell$ contains each sample point takes $n^{1+o(1)}$ time. 
Moreover, $|N_\delta| \le (3/\delta)^d = (3L^2)^d$. Therefore, the total running time is $n^{1+o(1)}\,d\,(3L^2)^d$. Since $\omega_n\le\log n$, this is $n^{1+o(1)}d^{\Theta(d)}$, and it is $n^{1+o(1)}$ when $d=o(\log n/\log\log n)$.

\paragraph{Analysis under $\mathcal{H}_1$}
We first show that for any constant $0<c<1,$
    \begin{equation}\label{eq:high-d-testing-statistic-H1}
            \mathbb P_1\left[ \ASWT<c\sqrt{\frac{L}{n}} \right]=1-o(1)\,.
    \end{equation}
    Without loss of generality, we condition on $\pi_*=\mathsf{Id}$ and $\beta_*=\beta_0$ for a fixed $\beta_0\in\mathbb S^{d-1}$. 
    By the definition of $N_\delta$, there exists a $\widehat\beta$ such that $\|\beta_0-\widehat\beta\|_2\leq 1/L^2$. It suffices to show 
\begin{equation}
    \label{eq-high-d-testing-statistic-H1-goal1}
    \mathbb E_1\left[\mathsf d(\mu_{X \widehat\beta},\mu_Y)\right] = o\left(\sqrt{\frac{L}{n}}\right) \,.
\end{equation}
The desired \eqref{eq:high-d-testing-statistic-H1} then follows from Markov's inequality. To prove~\eqref{eq-high-d-testing-statistic-H1-goal1}, let $ Z_i=\langle X_i,\widehat\beta\rangle$ and $W_i=Y_i/\sqrt{1+\sigma^2}$. Then  
   $$\left\{\left(Z_i, W_i \right)\right\}_{i=1}^n\overset{\mathrm{iid}}{\sim}N \left( \mathbf{0},
   \begin{bmatrix}
       1 & \rho \\
       \rho & 1\,
   \end{bmatrix}
   \right),$$
   where
   $\rho =\langle \hat\beta,\beta_0\rangle/\sqrt{1+\sigma^2} \ge 1-\sigma^2-1/L^4$.
   Let $\mu_{ Z}$ and $\mu_{ W}$
   denote the empirical distributions of $Z_i$ and $W_i$, respectively. By definition,
   $$
   \mathbb E_1\left[\mathsf d(\mu_{X \widehat\beta},\mu_Y)\right]  = \sum_{\ell \in [L]} \mathbb E\left[ | \mu_{ Z}(I_\ell) - \mu_{ W}(I_\ell) | \right].
   $$
   For each $\ell \in [L]$, we have
    \begin{equation}\label{eq-est-each-ell-mudiff}
    \begin{split}
        &\mathbb E\left[ | \mu_{ Z}(I_\ell) - \mu_{ W}(I_\ell) | \right] \leq \sqrt{\mathbb E\left[ ( \mu_{ Z}(I_\ell) - \mu_{W}(I_\ell) )^2 \right]}\\
        = & \frac{1}{n} \sqrt{\mathbb E\left[ \left( \sum_{i=1}^n\mathbf{1}_{ Z_i\in I_\ell} - \sum_{i=1}^n \mathbf{1}_{W_i\in I_\ell}  \right)^2 \right]} = \frac{1}{\sqrt{n}}\sqrt{\mathbb E\left[( \mathbf{1}_{ Z_i\in I_\ell}- \mathbf{1}_{ W_i\in I_\ell})^2\right]}\\
        = & \frac{1}{\sqrt{n}}\sqrt{\mathbb P[ Z_i\in I_\ell, W_i\notin I_\ell] + \mathbb P[ Z_i\notin I_\ell,  W_i\in I_\ell]}\,.
    \end{split}
    \end{equation}
To proceed, we need the following auxiliary result to bound 
     $P[ Z_i\in I_\ell,  W_i\notin I_\ell]$.
     
\begin{lemma}
\label{lem:within-without-Iell-estimate}
    Fix any interval $I = [a,a+1/L]$ with $-1/2\leq a \leq 1/2$ and assume that $\widehat\sigma\triangleq\sigma\vee (1/L^2) \le 1/(8L)$. For two correlated standard Gaussian variables $ Z, W$ with correlation $\rho\geq 1-2\widehat\sigma^2$, 
    \begin{equation}
        \mathbb P\left[  Z\in I, W\notin I \right] \leq \frac{24}{\sqrt{\pi}}\widehat{\sigma}\,.\notag
    \end{equation}
\end{lemma}

\begin{proof}
     Since $(-Z,-W)$ has the same law as $(Z,W)$, and the event $\{Z\in I, W\notin I\}$ becomes $\{-Z\in -I,-W\notin -I\}$ with $-I=[-a-1/L,-a]$, whose left endpoint $-a-1/L$ is positive when $a<-1/L$, we may assume without loss of generality that $a\ge -1/L$.
    We can split the interval $I$ into two parts $$
    I_L=[a,a+8\hat{\sigma}] , \quad \text{ and } \quad  
    I_R=(a+8\hat{\sigma},a+1/L],$$ 
    so that 
        \begin{align*}
         \mathbb P\left[ Z\in I,  W\notin I \right]  
           = 
          \mathbb P\left[ Z\in I_L,  W\notin I \right]
       + \mathbb P\left[ Z\in I_R,W\notin I \right] \,.
       \end{align*}
    Since the probability density function of $Z$ is bounded by $1/\sqrt{2\pi}$,
    $$
    \mathbb P\left[ Z\in I_L,  W\notin I \right]
    \le \mathbb P\left[ Z\in I_L \right]
        \le \frac{8\widehat \sigma}{\sqrt{2\pi}}.
    $$
It remains to bound 
 $ \mathbb P\left[ Z\in I_R,  W\notin I \right]$. Note that 
 $$
 \mathbb P\left[ Z\in I_R,  W\notin I \right]
= \mathbb P\left[ Z\in I_R,  W<a \right]
+ \mathbb P\left[ Z\in I_R,  W>a+1/L \right].
 $$
 Write $W=\rho Z + \sqrt{1-\rho^2} Z'$ for a standard normal $Z'$ independent of $Z$. We have
 $$
 \mathbb P\left[ Z\in I_R,  W<a \right]
 = \mathbb P\left[ Z\in I_R,  Z'< \frac{a - \rho Z }{\sqrt{1-\rho^2}}\right]
 = \expect{ \indc{Z \in I_R} Q \left( \frac{\rho Z -a }{\sqrt{1-\rho^2}}\right)},
 $$
 where $Q(t)=\prob{Z' \ge t}$. Note that  $Q(t) \le \exp(-t^2/2)$ when $t \ge 0$, and $ \rho z \ge a $ for all $z \in I_R$: if $z\le 0$, then $\rho z\ge z>a$; if $z>0$, then $z\le a+1/L\le 1$ and hence $\rho z\ge z-2\widehat\sigma^2 z> a+8\widehat\sigma-2\widehat\sigma^2\ge a$. 
 Then it follows that 
 $$
 \expect{ \indc{Z \in I_R} Q \left( \frac{\rho Z -a }{\sqrt{1-\rho^2}}\right)}
 \le \expect{ \exp \left( - \frac{(\rho Z -a)^2 }{2(1-\rho^2)}\right) }
 \le \sqrt{1-\rho^2},
 $$
 where the last inequality holds due to 
$
\mathbb{E}\left[\exp\left(-\frac{(aZ - b)^2}{2c}\right)\right]
= \sqrt{\frac{c}{a^2 + c}} \;
\exp\left(-\frac{b^2}{2(a^2 + c)}\right)$ for $c>0.$
Similarly,
\begin{align*}
 \mathbb P\left[ Z\in I_R,  W>a+1/L \right]
 &= \mathbb P\left[ Z\in I_R,  Z'> \frac{a +1/L - \rho Z }{\sqrt{1-\rho^2}}\right] \\
 & = \expect{ \indc{Z \in I_R} Q \left( \frac{a+1/L-\rho Z }{\sqrt{1-\rho^2}}\right)} \\
 & \le \expect{ \exp \left( - \frac{(a+1/L-\rho Z)^2 }{2(1-\rho^2)}\right) }
 \le \sqrt{1-\rho^2},
 \end{align*} 
 where the first inequality uses $Q(t)\le e^{-t^2/2}$ for $t\ge0$, which applies since $\rho z\le a+1/L$ for all $z\in I_R$: if $z\ge 0$, then $\rho z\le z\le a+1/L$; if $z<0$, then $\rho z<0\le a+1/L$ because $a\ge -1/L$.
Combining the above bounds and using $1-\rho^2\le 2(1-\rho)\le 4\widehat\sigma^2$, we obtain $\mathbb P[Z\in I,W\notin I]\le \frac{8\widehat\sigma}{\sqrt{2\pi}}+2\sqrt{1-\rho^2}\le \big(\frac{8}{\sqrt{2\pi}}+4\big)\widehat\sigma\le \frac{24}{\sqrt\pi}\widehat\sigma$, which completes the proof of this lemma.
\end{proof}

Next, we apply Lemma~\ref{lem:within-without-Iell-estimate}. 
Note that $\widehat\sigma \ll L^{-1}$ from the choice of $L$ in Theorem~\ref{thm:ASWT}.
Moreover, 
\begin{equation*}
     \rho  \geq 1-\sigma^2-1/L^4 \geq 1-2\widehat\sigma^2 \,,
 \end{equation*}
 since $\widehat\sigma=\sigma\vee L^{-2}$, so the hypothesis of Lemma~\ref{lem:within-without-Iell-estimate} holds for every bin $I_\ell$.
By Lemma~\ref{lem:within-without-Iell-estimate},
    \begin{equation}
        \label{eq-high-d-testing-statistic-H1-goal2}
        \mathbb P\left[ \mathsf Z_i\in I_\ell, \mathsf W_i\notin I_\ell \right] = \mathbb P[\mathsf Z_i\notin I_\ell, \mathsf W_i\in I_\ell]\leq \frac{24}{\sqrt{\pi}}\widehat\sigma\,.
    \end{equation}
    Plugging \eqref{eq-high-d-testing-statistic-H1-goal2} into \eqref{eq-est-each-ell-mudiff} yields
    \begin{equation}
        \mathbb E\left[ | \mu_{\mathsf Z}(I_\ell) - \mu_{\mathsf W}(I_\ell) | \right] \leq \frac{4\sqrt{3}}{\pi^{1/4}}\sqrt{\frac{\widehat\sigma}{n}}\,.\notag
    \end{equation}
    Recall that $\widehat\sigma \ll L^{-1}$.
    Summing over $\ell \in [L]$ yields \eqref{eq-high-d-testing-statistic-H1-goal1}.

\paragraph{Analysis under $\mathcal{H}_0$} 
We show that for $\bar{\mathsf c}\triangleq\frac{e^{-1/16}}{8(2\pi)^{1/4}}$ and for each $\beta\in N_\delta$, 
\begin{equation}\label{eq:high-d-testing-statistic-H0}
        \mathbb P_0\left[ \mathsf{d}(\mu_{X\beta},\mu_{Y}) \leq \bar{\mathsf c}\sqrt{\frac{L}{n}} \right]\leq\exp\left(-\frac{\bar{\mathsf c}^2}{4}L \right)\,.
\end{equation}
Let $Z$ and $W$ denote independent standard Gaussian vectors and $\mu_Z=\frac{1}{n}\sum_{i=1}^n \delta_{Z_i}$ (similarly for $\mu_W$).
Let
$
    f(z_1,\dots,z_n,w_1,\dots,w_n) \triangleq \mathsf{d}(\mu_z,\mu_w)\,.\notag
$
There exist $\ell,\ell'\in[L]$ such that we have 
\begin{equation}
\begin{split}
   &| f(z_1,\dots,z_i,\dots,z_n,w_1,\dots,w_n) - f(z_1,\dots,z'_i,\dots,z_n,w_1,\dots,w_n) |\\
   \leq & \left|\mu_z(I_\ell)-\mu_w(I_\ell)+1/n \right|+ \left|\mu_z(I_{\ell'})-\mu_w(I_{\ell'})-1/n\right| - \left|\mu_z(I_\ell)-\mu_w(I_\ell)\right|-\left|\mu_z(I_{\ell'})-\mu_w(I_{\ell'})\right|\\
   \leq & \frac{2}{n}\,.
\end{split}\notag
\end{equation}
By McDiarmid's inequality, we obtain
\begin{equation}
    \mathbb P\left[ \mathsf{d}(\mu_{Z},\mu_{W}) \leq \mathbb E\mathsf{d}(\mu_{Z},\mu_{W})-\sqrt{\frac{\Delta}{n}} \right]\leq \exp(-\frac{2\Delta/n}{2n (2/n)^2})=\exp(-\Delta/4)\,. \label{eq-high-d-testing-statistic-H0-goal-1}
\end{equation}
 Then we show 
\begin{equation}\label{eq-high-d-testing-statistic-H0-goal-2}
    \mathbb E\left[ \mathsf{d}(\mu_Z,\mu_W) \right] \geq 2\bar{\mathsf c}\sqrt{\frac{L}{n}}\,.
\end{equation}
For any $\ell\in[L]$, 
\begin{equation}\label{eq-high-d-muz-muw-each-l-H0}
\begin{split}
    & \mathbb E\left[ |\mu_Z(I_\ell)-\mu_W(I_\ell)| \right] = \mathbb E\left[\left|  \frac{1}{n}\sum_{i=1}^n\mathbf{1}_{Z_i\in I_\ell}-\frac{1}{n}\sum_{i=1}^n \mathbf{1}_{W_i\in I_\ell}  \right|  \right] \geq \inf_{x\in \mathbb R}\mathbb E\left[\left|  \frac{1}{n}\sum_{i=1}^n\mathbf{1}_{Z_i\in I_\ell}-x  \right|  \right]\\
    \geq & \frac{1}{n}\mathbb E\left[\left|  \sum_{i=1}^n\mathbf{1}_{Z_i\in I_\ell}-x_0  \right|  \right]\,,
\end{split}
\end{equation}
where $x_0$ is the median of $\sum_{i=1}^n\mathbf{1}_{Z_i\in I_\ell}$,
which is distributed as $\Binom(n,\alpha)$ with $\alpha=\prob{Z_i \in I_\ell}.$
Note that $\frac{e^{-\frac{1}{8}}}{\sqrt{2\pi}}\frac{1}{L} \leq \alpha \leq \frac{1}{\sqrt{2\pi}}\frac{1}{L}$. By \cite{kaas1980mean}, we have $|x_0-n\alpha|\leq 1$.
Using the estimate for the mean absolute deviation of a binomial distribution (see \cite[Theorem~1]{berend2013sharp}), we have
\begin{equation}
    \mathbb E\left[ | \Binom(n,\alpha)-n\alpha | \right]\geq \frac{\sqrt{n\alpha(1-\alpha)}}{\sqrt{2}}
   \geq \frac{1}{2}\sqrt{n\alpha} \geq \frac{e^{-1/16}}{2(2\pi)^{1/4}}\sqrt{\frac{n}{L}}\,.\label{eq:median-est-fluctuation}
\end{equation}
In addition, recalling the conditions on $d$ and $L$ in Theorem~\ref{thm:ASWT}, we have 
$$L\leq \omega_n^{1/3}d\log (d+1) \leq \omega_n^{1/3}\frac{n}{\omega_n\log n}\log n = o(n),$$
which implies $n/L\to\infty$. Therefore, plugging \eqref{eq:median-est-fluctuation} into \eqref{eq-high-d-muz-muw-each-l-H0} yields
$$
    \mathbb E\left[ |\mu_Z(I_\ell)-\mu_W(I_\ell)| \right]\geq \frac{2\bar{\mathsf c
    }}{\sqrt{nL}}
$$
with $\bar{\mathsf c}$ defined above.
Summing over all $\ell\in[L]$, we obtain \eqref{eq-high-d-testing-statistic-H0-goal-2}. Combining this with 
\eqref{eq-high-d-testing-statistic-H0-goal-1} and choosing $\Delta=\bar{\mathsf c}^2 L $, we obtain~\eqref{eq:high-d-testing-statistic-H0}.  Applying a union bound on $N_\delta$, we get that 
\begin{align*}
      \mathbb P_0\left[\ASWT\leq \bar{\mathsf c}\sqrt{\frac{L}{n}}\right]
     \leq  |N_\delta|\exp\left( -\frac{\bar{\mathsf c}^2}{4} L\right)\leq \exp\left( -\frac{\bar{\mathsf c}^2}{4} L+ 3d\log L\right)=o(1)\,,
\end{align*}
where the last inequality holds because $|N_\delta|\leq (3/\delta)^d$ and $\delta=1/L^2$; the last inequality follows from our choice of $L$.

Combining the above with \eqref{eq:high-d-testing-statistic-H1}, we conclude that for any constant $0<c<\bar{\mathsf c}$,
\begin{equation*}
    \mathbb P_0\left[\ASWT\leq c\sqrt{\frac{L}{n}}\right] + \mathbb P_1\left[\ASWT\geq c\sqrt{\frac{L}{n}}\right]=o(1)\,,
\end{equation*}
and thus $\ASWT$ achieves strong detection.

\subsection{Proof of Theorem~\ref{thm:lower-bound}}
\label{sec:lower-bound}

The proof for Theorem~\ref{thm:lower-bound} is inspired by \cite{jiao2025DetectionBroken} which  in turn builds on  the general strategy of~\cite{bai2005broken}. A key difference is that, unlike the broken sample problem in \cite{jiao2025DetectionBroken} where $\beta^*$ is known, our model involves two latent parameters, $\beta_*$ and $\pi_*$.
As a result, when computing the second moment of the likelihood ratio, $\mathbb{E}_0[L^2(X,Y)]$, we encounter two independent replicas $(\beta_1,\beta_2)$ and $(\pi_1,\pi_2)$, whereas \cite{jiao2025DetectionBroken} only involves the latter. 

To this end, we condition on $(\beta_1,\beta_2)$ and first take expectation over $(\pi_1,\pi_2)$, using cycle decomposition and P\'olya's theorem \cite[Theorem~7.7]{stanley2013algebraic}, analogous to the approach in \cite{jiao2025DetectionBroken}. 
This allows us to express the second moment $\mathbb{E}_0[L^2(X,Y)]$ as $\mathbb E_{\beta_1,\beta_2}\left[g(\beta_1,\beta_2)\right]$, where $g(\beta_1,\beta_2)$ is given by the $n$-th coefficients of the power series of $f(z)\triangleq \prod_{k=0}^{\infty}\frac{1}{1-z\lambda_k}$ with $\lambda_k=\iprod{\beta_1}{\beta_2}^k/(1+\sigma^2)^{k}$ and can be extracted via a contour integral of $f$ (see Lemma~\ref{lem:lower-bound-second-moment-before-integral} for details).

It remains to average over $(\beta_1,\beta_2)$. In the broken sample problem \cite{jiao2025DetectionBroken}, one effectively has $\beta_1=\beta_2$, so no such averaging is needed. Moreover, there it suffices to focus on the large-noise regime $\sigma^2=\Omega(1)$, in which case $|\lambda_k|$'s are uniformly bounded away from $1.$ In contrast, in the growing-dimension regime of our problem, we must handle noise as small as $\sigma^2 \gtrsim \frac{1}{d\log d}$. In this case,  $|\lambda_k|$'s can be close to $1$
when $|\iprod{\beta_1}{\beta_2}|$ is close to $1$, necessitating a tight control of the exponential growth $g(\beta_1,\beta_2)$ and a careful comparison with the density of $|\iprod{\beta_1}{\beta_2}|$ near $1$ (see~Lemma~\ref{lem:lower-bound-est-expectation} for details).

We begin the proof by rewriting the likelihood ratio $L(X,Y)$ as 
\begin{equation}\label{eq:def-likelihood-and-Lbeta}
    L(X,Y) = \frac{\mathbb P_1[X,Y]}{\mathbb P_0[X,Y]} = \mathbb E_{\beta} \underbrace{\mathbb E_{\Pi}\Big[ \frac{\mathbb P_1[X,Y\,|\,\beta,\Pi]}{\mathbb P_0[X,Y]} \Big]}_{\triangleq L_\beta(X,Y)}\,.
\end{equation}
Therefore, 
\begin{equation}\label{eq-second-moment-integral-over-betas}
    \mathbb E_0[L^2(X,Y)]  = \mathbb E_0\Big[ \mathbb E_{\beta_1}[L_{\beta_1}(X,Y)]\mathbb E_{\beta_2}[L_{\beta_2}(X,Y)] \Big]
    =\mathbb E_{\beta_1,\beta_2} \Big[ \underbrace{\mathbb E_0[L_{\beta_1}(X,Y)L_{\beta_2}(X,Y)]}_{\triangleq g(\beta_1,\beta_2)} \Big]\,.
\end{equation}
The following lemma provides an upper bound for $g(\beta_1,\beta_2)$ for a given pair $\beta_1,\beta_2$.

\begin{lemma}
    \label{lem:lower-bound-second-moment-before-integral}
    Denote $\lambda \triangleq \frac{|\langle\beta_1,\beta_2\rangle|}{1+\sigma^2}$. Then 
    \begin{equation}
   g(\beta_1,\beta_2)  \le \left(1+\left(\frac{1+\lambda}{2}\right)^n\left(\frac{1+\lambda}{1-\lambda}\right)^2\right) 
        \exp \left( \frac{1}{1-\lambda} \min \left\{ \frac{\lambda}{1-\lambda}, \frac{\pi^2}{6} \right\}\right) 
    \,. \label{eq:g_beta_bound}
    \end{equation}
\end{lemma}

\begin{proof}
Letting 
$\rho=1/\sqrt{1+\sigma^2}$ and $\overline{Y}_i= \rho Y_i,$ we have that for  each fixed $\Pi$ and $\beta$, 
\begin{equation}
    \frac{\mathbb P_1[X,Y\,|\,\beta,\Pi]}{\mathbb P_0[X,Y]}=
    \prod_{i=1}^n \ell\left(\iprod{X_{\pi(i)}}{\beta}, \overline{Y}_i\right),
\end{equation}
where $\ell(x,y)$ is known as Mehler's kernel and admits a spectral decomposition via Hermite polynomials:
$$
\ell(x,y) \triangleq \frac{1}{\sqrt{1-\rho^2}} \exp\left( 
-\frac{(y-\rho x)^2}{2(1-\rho^2)} + \frac{y^2}{2} \right) =  \sum_{k=0}^\infty \rho^k h_k(x) h_k(y).
$$
It follows that 
\begin{equation}
\begin{split}\label{eq-second-moment-expression-1}
    g(\beta_1,\beta_2) = & \mathbb E_{X,Y}\left[\mathbb E_{\Pi_1}\Big[\prod_{i=1}^n \ell\left( \iprod{X_{\pi_1(i)}}{\beta_1}, \overline{Y}_i \right)\Big]\times\mathbb E_{\Pi_2}\Big[\prod_{i=1}^n \ell \left( \iprod{X_{\pi_2(i)}}{\beta_2},\overline{Y}_i \right)\Big]\right]\\
    =&\mathbb E_{\Pi}\mathbb E_{X,Y}\left[\prod_{i=1}^n\ell \left( \iprod{X_i}{\beta_1},\overline{Y}_i \right)\ell\left( \iprod{X_{\pi(i)}}{\beta_2},\overline{Y}_i \right)\right]\\
    =&\mathbb E_{\Pi}\left[\prod_{\mathcal O\in\mathfrak{O}(\pi)}\mathbb E_{X,Y}\Big[\prod_{i\in\mathcal{O}} \ell\left( \iprod{X_i}{\beta_1},\overline{Y}_i \right)\ell\left( \iprod{X_{\pi(i)}}{\beta_2},\overline{Y}_i \right)\Big]\right]\,,\\
\end{split}
\end{equation}
where $\mathfrak O(\pi)$ is the set of cycles in the cycle decomposition of permutation $\pi$ and $\mathcal O$ represents a cycle.
For each cycle in the decomposition, we have
\begin{equation}
\begin{split}\label{eq-second-moment-expression-2}
    &\mathbb E_{X,Y}\left[\prod_{i\in\mathcal{O}} \ell\left( \iprod{X_i}{\beta_1},\overline{Y}_i \right)\ell\left( \iprod{X_{\pi(i)}}{\beta_2},\overline{Y}_i \right)\right]\\
    =&\mathbb E_{X,Y} \left[\prod_{i\in \mathcal{O}}\left(\sum_{k,k'\geq 0}\rho^{k+k'} h_k\left(\langle X_{i},\beta_1\rangle\right)h_k\left(\overline{Y}_i\right)h_{k'}\left(\langle X_{\pi(i)},\beta_2\rangle\right)h_{k'}\left(\overline{Y}_i\right)\right)\right]\\
    =&\mathbb E_{X}\left[\prod_{i\in \mathcal{O}}\left(\sum_{k=0}^\infty\rho^{2k} h_k\left(\langle X_{i},\beta_1\rangle\right)h_k\left(\langle X_{\pi(i)},\beta_2\rangle\right)\right)\right]=
    \sum_{k=0}^\infty \left(\rho^2 \langle \beta_1,\beta_2\rangle\right)^{k|\mathcal{O}|}\,,
\end{split}
\end{equation}
where the second equality follows from taking expectation with respect to $Y$ and the third equality follows from expanding the summation within the parentheses. Plugging \eqref{eq-second-moment-expression-2} into \eqref{eq-second-moment-expression-1}, we have
\begin{equation}
\begin{split}
    g(\beta_1,\beta_2)=\mathbb E_{\Pi} \left[\prod_{\mathcal{O}\in\mathfrak{O}(\pi)}\left(\sum_{k=0}^\infty \left(\rho^2\langle \beta_1,\beta_2\rangle \right)^{k|\mathcal{O}|}\right)\right]
    =  \mathbb E_{\Pi}\left[\prod_{i=1}^n\left(\sum_{k=0}^\infty \left(\rho^2 \langle \beta_1,\beta_2\rangle\right)^{ki}\right)^{N_i(\Pi)}\right]\,,
\end{split}\notag
\end{equation}
where $N_i(\pi)$ is the number of cycles of length $i$ in the cycle decomposition of $\pi$. 
Denoting $\lambda_k=\rho^{2k} \langle \beta_1,\beta_2\rangle^k$, 
by P\'olya's theorem in \cite[Theorem~7.7]{stanley2013algebraic}, $
    g(\beta_1,\beta_2)$ equals the $n$-th coefficients of the power series of $f(z)\triangleq \prod_{k=0}^{\infty}\frac{1}{1-z\lambda_k}$.
We can extract the $n$-th coefficient via the contour integral. Specifically, choose
$r=\frac{2}{1+|\lambda_1|} \in (1, 1/|\lambda_1|)$. 
Applying the residue theorem yields
\begin{equation}\label{eq-residue-application}
\begin{split}
    I\triangleq \frac{1}{2\pi i}\oint_{|z|=r}\frac{1}{z^{n+1}}f(z)dz 
    =\frac{(f(z))^{(n)}}{n!}\bigg|_{z=0}-\frac{(1-z)f(z)}{z^{n+1}}\bigg|_{z=1}=g(\beta_1,\beta_2)-\prod_{k=1}^\infty\frac{1}{1-\lambda_k}\,,
\end{split}\notag
\end{equation}
where the second equality holds because $f(z)/z^{n+1}$
has only two
singularities, $0$ and $1$, inside the contour; the other singularities, $\{1/\lambda_k, k \ge 1\}$, lie strictly outside the contour. Furthermore, 
\begin{align*}
    |I| & \leq \frac{1}{2\pi i}\oint_{|z|=r}\frac{1}{| z^{n+1} |} |f(z)| dz  \leq  
    \frac{1}{r^n(r-1)}\prod_{k=1}^\infty \frac{1}{1-r|\lambda_k|}  \\
    & =\left(\frac{1+\lambda}{2}\right)^n\left(\frac{1+\lambda}{1-\lambda}\right)^2  \cdot \prod_{k = 2}^{\infty}\frac{1}{1-\frac{2\lambda}{1+\lambda}\lambda^{k-1}} 
\leq \left(\frac{1+\lambda}{2}\right)^n\left(\frac{1+\lambda}{1-\lambda}\right)^2 \cdot \prod_{k=1}^\infty\frac{1}{1-\lambda^{k}}\,,
        \label{eq-lower-bound-control-I}
\end{align*}
where the equality follows by plugging $r=2/(1+|\lambda_1|)$ and noting that $|\lambda_k|=\lambda^k.$  Combining the last two displayed equations yields that 
$$
g(\beta_1, \beta_2) \le \left(1+\left(\frac{1+\lambda}{2}\right)^n\left(\frac{1+\lambda}{1-\lambda}\right)^2\right)\prod_{k=1}^\infty \frac{1}{1-\lambda^k}.
$$
It remains to bound from above $\prod_{k=1}^\infty\frac{1}{1-\lambda^k}$. Note that 
\begin{equation}\label{eq-leq-iota_n-est}
    \begin{split}
       \log  \prod_{k=1}^\infty \frac{1}{1-\lambda^k} = -\sum_{k=1}^\infty\log(1-\lambda^k)\leq \sum_{k=1}^\infty\frac{\lambda^k}{1-\lambda^k} 
        \leq  \frac{\sum_{k=1}^\infty\lambda^k}{1-\lambda}= \frac{\lambda}{(1-\lambda)^2} .
    \end{split}
\end{equation}

Alternatively, we have
\[
\log \prod_{k=1}^{\infty} \frac{1}{1-\lambda^k}
=
-\sum_{k=1}^\infty\log(1-\lambda^k)
=\sum_{k=1}^{\infty} \sum_{m=1}^{\infty} \frac{ \lambda^{mk}}{m} 
\le
\sum_{m=1}^{\infty} \frac{1}{m^2 \log (1/\lambda) } 
\le \frac{\pi^2}{6(1-\lambda)},
\]
where the second equality holds due to 
$-\log(1-x) = \sum_{m=1}^{\infty} \frac{x^m}{m}$ for $ |x|<1$;
the first inequality holds from
$\sum_{k=1}^\infty \lambda^{mk} \le \int_0^\infty \lambda^{xm} dx
=\frac{1}{m\log(1/\lambda)}$; the last inequality holds due to
$\log\lambda \le \lambda-1.$
Taking the minimum of the above two upper bounds completes the proof. 
\end{proof}
The next lemma is crucial for bounding the truncated mean of $g(\beta_1,\beta_2)$ in growing-dimension regime. 
\begin{lemma}
    \label{lem:lower-bound-est-expectation}
    Denote $\lambda \triangleq \frac{|\langle\beta_1,\beta_2\rangle|}{1+\sigma^2}$. Assume that $\sigma^2\geq \frac{(1+\varepsilon)\pi^2}{3d\log d}$ for some constant $\epsilon>0.$ There exists a constant $c>0$ such that for all sufficiently large $d,$
    \begin{equation}
        \mathbb E_{\beta_1,\beta_2}\left[\exp \left( 
        \frac{\pi^2}{6(1-\lambda)}\right) \mathbf{1}_{|\langle\beta_1,\beta_2\rangle|>d^{-1/4}}\right] \le \exp 
        \left(-c \sqrt{d} \right) \,.\label{eq:g_truncated_mean}
    \end{equation}
    where $\beta_1,\beta_2$ are independently and  uniformly drawn from the unit sphere $\mathbb S^{d-1}$.
\end{lemma}
\begin{proof}

Since $\langle\beta_1,\beta_2\rangle$ is the sub-Gaussian with sub-Gaussian norm $1/\sqrt{d},$ it follows that 
    \begin{equation}\label{eq-tail-angle-of-betas}
        \mathbb P \Big[|\langle\beta_1,\beta_2\rangle|\geq d^{-1/4}\Big] \leq \exp(-c'\sqrt{d}),
    \end{equation}
    where $c'>0 $ is a constant. 
If $\sigma^2\geq 1/4$, then $\lambda \le 1/(1+\sigma^2) \le 4/5$ and thus
\begin{align*}
 \mathbb E_{\beta_1,\beta_2}\left[\exp \left( 
        \frac{\pi^2}{6(1-\lambda)}\right) \mathbf{1}_{|\langle\beta_1,\beta_2\rangle|>d^{-1/4}}\right] & \le \exp \left( \frac{5\pi^2}{6}\right)
        \prob{|\langle\beta_1,\beta_2\rangle|>d^{-1/4} } \\
        & \le \exp \left( \frac{5\pi^2}{6} -c'\sqrt{d} \right)
        \le \exp\left(-c\sqrt{d}\right).
\end{align*}
Henceforth, we focus on $\sigma^2\leq 1/4$. Similarly as above, if
$|\langle\beta_1,\beta_2\rangle| \le 4/5,$ then  
$\lambda =|\langle\beta_1,\beta_2\rangle|/(1+\sigma^2) \le 4/5$
and thus
\begin{align*}
 \mathbb E_{\beta_1,\beta_2}\left[\exp \left( 
        \frac{\pi^2}{6(1-\lambda)}\right) \mathbf{1}_{d^{-1/4} \leq |\langle\beta_1,\beta_2\rangle|\leq 4/5}\right]  \le \exp \left( \frac{5\pi^2}{6}\right)
        \prob{|\langle\beta_1,\beta_2\rangle|>d^{-1/4} } \le \exp\left(-c\sqrt{d}\right).
\end{align*}
Therefore, it remains  to show 
\begin{equation}\label{eq:final-goal-lower-bound}
    \mathbb E_{\beta_1,\beta_2}\left[\exp \left( 
        \frac{\pi^2}{6(1-\lambda)}\right) \mathbf{1}_{|\langle\beta_1,\beta_2\rangle|>4/5}\right] \le \exp\left(-c \sqrt{d}\right)\,.
\end{equation}
Recall that the probability density function of $\langle \beta_1,\beta_2 \rangle$ is given by 
\begin{equation}
    f(t)\triangleq c_d (1-t^2)^{D}\mathbf{1}_{[-1,1]}, \quad \text{ where } c_d = \frac{1}{B\left(\frac{d-1}{2},\frac{1}{2}\right)} \text{ and } D=\frac{d-3}{2} \,.\notag
\end{equation}
Therefore, 
\begin{align*}
\mathbb E_{\beta_1,\beta_2}\left[\exp \left( 
        \frac{\pi^2}{6(1-\lambda)}\right) \mathbf{1}_{|\langle\beta_1,\beta_2\rangle|>4/5}\right]
        & = 2\int_{4/5}^{1}
        \exp \left( 
        \frac{\pi^2}{6(1-t/(1+\sigma^2))}\right)
        f(t) \diff t \\
        & = 2\int_{0}^{1/5}
        \exp \left( 
        \frac{c_\sigma }{\sigma^2+s}\right)
        f(1-s) \diff s,
\end{align*}
where $c_\sigma= \pi^2 (1+\sigma^2)/6.$
Note that $f(1-s) \le c_d (2s)^{D}$. Let $u = 1/(\sigma^2+s)$. Then $s=1/u-\sigma^2 \le 1/u$ and $\diff s = - u^{-2}\diff u$. It follows that 
\begin{align*}
\int_{0}^{1/5}
        \exp \left( 
        \frac{c_\sigma}{\sigma^2+s}\right)
        f(1-s) \diff s
&\leq c_d \int_{1/(\sigma^2+1/5)}^{1/\sigma^2}
\exp \left( c_\sigma u \right) (2/u)^{D} u^{-2} \diff u \\
& \le c_d  \exp \left( \max_{20/9 \le u \le 1/\sigma^2}  h(u) \right)\int_{20/9}^{1/\sigma^2}  u^{-2} \diff u \\
& \le c_d \exp \left( \max\{h(20/9), h(1/\sigma^2)\} \right) \cdot (9/20),
\end{align*}
where the second inequality holds due to $1/(\sigma^2+1/5) \ge 20/9$ as $\sigma^2 \le 1/4$ and $h(u) \triangleq c_\sigma u + D \log (2/u)$; the last inequality holds because $h(u)$ is convex and thus has its maximum at the endpoints. Note that 
\begin{align*}
h(20/9) &= (20/9) c_\sigma + D \log (9/10) \le - c' d, \\
h(1/\sigma^2) & = c_\sigma /\sigma^2 + D \log (2\sigma^2) \le -c'd,
\end{align*}
where $c'>0$ is some constant and inequalities hold for all sufficiently large $d$, due to the assumption that $\frac{(1+\epsilon)\pi^2}{3d\log d} \leq \sigma^2 \leq 1/4. $
For the second inequality, note that $s\mapsto \pi^2(1+s)/(6s)+D\log(2s)$ is decreasing on $(0,\pi^2/(6D)]$ and increasing afterwards, so it suffices to check the endpoints $\sigma^2=\frac{(1+\epsilon)\pi^2}{3d\log d}$ and $\sigma^2=1/4$.
Combining the last three displayed equations and noting that $c_d =O(\sqrt{d})$ yields the desired bound~\eqref{eq:final-goal-lower-bound}.

\end{proof}

Now we are ready to finish the proof of Theorem~\ref{thm:lower-bound}.
\begin{proof}[Proof of Theorem~\ref{thm:lower-bound}]
Recalling \eqref{eq-second-moment-integral-over-betas} and applying Lemma~\ref{lem:lower-bound-second-moment-before-integral}, we have the following:
\begin{itemize}
    \item For $d=O(1)$, for any constant $\epsilon$, if $\sigma^2\geq \epsilon$, then $\lambda \le 1/(1+\sigma^2) \le 1/(1+\epsilon)$. Therefore, $\left( \frac{1+\lambda}{2} \right)^n\left(\frac{1+\lambda}{1-\lambda}\right)^2\to 0$ and hence by~\prettyref{eq:g_beta_bound}, $
    \mathbb E_0L^2=O(1)$.  If $\sigma^2\to\infty$, then $\lambda\leq 1/(1+\sigma^2)=o(1)$. Again, by~\prettyref{eq:g_beta_bound},  $\mathbb E_0[L^2]=1+o(1)$.

\item For $ \omega(1) \leq d\leq \frac{n}{(\log n)^3}$ and  $\sigma^2 \geq  \frac{(1+\varepsilon)\pi^2}{3d\log d}$, we have
\begin{equation}
    \left( \frac{1+\lambda}{2} \right)^n\left(\frac{1+\lambda}{1-\lambda}\right)^2\leq \exp\left(-\frac{\sigma^2n}{2(1+\sigma^2)}\right)\left(\frac{2(1+\sigma^2)}{\sigma^2}\right)^2 \le o(1)\,.\label{eq:lambda-factors}
\end{equation}
Therefore, applying  Lemma~\ref{lem:lower-bound-second-moment-before-integral} and Lemma~\ref{lem:lower-bound-est-expectation}, we obtain 
\begin{align*}
    \mathbb E_0L^2 \leq  &\left(1+o(1) \right)
   \mathbb E_{\beta_1,\beta_2}\left[\exp \left( 
        \frac{\lambda}{(1-\lambda)^2}\right) \mathbf{1}_{|\langle\beta_1,\beta_2\rangle|\le d^{-1/4}}\right]  \\
         & +  (1+o(1)) \mathbb E_{\beta_1,\beta_2}\left[\exp \left( 
        \frac{\pi^2}{6(1-\lambda)}\right) \mathbf{1}_{|\langle\beta_1,\beta_2\rangle|>d^{-1/4}}\right]  \\
    & \le \left(1+o(1) \right)   \exp \left(O(d^{-1/4}) \right) +o(1) 
    =1+o(1)\,.\notag
\end{align*}
\end{itemize}
\end{proof}

\subsection{Proof of Theorem~\ref{thm:beta-recovery}}
\label{subsec-beta-recovery}
   Now we construct an estimator $\tilde{\beta}$ such that $\|\tilde{\beta}-\beta_*\|_2\to0$. Recall the definition of $\SMT(\beta)$ above \eqref{eq-final-goal-H1-d-large}. Motivated by \emph{sliced-moment test}, we define 
    \begin{equation}\label{eq-def-tilde-beta}
        \tilde{\beta}\triangleq \argmax_{\beta\in N_\delta} \SMT(\beta)\,,
    \end{equation}
    where the choice of parameters is the same as Theorem~\ref{thm:SMT}. We continue to use $\widehat\beta$ to denote the closest element to $\beta$   in $N_\delta$. 
   Recall in \eqref{eq-final-goal-H1-d-large}, we have proved
   \begin{equation}
   \forall \beta,\,\,
        \mathbb P_{1,\beta}\left[ \SMT(\widehat\beta) \leq \mathsf c_*M \right]\leq o(1)\,,
    \end{equation}
    where $\mathsf c_*=e^{-1}-e^{-\frac{3}{2}}$.

  It remains to show that with probability $1-o(1)$,
    \begin{equation}\label{eq-goal-beta-recovery}
       \SMT(\beta)<\mathsf c_* M, \quad \forall \beta \in N_\delta, \,\|\beta-\beta_*\|_2> \frac{10}{\sqrt{M}},
   \end{equation}
 which immediately implies that $\|\tilde\beta - \beta_*\|_2\leq 10/\sqrt{M}$.

   We first consider the case $d=1$. Then $N_\delta=\{\pm1\}$, so~\eqref{eq-goal-beta-recovery} only concerns $\beta=-\beta_*$. The derivation of~\eqref{eq-expansion-second-moment-H1-d-large} applies to any fixed direction in place of $\widehat\beta$, so $\operatorname{Var}_1[\SMT(-\beta_*)]\le 2M+\exp(O(M))/n$, while $\mathbb E_1[\SMT(-\beta_*)]=\sum_{k=1}^M\big(-\frac{1}{1+\sigma^2}\big)^k\le 0$. By Chebyshev's inequality, $\SMT(-\beta_*)<\mathsf c_*M$ with probability $1-o(1)$, which proves~\eqref{eq-goal-beta-recovery} for $d=1$.
   
  For $d\geq2$, to prove~\eqref{eq-goal-beta-recovery}, we first approximate $S_k(\beta)$ and $T_k$ by Gaussian random variables, using Lemma~\ref{lem:gaussian-approximation-hermite-w/correlation}. However, since the upper bound in Lemma~\ref{lem:gaussian-approximation-hermite-w/correlation} contains the factor $(1-|\rho|)^{-1}$, we need to ensure that $\beta$ is not very close to $\pm \beta^*$ for all $\beta \in N_\delta:$
   \begin{equation}
       \label{eq-est-distance-to-betastar}
       \begin{split}
      & \mathbb P\left[ \|\beta-\beta_*\|_2\geq n^{-0.01}, \forall \beta\in N_\delta\right] = 1-o(1)\,,\\
      & \mathbb P\left[ \|\beta+\beta_*\|_2\geq n^{-0.01}, \forall \beta\in N_\delta\right] = 1-o(1)\,.
      \end{split}
   \end{equation}
   We prove the first line; the second line follows in the same way. Let  $\angle(\beta,\beta_*)$ denote the angle between $\beta,\beta_*$. If $\|\beta-\beta_*\|_2\leq n^{-0.01}$, then $\angle(\beta,\beta_*)\leq 10n^{-0.01}$. 
   Since $\beta_*$ is uniformly distributed on the unit sphere in $\mathbb R^d$, for each $\beta\in N_\delta$, recalling that $d=o\left(\frac{\log n}{\log\log n}\right)$, we have
   \begin{equation}
   \begin{split}
       & \mathbb P[\|\beta-\beta_*\|_2\leq n^{-0.01}]\leq \mathbb P[ \angle(\beta,\beta_*)\leq 10n^{-0.01} ]\\
       = & \frac{\Gamma(\frac{d}{2})}{\sqrt{\pi}\Gamma(\frac{d-1}{2})}\int_0^{10n^{-0.01}} \sin^{d-2}(\varphi)\mathrm{d}\varphi \leq  n^{-0.005}\,.
    \end{split}
   \end{equation}
   Also note that $|N_\delta|=n^{o(1)}$. Applying a union bound yields \eqref{eq-est-distance-to-betastar}. 
   
   Now fix $\beta_*$ satisfying \eqref{eq-est-distance-to-betastar}. Since $\|\beta\pm\beta_*\|_2^2=2\pm 2\langle \beta,\beta_*\rangle$, we have $|\langle \beta,\beta_*\rangle| \le 1-\frac12 n^{-0.02}$ for every $\beta\in N_\delta$.
   For ease of notation, for any $\beta\in N_\delta$, we let $T_{M+k} = (1+\sigma^2)^{\frac{k}{2}}S_k(\beta)$ and $\{(\mathsf U_k,\mathsf U_{M+k})\}_{1\le k\le M}$ be independent pairs of correlated standard Gaussian variables with correlation $\left(\frac{\langle\beta,\beta_*\rangle}{\sqrt{1+\sigma^2}}\right)^k$, by Lemma~\ref{lem:gaussian-approximation-hermite-w/correlation}, we get that 
   \begin{equation}
       \mathcal{W}_1\left\{[T_k]_{k=1}^{2M},[\mathsf U_k]_{k=1}^{2M}\right\}\leq n^{-0.5+0.02+o(1)}=n^{-0.48+o(1)}\,,
   \end{equation}
   where $\mathcal{W}_1$ denotes the Wasserstein-1 distance. Thus, there exists a coupling of $\mathsf T = [T_k]_{k=1}^{2M}$ and $\mathsf U=[\mathsf U_k]_{k=1}^{2M}$ such that
   \[
   \mathbb E\{\|\mathsf T-\mathsf U\|_2\}\leq n^{-0.47}.
   \]
   Let $\mathsf \Delta= \mathsf T-\mathsf U$ and define the coupling-good event $\mathcal{G}_{\mathsf{cp}}$ as $\|\mathsf \Delta\|_2\leq n^{-0.3}$.
   By Markov's inequality, we obtain $\mathbb P[\mathcal{G}_{\mathsf{cp}}] \geq 1-n^{-0.17}$.
     Now we derive the concentration property of $\SMT(\beta)$. Define $$\mathsf E(\beta)\triangleq\sum_{k=1}^M (1+\sigma^2)^{-\frac k2}\mathbb E[\mathsf U_k\mathsf U_{M+k}]=\sum_{k=1}^M \left(\frac{\langle\beta,\beta_*\rangle}{1+\sigma^2}\right)^k.$$ Then for any  $\epsilon>0$,
    \begin{equation}
    \begin{split}
        &\mathbb P_1\left[ |\SMT(\beta) - \mathsf E(\beta)| \geq \epsilon M \right] \leq \mathbb P_{\mathsf{cp}}\left[ |\SMT(\beta) - \mathsf E(\beta)| \geq \epsilon M,\mathcal{G}_{\mathsf{cp}} \right] + \mathbb P[\mathcal{G}_{\sf{cp}}^c]\\
        \leq & \mathbb P_{\mathsf{cp}}\left[ \Big|\sum_{k=1}^M {(1+\sigma^2)^{-\frac{k}{2}}} \left( \mathsf U_k+\mathsf\Delta_k \right)\left( \mathsf U_{M+k}+\mathsf\Delta_{M+k} \right)-\mathsf E(\beta) \Big|\geq \epsilon M , \mathcal{G}_{\mathsf{cp}}\right]+{n^{-0.17}}\\
        \leq & \mathbb P_{\mathsf{cp}}\left[ \Big|\sum_{k=1}^M {(1+\sigma^2)^{-\frac{k}{2}}}  \mathsf U_k  \mathsf U_{M+k} -\mathsf E(\beta) \Big|\geq \epsilon M -{n^{-0.6}}- {n^{-0.3}}\|\mathsf U\|_2\right]+{n^{-0.17}}\,,
    \end{split}\label{eq-H1-SWT-recover-beta-1}
    \end{equation}
    where in the last inequality we apply the Cauchy-Schwarz inequality to bound $\sum_{k=1}^M \mathsf \Delta_k\mathsf\Delta_{M+k}\leq \|\mathsf\Delta\|_2^2$ and $\sum_{k=1}^M \mathsf\Delta_k\mathsf U_{M+k}+\sum_{k=1}^M \mathsf\Delta_{M+k}\mathsf U_k \leq \|\mathsf\Delta\|_2\|\mathsf U\|_2$. By Bernstein inequality for the sum of sub-exponential variables, 
    \begin{align}
        &\mathbb P\left[ \| \mathsf{U}\|_2 \geq n^{0.1} \right]\leq \exp(-n^{0.1})\,.\notag\\
        & \mathbb P\left[  \Big|\sum_{k=1}^M (1+\sigma^2)^{-\frac{k}{2}}\mathsf U_k\mathsf U_{k+M} -\mathsf E(\beta) \Big| \geq \epsilon M/2 \right]\leq \exp(-\Theta(M))\notag\,,
    \end{align}
    where the first line follows from $M=n^{o(1)}$ by \eqref{eq-def-M-d-large}, and the second line uses the fact that $\{\mathsf U_k \mathsf U_{k+M}\}$ are independent sub-exponential random variables. Plugging the above tail probability into \eqref{eq-H1-SWT-recover-beta-1}, we get that
    \begin{equation*}
        \mathbb P_1\left[ |\SMT(\beta) - \mathsf E(\beta)| \geq \epsilon M \right] \leq \exp(-\Theta(M)) + \exp(-n^{0.1}) + n^{-0.17}=\exp(-\Theta(M)),
    \end{equation*}
    where the last inequality holds due to $M=o(\log n)$. Applying a union bound yields
    \begin{equation}\label{eq:concentration-SMTbeta}
        \mathbb P_1\left[ \forall \beta\in N_\delta\,,\, |\SMT(\beta) - \mathsf E(\beta)| \leq \epsilon M \right]= 1-o(1).
    \end{equation}
    It remains to compare the deterministic profiles $\{\mathsf E(\beta):\,\beta\in N_{\delta,2}\} $ with $\mathsf c_*M$. Here $N_{\delta,2}\triangleq\{\beta\in N_\delta:\|\beta-\beta_*\|_2>10/\sqrt M\}$.
    For any $\beta\in N_{\delta,2}$, $\langle\beta,\beta_*\rangle\leq 1- \frac{50}{M}$. Therefore,
    \begin{equation}
        \mathsf E(\beta) = \sum_{k=1}^M \left( \frac{\langle \beta,\beta_*\rangle}{1+\sigma^2} \right)^k = \frac{\frac{\langle \beta,\beta_*\rangle}{1+\sigma^2} \left( 1-\left(\frac{\langle \beta,\beta_*\rangle}{1+\sigma^2}\right)^M \right) }{1-\frac{\langle \beta,\beta_*\rangle}{1+\sigma^2}} \leq \frac{ \left(1-\frac{50}{M}\right) \left( 1-\left(\frac{\langle \beta,\beta_*\rangle}{1+\sigma^2}\right)^M \right) }{\frac{50}{M}}\leq \frac{M}{50}\,.\notag
    \end{equation}
    As a result, $c_*M-\mathsf E(\beta)\geq \frac{M}{10}$ for any $\beta\in N_{\delta,2}$. Taking $\epsilon = 1/20$ in \eqref{eq:concentration-SMTbeta}, we get that
    \begin{equation}
         \SMT(\beta) \leq \mathsf c_* M \text{ for any }\beta\in N_{\delta,2}\notag
    \end{equation}
    holds with probability $1-o(1)$. 
    Consequently, $\|\tilde\beta-\beta_*\|_2\le 10/\sqrt{M}$ with probability $1-o(1)$. Since $\tilde\beta$ and $\beta_*$ are unit vectors, $\|\tilde\beta-\beta_*\|_2^2\le 4$, and therefore $\mathbb E\|\tilde\beta-\beta_*\|_2^2\le 100/M+4\,\mathbb P[\|\tilde\beta-\beta_*\|_2>10/\sqrt{M}]=o(1)$.

\subsection{Proof of Theorem~\ref{thm:impossibility-beta-recovery}}
\label{sec:impossibility-beta-recovery}

We prove the result by analyzing the posterior distribution of $\beta_*$. Let $\nu$ denote the uniform measure on the unit sphere, and define
\begin{equation}
    \label{eq-def-posterior-beta*}
    \mu_{X,Y}(\mathrm d\beta)\triangleq\mathbb P_1[\beta_*\in \mathrm d\beta\mid X,Y]\,.
\end{equation}
Let $\widetilde\beta$ be an independent sample from the posterior $\mu_{X,Y}$. We first show that, for $\iota_n=d^{-1/4}$,
\begin{equation}\label{eq-final-goal-recovery-impossibility}
    \mathbb P_1\Big[\langle\beta_*,\widetilde\beta\rangle\geq \iota_n\Big] = o(1)\,.
\end{equation}
We begin by writing the posterior density with respect to $\nu$. By Bayes' formula, for $\nu$-a.e. $\beta$,
\begin{equation}\label{eq-expression-posterior}
\frac{\mathrm d\mu_{X,Y}}{\mathrm d\nu}(\beta)
=
\frac{L_\beta(X,Y)}{L(X,Y)},
\end{equation}
where we recall
\[L_\beta(X,Y)\triangleq\mathbb E_{\Pi}\left[\frac{\mathbb P_1[X,Y\mid \beta,\Pi]}{\mathbb P_0[X,Y]}\right], \qquad L(X,Y)=\mathbb E_{\beta\sim\nu}L_\beta(X,Y).
\]
Here $L(X,Y)$ is the likelihood ratio defined in \eqref{eq:def-likelihood-and-Lbeta}.
In addition, for any positive sequence $\epsilon_n$,
\begin{equation}\label{eq-small-likelihood-ratio}
    \mathbb P_1[L(X,Y)\leq\epsilon_n] = \mathbb E_0[ L\mathbf{1}_{\{ L(X,Y)\leq\epsilon_n\}} ]\leq \epsilon_n\,.
\end{equation}

For fixed $\beta_*$, let $\mathcal A(\beta_*)$ denote the spherical cap consisting of all unit vectors $\beta$ with $\langle \beta,\beta_*\rangle\ge \iota_n$. Using \eqref{eq-expression-posterior}, we obtain
\begin{equation}
    \begin{split}
       & \mu_{X,Y}(\mathcal{A}(\beta_*))  =\mu_{X,Y}(\mathcal{A}(\beta_*))\mathbf{1}_{\{L(X,Y)\leq\epsilon_n\}}+\mu_{X,Y}(\mathcal{A}(\beta_*))\mathbf{1}_{\{L(X,Y)>\epsilon_n\}}\\
       \leq & \mathbf{1}_{\{L(X,Y)\leq\epsilon_n\}} + \frac{1}{\epsilon_n} \mathbb E_{\beta\sim\nu}\left[ L_\beta(X,Y) \mathbf{1}_{\{\beta\in\mathcal{A}(\beta_*)\}} \right] \,,\notag
    \end{split}
\end{equation}
where on the event $\{L(X,Y)>\epsilon_n\}$ we used $L(X,Y)^{-1}\le \epsilon_n^{-1}$, and on its complement we used $\mu_{X,Y}(\mathcal A(\beta_*))\le 1$. 
Taking expectation under $\mathbb P_1$ and using \eqref{eq-small-likelihood-ratio}, we get
\begin{align}
\mathbb P_1\left[\langle \widetilde\beta,\beta_*\rangle\ge \iota_n\right] & \le\epsilon_n+\frac{1}{\epsilon_n}\mathbb E_{\beta,\beta_*\sim\nu}\left[g(\beta_*,\beta)\mathbf 1_{\{\langle\beta,\beta_*\rangle\ge \iota_n\}}\right] \nonumber \\
& \le 2 \sqrt{\mathbb E_{\beta,\beta_*\sim\nu}\left[g(\beta_*,\beta)\mathbf 1_{\{\langle\beta,\beta_*\rangle\ge \iota_n\}}\right]}
\le 2 \exp\left(-\Omega(\sqrt{d}) \right),
\end{align}
where the first inequality holds by recalling $g(\beta_*,\beta)=\mathbb E_0[L_{\beta_*}(X,Y)L_\beta(X,Y)]$; the second inequality holds by optimizing $\epsilon_n$;
and the last inequality follows from Lemmas~\ref{lem:lower-bound-second-moment-before-integral}--\ref{lem:lower-bound-est-expectation}, 
together with $\big(\frac{1+\lambda}{2}\big)^n\big(\frac{1+\lambda}{1-\lambda}\big)^2\le\big(\frac{2(1+\sigma^2)}{\sigma^2}\big)^2=O((d\log d)^2)$, which is absorbed into $e^{-c\sqrt d}$.
It follows that, as $d \to \infty,$
 \begin{equation}\label{eq-expect-inner-product-independent-posterior}
 \mathbb E_1\left[\langle \widetilde\beta,\beta_*\rangle\right]\le\iota_n+\mathbb P_1\left[\langle \widetilde\beta,\beta_*\rangle\ge \iota_n\right]=o(1).
\end{equation}

Finally, recall that MMSE is achieved by the posterior mean. Therefore,  
\begin{equation}
\begin{split}
    \mathsf{MMSE}=
    \mathbb E_1\left[\left\|\beta_*-\mathbb E_{\beta\sim\mu_{X,Y}}[\beta] \right\|_2^2\right] 
    &=1-\mathbb E_1\left[\langle \widetilde\beta,\beta_*\rangle\right]=1-o(1),
\end{split}
\notag
\end{equation}
where the second equality holds because
$$\mathbb E_1\left[\langle \widetilde\beta,\beta_*\rangle\right]
= \mathbb E_1\left[ \langle \beta_*, 
\mathbb E_{\beta\sim\mu_{X,Y}}[\beta]  \rangle \right]
=\mathbb E_1\left[ \langle \mathbb{E}_1[\beta_* \mid X,Y], 
\mathbb E_{\beta\sim\mu_{X,Y}}[\beta] \rangle  \right]
=\mathbb E_1\left[\norm{\mathbb E_{\beta\sim\mu_{X,Y}}[\beta]}^2\right].
$$

\subsection{Proof of Theorem~\ref{thm:sparse-beta}}
\label{sec:sparse-beta-test}

\paragraph{Positive part.}

The construction of the \emph{sparse sliced-moment test} follows the same idea as the \emph{sliced-moment test} in Theorem~\ref{thm:SMT}. The only difference is that the search space is now the discrete sparse set $\mathbb A_{\sf{sparse}}$, rather than the unit sphere in $\mathbb R^d$. Therefore, we only sketch the proof and highlight the necessary modifications.
The proofs of \eqref{eq-final-goal-H1-d-large} and \eqref{lem:eachbeta-analysis-H0-d-large} use $M=o(\log n)$ only through $e^{O(M)}\le n^{0.01}$. Here
$M'\le C_0\log\binom ds\le 2C_0\, s\log d\le 2C_0c_0\log n$, so this holds once $c_0$ is small enough.

For the analysis under $\mathcal H_1$, as in the proof of Theorem~\ref{thm:SMT}, we evaluate the statistic at the closest element to $\beta_*$ in the search space. In the sparse setting, $\mathbb A_{\sf{sparse}}$ is discrete and contains $\beta_*$, so this closest element is simply $\widehat\beta=\beta_*$. Define $\SSMT(\beta)
\triangleq
\sum_{k=1}^{M'} T_k S_k(\beta).$
Then
\[
\mathbb E_1\left[\SSMT(\beta_*)\right]
=
\sum_{k=1}^{M'}
\mathbb E_1\left[T_k S_k(\beta_*)\right]
=
\sum_{k=1}^{M'}
\frac{1}{(1+\sigma^2)^k}
\ge
\frac{M'}{(1+\sigma^2)^{M'}}
\ge
e^{-1}M',
\]
where the last inequality follows from the choice of $M'$. By the same second-moment argument used to derive \eqref{eq-final-goal-H1-d-large}, we obtain
\begin{equation}\label{eq:SSMT-H1}
    \mathbb P_1\left[
    \SSMT \le \mathsf c' M'
    \right]
    \le
    \mathbb P_1\left[
    \SSMT(\beta_*) \le \mathsf c' M'
    \right]
    =
    o(1),
\end{equation}
with $\mathsf c'=(2e)^{-1}$.
It remains to control the statistic under $\mathcal H_0$. Following the same proof of \eqref{lem:eachbeta-analysis-H0-d-large}, for any fixed $\beta\in\mathbb A_{\sf{sparse}}$, we have
\begin{equation}\label{lem:eachbeta-analysis-H0-d-large-sparse}
    \mathbb P_0\left[
    \SSMT(\beta)\ge \mathsf c'M'
    \right]
    \le
    e^{-\kappa M'}+3n^{-0.19},
\end{equation}
where $\kappa>0$ is the constant given by Bernstein's inequality in the proof of \eqref{lem:eachbeta-analysis-H0-d-large}, which depends only on $\mathsf c'$.
Applying a union bound over $\mathbb A_{\sf{sparse}}$ gives
\begin{align}
    \mathbb P_0\left[\SSMT\ge \mathsf c'M'\right]
    &\le \sum_{\beta\in\mathbb A_{\sf{sparse}}}
    \mathbb P_0\left[\SSMT(\beta)\ge \mathsf c'M'\right] \notag\\
    &\le \binom ds\left(e^{-\kappa M'}+3n^{-0.19}\right) \notag\\
    &\le e^\kappa \binom ds^{1-\kappa C_0}+3n^{2c_0-0.19}=o(1),
    \label{eq:SSMT-H0}
\end{align}
where we used $|\mathbb A_{\sf{sparse}}|=\binom ds\le e^{2s\log d}\le n^{2c_0}$ and
$M'=\lfloor C_0\log\binom ds\rfloor$, which holds since $\sigma^2\le c/\log\binom ds$ with $c\le 1/C_0$.
The last step holds once $C_0\ge 2/\kappa$ (so that $e^\kappa\binom ds^{1-\kappa C_0}\le e^\kappa\binom ds^{-1}\to0$) and $c_0<0.09$.
Combining \eqref{eq:SSMT-H1} and \eqref{eq:SSMT-H0} completes the proof.

\paragraph{Negative part.}
We denote by $\mathbf{HG}(N,M,K)$ the hypergeometric distribution, i.e.,
\begin{equation}
    \mathbb P\left[ \mathbf{HG}(N,M,K)=k \right]=\frac{\binom{K}{k}\binom{N-K}{M-k}}{\binom{N}{M}}\,.\notag
\end{equation}
Recalling the definition of $\mathbb A_{\mathsf{sparse}}$ in Definition~\ref{def:slr_sparse}, we obtain that the distribution of $\langle\beta_1,\beta_2\rangle$ satisfies
\begin{equation}
    s\cdot\langle\beta_1,\beta_2\rangle \sim  \mathbf{HG}(d,s,s)\,.\notag
\end{equation}
The mean of $\mathbf{HG}(d,s,s)$ is given by
$\mu=\frac{s^2}{d}$.
Moreover, for any $t\ge 1$, the standard hypergeometric tail bound yields
\[
\mathbb P\left[\left\langle \beta_1,\beta_2\right\rangle \ge \frac{t}{s}\right]
=
\mathbb P\left[s\left\langle \beta_1,\beta_2\right\rangle \ge t\right]
\le
\left(\frac{e s^2}{dt}\right)^t.
\]
 Since $s\langle\beta_1,\beta_2\rangle$ is integer-valued, for $0<t\le 1$ the left-hand side equals $\mathbb P[s\langle\beta_1,\beta_2\rangle\ge 1]\le es^2/d$.
Consequently, uniformly over all regimes $s=o(d)$, by taking $t$ as below
\begin{equation}
    \label{eq:condition-t}
    \begin{split}
    & \frac{s^2}{d}\ll t\ll s\text{ if }s\to\infty\,,\\
    & t=1/2 \text{ if }s=O(1)\,,
    \end{split}
\end{equation}
we obtain
\begin{equation}\label{eq:inner-prod-betas-sparse}
    \mathbb P\left[\left\langle \beta_1,\beta_2\right\rangle \ge t/s\right]=o(1).
\end{equation}

The following lemma is crucial for controlling the truncated mean of $g(\beta_1,\beta_2)$. The proof is similar to that of Lemma~\ref{lem:lower-bound-est-expectation}. The main difference is that, in
the sparse setting, $s\langle\beta_1,\beta_2\rangle$ follows a hypergeometric distribution. As before, the most delicate contribution to \eqref{eq:sparse-truncated-mean-g} comes from the event $\langle\beta_1,\beta_2\rangle\approx1$.  Requiring this contribution to be $o(1)$
is precisely what leads to the lower-bound condition on the noise level.

\begin{lemma}
    \label{lem:sparse-truncated-mean-g}
   Recall $\lambda=\frac{|\langle\beta_1,\beta_2\rangle|}{1+\sigma^2}$. For $t$ defined in \eqref{eq:condition-t}. Assume that $\sigma^2\geq \frac{(1+\epsilon)\pi^2}{6\log\binom{d}{s}}$ for some constant $\epsilon>0$. 
    \begin{equation}
        \label{eq:sparse-truncated-mean-g}
        \mathbb E_{\beta_1,\beta_2}\left[ \exp\left( \frac{\pi^2}{6(1-\lambda)} \right)\mathbf 1_{|\langle\beta_1,\beta_2\rangle|>t/s} \right] = o(1)\,,
    \end{equation}
    where $\beta_1,\beta_2\overset{\mathrm{iid}}\sim \operatorname{Uniform}(\mathbb A_{\sf{sparse}})$.
\end{lemma}
\begin{proof}
    If $\sigma^2 \geq 1/4$, then $\lambda\leq 4/5$ and thus 
    \begin{equation}
        \begin{split}
             \mathbb E_{\beta_1,\beta_2}\left[ \exp\left( \frac{\pi^2}{6(1-\lambda)} \right)\mathbf 1_{|\langle\beta_1,\beta_2\rangle|>t/s} \right] &\leq \exp\left(\frac{5\pi^2}{6}\right) \mathbb P\left[ |\langle\beta_1,\beta_2\rangle|\geq t/s \right]\overset{\eqref{eq:inner-prod-betas-sparse}}{=}o(1).
        \end{split}\notag
    \end{equation}
    Thus we focus on $\sigma^2\leq 1/4$. Similarly if $|\langle\beta_1,\beta_2\rangle|\leq 4/5$, then $\lambda=|\langle\beta_1,\beta_2\rangle|/(1+\sigma^2)\leq 4/5$ and hence
    \begin{equation}
         \mathbb E_{\beta_1,\beta_2}\left[ \exp\left( \frac{\pi^2}{6(1-\lambda)} \right)\mathbf 1_{ t/s< |\langle\beta_1,\beta_2\rangle|\le 4/5} \right] \le \exp\left(\frac{5\pi^2}{6}\right)\mathbb P\left[ |\langle\beta_1,\beta_2\rangle|> t/s \right]=o(1).\notag
    \end{equation}
    Therefore, it remains to show 
    \begin{equation}\label{eq:final-goal-truncated-mean}
         \mathbb E_{\beta_1,\beta_2}\left[ \exp\left( \frac{\pi^2}{6(1-\lambda)} \right)\mathbf 1_{|\langle\beta_1,\beta_2\rangle|>4/5} \right] = o(1)\,.
    \end{equation}
    We write $X\triangleq s\langle\beta_1,\beta_2\rangle$, therefore $X\sim \mathbf{HG}(d,s,s)$ with
    \begin{equation}
        \mathbb P\left[ X = s-r\right] = \frac{\binom{s}{r}\binom{d-s}{r}}{\binom{d}{s}}\,.
    \end{equation}
    Let $\mathsf C=\pi^2/6$. Hence,
    \begin{equation}
        \mathbb E_{\beta_1,\beta_2}\left[ \exp\left( \frac{\pi^2}{6(1-\lambda)} \right)\mathbf 1_{|\langle\beta_1,\beta_2\rangle|>4/5} \right] = \sum_{ 0\leq  r < s/5} \exp\left( \frac{\mathsf C (1+\sigma^2)}{\sigma^2+\frac{r}{s} } \right)\frac{\binom{s}{r}\binom{d-s}{r}}{\binom{d}{s}}\,.
    \end{equation}
   If $\sigma^2\geq (1+\epsilon)\pi^2\Big/6\log\binom{d}{s}$, the $r=0$ term satisfies 
    \begin{equation}
        \exp\left( \frac{\mathsf C(1+\sigma^2)}{\sigma^2} - \log\binom{d}{s} \right) = o(1)\,.
    \end{equation}
    For $1\leq r\leq s/5$, we can bound binomial coefficients as 
    \begin{equation}
        \binom{s}{r}\binom{d-s}{r}\leq \left(\frac{e^2sd}{r^2}\right)^r\text{ and }\binom{d}{s}\geq \left(\frac{d}{s}\right)^s\,.\notag
    \end{equation}
    Recall that we assumed $\sigma^2\leq 1/4$ and $s=o(d)$. Therefore,
    \begin{equation}
        \exp\left( \frac{\mathsf C (1+\sigma^2)}{\sigma^2+\frac{r}{s} } \right)\frac{\binom{s}{r}\binom{d-s}{r}}{\binom{d}{s}} \leq \exp\left(\frac{2\mathsf Cs}{r}+r\log \left(\frac{e^2sd}{r^2}\right)- s\log\left(\frac{d}{s}\right)\right)\leq \exp\left( -\frac{s}{2}\log\left(\frac{d}{s}\right) \right).\notag
    \end{equation}
    For the second inequality, denote the middle exponent by $f(r)$. Then $f''(r)=4\mathsf Cs/r^3-2/r>0$ for $r<\sqrt{2\mathsf Cs}$, and $f'(r)=-2\mathsf Cs/r^2+\log(sd/r^2)>0$ for $\sqrt{2\mathsf Cs}\le r\le s/5$ once $d/s$ is large. Hence $\max_{1\le r\le s/5}f(r)=\max\{f(1),f(s/5)\}$, and both values are at most $-\frac s2\log\frac ds$ once $d/s$ is sufficiently large, uniformly in $s\ge5$.
    Summing the above expression over $1\leq r\leq s/5$ 
    \begin{equation}
       \sum_{1\leq r\leq s/5} \exp\left( \frac{\mathsf C (1+\sigma^2)}{\sigma^2+\frac{r}{s} } \right)\frac{\binom{s}{r}\binom{d-s}{r}}{\binom{d}{s}}= o(1)\,.\notag
    \end{equation}
    This shows \eqref{eq:final-goal-truncated-mean}.
\end{proof}
Now we are ready to present the proof of second Item of Theorem~\ref{thm:sparse-beta}.
\begin{proof}
    Suppose that $d\to\infty$, $s\log d\leq c_0\log n$ and $\sigma^2 \geq  \frac{(1+\varepsilon)\pi^2}{6\log\binom{d}{s}}$. Then
\begin{equation}
    \left( \frac{1+\lambda}{2} \right)^n\left(\frac{1+\lambda}{1-\lambda}\right)^2\leq \exp\left(-\frac{\sigma^2n}{2(1+\sigma^2)}\right)\left(\frac{2(1+\sigma^2)}{\sigma^2}\right)^2 = o(1)\,,\notag
\end{equation}
where the last step holds since $\log\binom ds\le 2s\log d\le 2c_0\log n$ implies $\sigma^2\gtrsim 1/\log n$.
Therefore, applying Lemma~\ref{lem:lower-bound-second-moment-before-integral} and Lemma~\ref{lem:sparse-truncated-mean-g}, we obtain 
\begin{align*}
    \mathbb E_0L_{\sf{sparse}}^2 \leq  &\left(1+o(1) \right)
   \mathbb E_{\beta_1,\beta_2}\left[\exp \left( 
        \frac{\lambda}{(1-\lambda)^2}\right) \mathbf{1}_{|\langle\beta_1,\beta_2\rangle|\le t/s}\right]  \\
         & +  (1+o(1)) \mathbb E_{\beta_1,\beta_2}\left[\exp \left( 
        \frac{\pi^2}{6(1-\lambda)}\right) \mathbf{1}_{|\langle\beta_1,\beta_2\rangle|>t/s}\right]  \\
    & \le 1+o(1)\,,\notag
\end{align*}
where the last inequality holds because $\lambda\leq t/s\to 0$ when $s\to\infty$ and $\lambda=0$ when $s=O(1)$. 
\end{proof}

\appendix

\section{Computation of the cross-moment test}
\label{sec:computation-CMT}
The cross-moment test \eqref{eq-testing-statistic(I)} requires $T_A$ and $S_A=\mathbb{E}_1[T_A\mid X]$ for $A\in\binom{[K]}{D}$. The coefficients
$T_1,\dots,T_K$ can be computed in $O(nK)$ time using the three-term recurrence
$H_{k+1}(x)=xH_k(x)-kH_{k-1}(x)$, so the only nontrivial step is $S_A$.
More generally, for positive integers $m=(m_1,\dots,m_D)$, 
\[
S_m \triangleq \mathbb{E}_1\Big[\prod_{j=1}^D T_{m_j}\,\Big|\,X\Big],
\]
so that $S_A=S_m$ when $m$ lists the elements of $A$. Throughout this section,
$\rho=(1+\sigma^2)^{-1/2}$, and for a finite set $B$ we write $\Pi(B)$ for the
set of partitions of $B$, i.e., collections of pairwise disjoint nonempty subsets of $B$
whose union is $B$.

\begin{lemma}\label{lem:SA-formula}
Let $m=(m_1,\dots,m_D)$ be positive integers and $M=\sum_{j=1}^D m_j$. For
nonempty $B\subseteq[D]$ and $u\in\mathbb{R}$, let $\kappa_B(u)$ be the joint
cumulant of the random variables $h_{m_j}(\rho(u+\sigma Z))$, $j\in B$, where
$Z\sim N(0,1)$. By the moment--cumulant formula (see, e.g.,
\cite{leonov1959method,peccati2011wiener}), these cumulants satisfy
\begin{equation}\label{eq:moment-cumulant}
\mathbb{E}_Z\Big[\prod_{j\in B} h_{m_j}\big(\rho(u+\sigma Z)\big)\Big]
=\sum_{\pi\in\Pi(B)}\ \prod_{C\in\pi}\kappa_C(u),
\qquad \emptyset\neq B\subseteq[D],
\end{equation}
where $\mathbb{E}_Z$ denotes the expectation over $Z$ with $u$ fixed.
Then each $\kappa_B$ is a polynomial of degree at most $\sum_{j\in B}m_j$, and
\begin{equation}\label{eq:SA-formula}
S_m=\frac{1}{n^{D/2}}\sum_{\pi\in\Pi([D])}
\mathbb{E}_\beta\Big[\prod_{B\in\pi}\ \sum_{i=1}^n\kappa_B(\langle X_i,\beta\rangle)\Big],
\end{equation}
where $\beta$ is uniform on $\mathbb{S}^{d-1}$. Moreover, each $S_m$ can be computed
in $O\big(nN+D^{D+1}dN^2\big)$ time, where
$N=\binom{M+d}{d}$.
\end{lemma}

\begin{proof}
We first remove the permutation. Under $\mathcal{H}_1$, we have
$Y_{\pi_*^{-1}(i)}=\langle X_i,\beta_*\rangle+\sigma Z_{\pi_*^{-1}(i)}$ for every
$i\in[n]$, and $(Z_{\pi_*^{-1}(i)})_{i\in[n]}$ are again i.i.d.\ $N(0,1)$ and
independent of $(X,\beta_*)$. Since each $T_k$ depends on $Y$ only through the
multiset $\{Y_1,\dots,Y_n\}$, and $Y_i/\sqrt{1+\sigma^2}=\rho Y_i$, it follows
that, conditionally on $X$,
\[
\prod_{j=1}^D T_{m_j}\ \overset{d}{=}\ \frac{1}{n^{D/2}}\prod_{j=1}^D\sum_{i=1}^n
h_{m_j}\big(\rho(u_i+\sigma Z_i)\big),\qquad u_i\triangleq\langle X_i,\beta\rangle,
\]
where $\beta\sim\mathrm{Unif}(\mathbb{S}^{d-1})$ and $Z_1,\dots,Z_n\overset{iid}{\sim}N(0,1)$
are independent. In particular the permutation plays no role in $S_m$.

Next we take the expectation over the noise. Fix $\beta$. Expanding the product,
\begin{equation}\label{eq:SA-expand}
\mathbb{E}_Z\Big[\prod_{j=1}^D\sum_{i=1}^n h_{m_j}\big(\rho(u_i+\sigma Z_i)\big)\Big]
=\sum_{(i_1,\dots,i_D)\in[n]^D}\mathbb{E}_Z\Big[\prod_{j=1}^D h_{m_j}\big(\rho(u_{i_j}+\sigma Z_{i_j})\big)\Big].
\end{equation}
For a fixed tuple $(i_1,\dots,i_D)$, let $B_i=\{j: i_j=i\}$. Factors with the
same index $i$ share the noise $Z_i$, while factors with different indices are
independent. Therefore, by \eqref{eq:moment-cumulant},
\[
\mathbb{E}_Z\Big[\prod_{j=1}^D h_{m_j}\big(\rho(u_{i_j}+\sigma Z_{i_j})\big)\Big]
=\prod_{i:\,B_i\neq\emptyset}\mathbb{E}_Z\Big[\prod_{j\in B_i}h_{m_j}\big(\rho(u_i+\sigma Z)\big)\Big]
=\prod_{i:\,B_i\neq\emptyset}\ \sum_{\pi_i\in\Pi(B_i)}\ \prod_{C\in\pi_i}\kappa_C(u_i).
\]
Substituting this into \eqref{eq:SA-expand} and expanding, each term is
determined by the partition $\pi=\bigcup_i\pi_i$ of $[D]$ and the map
$\tau:\pi\to[n]$ given by $\tau(C)=i$ for $C\in\pi_i$, and every such pair
$(\pi,\tau)$ arises exactly once. Hence
\[
\text{\eqref{eq:SA-expand}}
=\sum_{\pi\in\Pi([D])}\ \sum_{\tau:\,\pi\to[n]}\ \prod_{C\in\pi}\kappa_C(u_{\tau(C)})
=\sum_{\pi\in\Pi([D])}\ \prod_{C\in\pi}\ \sum_{i=1}^n\kappa_C(u_i).
\]
Multiplying by $n^{-D/2}$ and averaging over $\beta$ yields \eqref{eq:SA-formula}.

It remains to show that each $\kappa_B$ is a polynomial and to count the
operations needed to evaluate \eqref{eq:SA-formula}. Let $q_B(u)$ denote the left-hand side of \eqref{eq:moment-cumulant}, and
$q_\emptyset=1$. Writing $\prod_{j\in B}h_{m_j}(w)=\sum_k\gamma_k w^k$ (with $\gamma_k$ depending on $B$) and using
$\mathbb{E}[Z^{2r}]=(2r-1)!!$,
\begin{equation}\label{eq:qB}
q_B(u)=\sum_k\gamma_k\,\rho^k\sum_{r\le k/2}\binom{k}{2r}(2r-1)!!\,\sigma^{2r}u^{k-2r},
\end{equation}
a polynomial of degree at most $\sum_{j\in B}m_j$. Fix $b\in B$. Each
partition of $B$ consists of the block $A\ni b$ together with a partition of
$B\setminus A$, so applying \eqref{eq:moment-cumulant} to $B\setminus A$ gives
\begin{equation}\label{eq:kappa-rec}
q_B=\sum_{A\subseteq B,\ A\ni b}\kappa_A\,q_{B\setminus A}.
\end{equation}
The term $A=B$ equals $\kappa_B$, since $q_\emptyset=1$, and all other terms
involve only $\kappa_A$ with $|A|<|B|$. Hence
\[
\kappa_B=q_B-\sum_{A\subseteq B,\ A\ni b,\ A\neq B}\kappa_A\,q_{B\setminus A},
\]
which determines all $\kappa_B$ recursively in increasing order of $|B|$. By
induction on $|B|$, each $\kappa_B$ is a polynomial of degree at most
$\sum_{j\in B}m_j$.

For a multi-index $\alpha\in\mathbb{Z}_{\ge0}^d$, write
$x^\alpha=\prod_{t=1}^d x_t^{\alpha_t}$ and $|\alpha|=\sum_t\alpha_t$. If
$\kappa_B(u)=\sum_r c_{B,r}u^r$, the multinomial theorem gives
\begin{equation}\label{eq:KB-expand}
\sum_{i=1}^n\kappa_B(\langle X_i,\beta\rangle)
=\sum_{|\alpha|\le \sum_{j\in B}m_j} c_{B,|\alpha|}\binom{|\alpha|}{\alpha}
\Big(\sum_{i=1}^n X_i^\alpha\Big)\beta^\alpha .
\end{equation}
Thus the expression inside $\mathbb{E}_\beta$ in \eqref{eq:SA-formula} is a
polynomial in $\beta$ of degree at most $M$, hence with at most $N$
coefficients, and $X$ enters only through the empirical moments
$\sum_{i}X_i^\alpha$, $|\alpha|\le M$. To average over $\beta$, let
$z\sim N(0,I_d)$; since $z=\|z\|\beta$ with $\|z\|$ independent of $\beta$, we
have $\mathbb{E}[z^\alpha]=\mathbb{E}\|z\|^{|\alpha|}\,\mathbb{E}[\beta^\alpha]$.
Applying moment formula for Gaussian and $\|z\|^2\sim\chi_d^2$, we get
\begin{equation}\label{eq:sphere-moments}
\mathbb{E}[\beta^\alpha]=
\begin{cases}
\dfrac{\prod_{t=1}^d(\alpha_t-1)!!}{\prod_{r=0}^{|\alpha|/2-1}(d+2r)}, & \text{all $\alpha_t$ even},\\[2ex]
0, & \text{otherwise},
\end{cases}
\end{equation}
with the convention $(-1)!!=1$.

We now address the time complexity. The only step 
in the preceding calculation that depends on $n$ is computing the
empirical moments $\sum_i X_i^\alpha$, $|\alpha|\le M$, which costs $O(nN)$,
since each monomial $X_i^\alpha$ is obtained from one of degree $|\alpha|-1$ by
a single multiplication. All remaining steps operate on polynomials with at most
$N$ coefficients, and their cost is negligible compared to $n$. First, all $q_B$ and
$\kappa_B$, $B\subseteq[D]$, are obtained from \eqref{eq:qB} and
\eqref{eq:kappa-rec} with $O(3^DM^2)$ operations, and then \eqref{eq:KB-expand}
gives the coefficients of each $\sum_i\kappa_B(\langle X_i,\beta\rangle)$ with
$O(N)$ operations. Next, since $\mathbb{E}_\beta$ is evaluated monomial by
monomial via \eqref{eq:sphere-moments}, for each partition $\pi$ in
\eqref{eq:SA-formula} we expand the product
$\prod_{B\in\pi}\sum_i\kappa_B(\langle X_i,\beta\rangle)$ into monomials in
$\beta$, multiplying in one factor at a time and combining like terms. Each
partial product has degree at most $M$, hence at most $N$ coefficients, so each
of the at most $D$ multiplications involves at most $N^2$ pairs of monomials,
each multiplication requiring the addition of two exponent vectors in $\mathbb{Z}^d$; this
costs $O(dN^2)$. Finally, \eqref{eq:sphere-moments} is applied to each of the at
most $N$ resulting coefficients. Since $|\Pi([D])|\le D^D$, $M\le N$ and
$3^D=O(D^{D+1})$, the total number of operations 
for evaluating \prettyref{eq:SA-formula} 
is $O(nN+D^{D+1}dN^2)$.
\end{proof}

\paragraph{Running time of the cross-moment test.}
Under \eqref{eq-choice-par-CMT}, we have $D=30d$, $K\le D^2(\log n)^{0.01}$ and
$d\le(\log n)^{0.3}$, so $\log D,\log K=O(\log\log n)$. Every
$A\in\binom{[K]}{D}$ has $M=\sum_{k\in A}k\le DK$, so in
Lemma~\ref{lem:SA-formula}, $N\le\binom{DK+d}{d}\le(e(30K+1))^d=n^{o(1)}$.
The empirical moments with $|\alpha|\le DK$ are shared by all $A$ and cost
$n^{1+o(1)}$ in total; the remaining cost of all $S_A$ is
$O\big(\binom{K}{D}D^{D+1}dN^2\big)=\exp\big(O(d\log\log n)\big)=n^{o(1)}$.
Together with the $O(nK)$ cost of $T_1,\dots,T_K$, the cross-moment test is
computed in $n^{1+o(1)}$ time.

\section{Useful facts on Hermite polynomials}
\label{sec:auxiliary-results}

We provide a useful lemma of Hermite polynomials, as below.
\begin{lemma}\label{lem:expectation-over-beta}
Fix $z,x \in \reals^d$.
Then
\[
\|z\|^m
H_m\pth{\frac{\Iprod{z}{x}}{\|z\|}}
=\sum_{m_1+\ldots+m_d =m} \binom{m}{m_1, \ldots, m_d} \prod_{i=1}^d z_i^{m_i} H_{m_i}(x_i).
\]
As a result, for any $\beta$ in the unit sphere of $\mathbb R^d$,
\[
H_m\pth{\Iprod{\beta}{x}}
=\sum_{m_1+\ldots+m_d =m} \binom{m}{m_1, \ldots, m_d} \prod_{i=1}^d \beta_i^{m_i} H_{m_i}(x_i).
\]
\end{lemma}
\begin{proof}
    The form of the identity suggests some application of multinomial expansion. 
    Indeed, recall the complex integral representation of Hermite polynomials:
    $H_m(x) = \Expect[(x+\mathbf i W)^m]$, where $W\sim N(0,1)$ and $\mathbf i = \sqrt{-1}$.
    Without loss of generality, assume that $\|z\|=1$. 
    Pick $W=(W_1,\ldots,W_d) \sim N(0,\mathbf{I}_d)$.
    Then $\Iprod{z}{W}\sim N(0,1)$ and 
    \[
    H_m(\Iprod{z}{x})
    = \Expect[(\Iprod{z}{x+\mathbf iW})^m]
    \]
    which yields the desired identity applying multinomial theorem and the integral representation for each term.    
\end{proof}

\begin{lemma}\label{lem:scaling-law-hermite}
For each non-negative integer $k$, we have
\begin{equation}
    H_k(\sqrt{a}x+\sqrt{1-a}X) = \sum_{\ell=0}^k \binom{k}{\ell}a^{\ell/2}(1-a)^{(k-\ell)/2}H_\ell(x)H_{k-\ell}(X)\,.\notag
\end{equation}
\end{lemma}
\begin{proof}
    Note that for $H_n(y)$, we have the following expansion
    \begin{equation}\label{eq-expansion-representation-hermite}
        e^{yt-\frac{1}{2}t^2}=\sum_{n\geq 0}H_n(y)\frac{t^n}{n!}\,.
    \end{equation}
    Let $y=\sqrt{a}x+\sqrt{1-a}X$. Then we have
    \begin{equation}\label{eq-expansion-hermite}
        \exp\left( \left( \sqrt{a}x+\sqrt{1-a}X \right)t-\frac{1}{2}t^2  \right) = \sum_{n\geq 0} H_n\left(\sqrt{a}x+\sqrt{1-a}X \right) \frac{t^n}{n!}\,.
    \end{equation}
    Moreover, we can rewrite left-hand side as
    \begin{equation}\label{eq-expansion-right-hand side}
    \begin{split}
        &\exp\left( (\sqrt{a}t)x-\frac{(\sqrt{a}t)^2}{2} \right)\exp\left( (\sqrt{1-a}t)X-\frac{(\sqrt{1-a}t)^2}{2} \right) \\
        = &\left(\sum_{p\geq 0}H_p(x)\frac{(\sqrt{a}t)^p}{p!}\right)\left(\sum_{q\geq 0}H_q(X)\frac{(\sqrt{1-a}t)^q}{q!}\right)\,,
    \end{split}
    \end{equation}
    where the equality follows from \eqref{eq-expansion-representation-hermite}. 
    Then the desired result follows by comparing the coefficients of $t^k$ on the right-hand sides of \eqref{eq-expansion-right-hand side} and \eqref{eq-expansion-hermite}.
\end{proof}

\begin{lemma}[Hermite polynomials of correlated Gaussians, {\cite[Lemma~16]{fan2026empirical}}]\label{lem:hermite_cov}
Let $H_k$ denote the monic degree-$k$ Hermite polynomials. Let $(Z_1,Z_2)$ 
be jointly normally distributed with $Z_1,Z_2\sim N(0,1)$ and $\Expect[Z_1Z_2]=\rho \in [-1,1]$. Then for $k,k'\geq 1$,
\begin{align}
\Cov(H_k(Z_1),H_{k'}(Z_2)) = \bm{1}_{k=k'}\rho^kk!,
\label{eq:cov-Hermite1}
\end{align}
and for any $\theta_1,\theta_2\in\reals$,
\begin{align}
\Cov(H_k(\theta_1+Z_1),H_{k}(\theta_2+Z_2)) = \sum_{\ell=1}^k {k\choose \ell}^2 \ell! \cdot \theta_1^{k-\ell}\theta_2^{k-\ell}\rho^\ell.
\label{eq:cov-Hermite2}
\end{align}
\end{lemma}
The following lemma gives an estimate for the $\ell$-th moment of the $k$-th normalized Hermite polynomial.
\begin{lemma}
    \label{lem:fourth-moment-hermite-estimate}
    Let $X\sim N(0,1)$, then we have the following estimate:
    \begin{equation}
        \mathbb E_X\left[ |h_k(X)|^{\ell} \right] \leq (2\sqrt{e})^{k\ell} \exp\left(\frac{k\ell}{2}\log (\ell)\right)\,.\notag
    \end{equation}
\end{lemma}
\begin{proof}
    For the monic Hermite polynomial $H_k(x)$, we have the complex integral representation: $H_k(x)=\mathbb E_{W}\left[ (x+\mathbf iW)^{k} \right]$ where $W\sim N(0,1)$. And also note that $h_k(X)=H_k(X)/\sqrt{k!}$. Therefore, for any positive integer $k$ and $\ell$
    \begin{equation}\label{eq-moments-Hermite-Holder}
        \begin{split}
            & \mathbb E_{X}\left[ |h_k(X)|^{\ell} \right] = \frac{1}{(k!)^{\frac{\ell}{2}}}\mathbb E_X\left[ \left| \mathbb E_W\left[ \left(X+\mathbf iW\right)^{k} \right] \right|^{\ell} \right] \\
            \leq & \frac{1}{ (k!)^{\frac{\ell}{2}} }\mathbb E_{X,W}\left[ |X+\mathbf iW|^{k\ell} \right]   = \frac{1}{ (k!)^{\frac{\ell}{2}} }\mathbb E_{X,W}\left[ (X^2+W^2)^{\frac{k\ell}{2}} \right]  \,,
        \end{split}
    \end{equation}
    where the first inequality follows from H\"older's inequality. 
    Note that $X^2+W^2$ is distributed as chi-square distribution $\chi^2(2)$, we have
    \begin{equation}
        \mathbb E_{X,W}\left[ (X^2+W^2)^{\frac{k\ell}{2}} \right] = 2^{\frac{k\ell}{2}} \frac{\Gamma(1+\frac{k\ell}{2})}{\Gamma(1)} \leq 2^{k\ell} \exp\left( \frac{k\ell}{2}\log\left(\frac{k\ell}{2}\right)  \right)\,.\notag
    \end{equation}
    Plugging the above estimate into \eqref{eq-moments-Hermite-Holder} and applying the Stirling's approximation, we obtain 
    \begin{equation}
        \mathbb E_X\left[ |h_k(X)|^{\ell} \right] \leq \frac{2^{k\ell} \exp\left( \frac{k\ell}{2}\log\left(\frac{k\ell}{2}\right)  \right)}{ \exp\left( \frac{k\ell}{2}\log\left( k \right) - \frac{\ell k}{2} +\frac{\ell}{4}\log k  \right) } \leq (2\sqrt{e})^{k\ell} \exp\left(\frac{k\ell}{2}\log (\ell)\right)\,.\notag
    \end{equation}
    The proof is complete.
\end{proof}

\section{Normal approximation of functions of Gaussian variables}
\label{sec:normal-approx}
The following lemma~\cite[Proposition~4.3]{nourdin2010multivariate} provides a quantitative normal approximation bound for functions of Gaussian variables in the Wasserstein-1 distance. 
\begin{lemma}\label{lem:CLT-Stein-lemma}
Let $Y \sim N(\mathbf 0,\mathbf K)$, where 
$\mathbf K $ is an $N \times N$ positive definite matrix. Consider absolutely continuous functions $f_a : \mathbb{R}^N\to \mathbb{R}$ and let $W_a=f_a(Y)$ for $a=1,\ldots,m$. Assume that each random variable $W_a$ has zero mean, and that each function $f_a$ and its partial derivatives have sub-exponential growth at infinity.  
Let $Z \sim N(\mathbf 0,\mathbf C)$, where $\mathbf C$ is an $m \times m$ positive definite covariance matrix.
Then
\[
\mathcal W_1(W,Z) \le 
\|\mathbf C^{-1}\|_{\operatorname{op}} \,\|\mathbf C\|_{\operatorname{op}}^{1/2} \, 
\sqrt{\sum_{a,b=1}^m \mathbb{E} \big[ (\mathbf C_{ab} - T_{ab}(Y))^2 \big]},
\]
where the functions $T_{ab}(\cdot)$ are defined as
\begin{align}
T_{ab}(y) = \int_0^1 \frac{1}{2\sqrt{t}} \sum_{i,j=1}^N \mathbf K_{ij} \, 
\mathbb{E}\Big[ 
\frac{\partial f_a}{\partial y_i}(y) \, 
\frac{\partial f_b}{\partial y_j}(\sqrt{t} y + \sqrt{1-t} Y) 
\Big] \, \mathrm dt. \label{eq:T_a_b}
\end{align}
\end{lemma}

 The next lemma specializes the above normal approximation results to Hermite polynomials of Gaussian variables. 
\begin{lemma}
    \label{lem:gaussian-approximation-hermite-w/correlation}
    Let $\{(\mathsf X_i,\mathsf Y_i)\}_{i=1}^n$ be iid pairs of Gaussian variables with covariance $\big(\begin{smallmatrix} 1 & \rho \\ \rho & 1 \end{smallmatrix}\big)$, and $\{(\mathsf U_k,\mathsf V_k)\}_{k=1}^M$ be independent pairs of Gaussian variables with covariance $\big(\begin{smallmatrix} 1 & \rho^k \\ \rho^k & 1 \end{smallmatrix}\big)$. Let $\mathsf P_k=\frac{1}{\sqrt{n}}\sum_{j=1}^n h_k(\mathsf X_j)$ and $\mathsf Q_k=\frac{1}{\sqrt{n}}\sum_{j=1}^n h_k(\mathsf Y_j)$. Then
    \begin{equation}
        \mathcal{W}_1\left( \{(\mathsf P_k,\mathsf Q_k)\}_{k=1}^M,\{(\mathsf U_k,\mathsf V_k)\}_{k=1}^M \right) \leq \frac{\sqrt{1+|\rho|}}{1-|\rho|}\frac{\exp(cM)}{\sqrt n}\,,\notag
    \end{equation}
    where $c$ is an absolute constant independent of $\rho$ and $n$.
\end{lemma}
\begin{proof}
    We use Lemma~\ref{lem:CLT-Stein-lemma}  to prove this lemma.  Let $\mathsf p_k$, $\mathsf q_k$ be the corresponding functions in defining $\mathsf P_k$ and $\mathsf Q_k$:
    \begin{align}
         \mathsf p_k(x_1,\dots,x_n)=\frac{1}{\sqrt{n}}\sum_{j=1}^n h_k( x_j) \quad \text{ and } \quad \mathsf q_k(y_1,\dots,y_n)=\frac{1}{\sqrt{n}}\sum_{j=1}^n h_k( y_j)\,.\notag
    \end{align}
    Note that $H'_k(y)=kH_{k-1}(y)$ and $h_k(y)=H_k(y)/\sqrt{k!}$, we have (similarly for $\mathsf q_k$)
    \begin{align}
        \frac{\partial \mathsf p_k}{\partial x_j}=\frac{\sqrt{k}}{\sqrt{n}}h_{k-1}(x_j) \quad \text{ and } \quad \frac{\partial \mathsf p_k}{\partial y_j}=0, \quad \text{ for all }1\leq j\leq n\,.\label{eq-partial-derivatives-pq}
    \end{align}
   Let $W=\{(\mathsf P_k,\mathsf Q_k)\}_{k=1}^M$, $Z=\{(\mathsf U_k,\mathsf V_k)\}_{k=1}^M$. Then $\{(\mathsf X_i,\mathsf Y_i)\}_{i=1}^n  \sim N(\mathbf 0,\mathbf K)$ 
   and $Z \sim  N(\mathbf 0,\mathbf C),$ where
    \begin{equation}
    \begin{split}
        \mathbf K = \mathrm{diag}\left\{
          \underbrace{\begin{pmatrix} 1 & \rho \\ \rho & 1 \end{pmatrix},\dots,
          \begin{pmatrix} 1 & \rho \\ \rho & 1 \end{pmatrix}}_{n\ \text{times}}
          \right\}\,, \quad
       \mathbf C  =\mathrm{diag}\left\{
          \underbrace{\begin{pmatrix} 1 & \rho \\ \rho & 1 \end{pmatrix},\begin{pmatrix} 1 & \rho^2 \\ \rho^2 & 1 \end{pmatrix},\dots,
          \begin{pmatrix} 1 & \rho^M \\ \rho^M & 1 \end{pmatrix}}_{M\ \text{times}}
          \right\}\,.
    \end{split}
    \end{equation}
    In view of~\prettyref{lem:CLT-Stein-lemma}, to bound the $W_1$ distance between $W$ and $Z$, it remains to  evaluate $T_{ab}$ as defined in~\prettyref{eq:T_a_b} for all $1\leq a, b \leq 2M$. We separately consider the following four cases: (i) $a,b$ are odd; (ii)  $a,b$ are even; (iii) $a$ is odd, $b$ is even; (iv) $a$ is even, $b$ is odd. Note that by Lemma~\ref{lem:scaling-law-hermite}, 
    \begin{equation}
    \begin{split}
        & \mathbb E\left[ \frac{\partial\mathsf p_r}{\partial x_j}\left(\sqrt{t}x+\sqrt{1-t}\mathsf X\right) \right]=\sqrt{\frac{r}{n}} \mathbb{E}\left[ h_{r-1} \left(\sqrt{t}x_j+\sqrt{1-t}\mathsf X_j\right)\right] 
        =  \sqrt{\frac{r}{n}} t^{\frac{r-1}{2}}h_{r-1}(x_j)\,.\notag
    \end{split}
    \end{equation} 
    For case (i), suppose $a=2k-1,b=2k+2s-1$, then 
    \begin{equation}\label{eq-case1}
    \begin{split}
        T_{ab}(x,y) &= \int_0^1 \frac{1}{2\sqrt{t}}\sum_{j=1}^n\sqrt{\frac{k}{n}}h_{k-1}(x_j) \sqrt{\frac{k+s}{n}} t^{\frac{k+s-1}{2}} h_{k+s-1}(x_j)\mathrm{d}t\\
        &=\sqrt{\frac{a+1}{b+1}}\frac{1}{n}\sum_{j=1}^n h_{(a+1)/2-1}(x_j)h_{(b+1)/2-1}(x_j)\,.
    \end{split}
    \end{equation}
   For case (ii), similarly we suppose $a=2k,b=2k+2s $, then 
    \begin{equation}\label{eq-case2}
    \begin{split}
        T_{ab}(x,y) &= \int_0^1 \frac{1}{2\sqrt{t}}\sum_{j=1}^n\sqrt{\frac{k}{n}}h_{k-1}(y_j) \sqrt{\frac{k+s}{n}} t^{\frac{k+s-1}{2}} h_{k+s-1}(y_j)\mathrm{d}t\\
        &=\sqrt{\frac{a}{b}}\frac{1}{n}\sum_{j=1}^n h_{a/2-1}(y_j)h_{b/2-1}(y_j)\,.
    \end{split}
    \end{equation}
    For case (iii), if $a=2k-1,b=2k+2s$,
    \begin{equation}
        \label{eq-case3}
        \begin{split}
            T_{ab}(x,y)&=\int_{0}^1\frac{1}{2\sqrt{t}} \sum_{j=1}^n \rho \sqrt{\frac{k}{n}}h_{k-1}(x_j)\times \sqrt{\frac{k+s}{n}}t^{\frac{k+s-1}2}h_{k+s-1} (y_j) \mathrm{d}t\\
            &=\frac{\rho}{n} \sqrt{\frac{a+1}{b}} \sum_{j=1}^n h_{(a+1)/2-1}(x_j)h_{b/2-1}(y_j)\,.
        \end{split}
    \end{equation}
    Similarly for case (iv), suppose $a=2k,b=2k+2s-1$,
    \begin{equation}
        \label{eq-case4}
        \begin{split}
            T_{a,b}(x,y)=\frac{\rho}{n}\sqrt{\frac{a}{b+1}} \sum_{j=1}^n h_{a/2-1}(y_j)h_{(b+1)/2-1}(x_j)
        \end{split}
    \end{equation}
    Combining \eqref{eq-case1}, \eqref{eq-case2}, \eqref{eq-case3} and \eqref{eq-case4} yields
    \begin{equation}
        \label{eq-Tab-final-formula}
        T_{ab}(x,y) = 
        \begin{cases}
            & \frac{1}{n}\sqrt{\frac{a+1}{b+1}}\sum_{j=1}^nh_{(a+1)/2-1}(x_j)h_{(b+1)/2-1}(x_j),  \text{ if }a,b\text{ are odd}\,;\\
             & \frac{1}{n}\sqrt{\frac{a}{b}}\sum_{j=1}^nh_{a/2-1}(y_j)h_{b/2-1}(y_j),  \text{ if }a,b\text{ are even}\,;\\
            & \frac{\rho}{n} \sqrt{\frac{a+1}{b}} \sum_{j=1}^n h_{(a+1)/2-1}(x_j)h_{b/2-1}(y_j),\text{ if $a$ is odd, $b$ is even}\,;\\
            & \frac{\rho}{n}\sqrt{\frac{a}{b+1}} \sum_{j=1}^n h_{a/2-1}(y_j)h_{(b+1)/2-1}(x_j),\text{ if $a$ is even, $b$ is odd}\,.
        \end{cases}
    \end{equation}
    Next we compute $\mathbb E\left( C(a,b) - T_{ab}(\mathsf X,\mathsf Y) \right)^2$ for each $1\leq a,b\leq 2M$. We first consider the case where $a=b$ and $a$ is odd. We have
    \begin{equation}
        \begin{split}
            & \mathbb E\left[ \left( C(a,a) - T_{aa}(\mathsf X,\mathsf Y)\right)^2 \right] = \mathbb E\left[ \left( 1 - \frac{1}{n}\sum_{j=1}^n h^2_{(a+1)/2-1}\left(\mathsf X_j
            \right)\right)^2\right]\\
          = & \frac{1}{n^2}\mathbb E\left[ \left( \sum_{j=1}^n \left(1-h^2_{(a+1)/2-1}(\mathsf X_j)\right) \right)^2 \right]=\frac{1}{n}\mathbb E\left[\left(1-h_{(a+1)/2-1}^2(\mathsf X_1)\right)^2\right]\\
          = & \frac{1}{n}\left(\mathbb E\left[h^4_{(a+1)/2-1}(\mathsf X_1)\right]-1\right)\,.
        \end{split}\notag
    \end{equation}
    Similarly, for $a=b$ and $a$ being even,
    \begin{equation}
        \mathbb E\left[ \left( C(a,a) - T_{aa}(\mathsf X,\mathsf Y)\right)^2 \right] = \frac{1}{n}\left(\mathbb E\left[h^4_{a/2-1}(\mathsf Y_1)\right]-1\right)\,.
    \end{equation}
    If $a\neq b$ are odd numbers,
    \begin{equation}
    \begin{split}
        & \mathbb E\left[ \left( C(a,b) - T_{ab}(\mathsf X,\mathsf Y)\right)^2 \right] = \mathbb E\left[ \left( T_{ab}(\mathsf X,\mathsf Y) \right)^2 \right]\\
      = &\mathbb E\left[ \left( \frac{1}{n}\sqrt{\frac{a+1}{b+1}}\sum_{j=1}^nh_{(a+1)/2-1}(\mathsf X_j)h_{(b+1)/2-1}(\mathsf X_j) \right)^2 \right]= \frac{1}{n^2}\frac{a+1}{b+1} \sum_{j=1}^n\mathbb E\left[h^2_{(a+1)/2-1}(\mathsf X_j)h^2_{(b+1)/2-1}(\mathsf X_j)\right]\\
      = & \frac{1}{n}\frac{a+1}{b+1}\mathbb E\left[ h^2_{(a+1)/2-1}(\mathsf X_1)h^2_{(b+1)/2-1}(\mathsf X_1) \right] \leq \frac{1}{n}\frac{a+1}{b+1}\left(\mathbb E\left[ h^4_{(a+1)/2-1}(\mathsf X_1) \right]\right)^{\frac{1}{2}}\left(\mathbb E\left[ h^4_{(b+1)/2-1}(\mathsf X_1) \right]\right)^{\frac{1}{2}}\,.
    \end{split}\notag
    \end{equation}
    Similarly if $a\neq b$ are even numbers,
    \begin{equation}
        \mathbb E\left[ \left( C(a,b) - T_{ab}(\mathsf X,\mathsf Y)\right)^2 \right] \leq  \frac{1}{n}\frac{a}{b}\left(\mathbb E\left[ h^4_{a/2-1}(\mathsf Y_1) \right]\right)^{\frac{1}{2}}\left(\mathbb E\left[ h^4_{(b/2-1}(\mathsf Y_1) \right]\right)^{\frac{1}{2}}\,.
    \end{equation}
    When $a=2k-1$ and $b=2k$, we have
    \begin{equation}
        \begin{split}
           & \mathbb E\left[ \left( C(a,b) - T_{ab}(\mathsf X,\mathsf Y)\right)^2 \right] = \mathbb E\left[ \left(\rho^{k} - \frac{\rho}{n} \sum_{j=1}^n h_{k-1}(\mathsf X_j)h_{k-1}(\mathsf Y_j)\right)^2\right]\\
           = & \frac{1}{n^2}\Bigg(\sum_{i\neq j}\mathbb E[\left(\rho^{k}-\rho h_{k-1}(\mathsf X_i)h_{k-1}(\mathsf Y_i)\right)]\mathbb E[\left(\rho^{k}-\rho h_{k-1}(\mathsf X_j)h_{k-1}(\mathsf Y_j)\right)]\\
           &\qquad \quad+\sum_{j=1}^n\mathbb E\left[ \left(\rho^{k}-\rho h_{k-1}(\mathsf X_j)h_{k-1}(\mathsf Y_j)\right)^2 \right]\Bigg) \\
           = &\frac{1}{n^2}\sum_{j=1}^n \left(\rho^2\mathbb E[h^2_{k-1}(\mathsf X_j)h^2_{k-1}(\mathsf Y_j)] - 2\rho^{k+1}\cdot \rho^{k-1}+\rho^{2k}\right)\\
           \leq & \frac{1}{n}\left(\mathbb E[h^4_{( a+1)/2-1}(\mathsf X_1)]\right)^{1/2}\left(\mathbb E[h^4_{(a+1)/2-1}(\mathsf Y_1)]\right)^{1/2}\,.
        \end{split}
    \end{equation}
    And we also have that for $ a=2k$ and $b=2k-1$,
    \begin{equation}
        \mathbb E\left[ \left( C(a,b) - T_{ab}(\mathsf X,\mathsf Y)\right)^2 \right]\leq \frac{1}{n}\left(\mathbb E[h^4_{a/2-1}(\mathsf X_1)]\right)^{1/2}\left(\mathbb E[h^4_{a/2-1}(\mathsf Y_1)]\right)^{1/2}\,.
    \end{equation}
    If $a=2k-1,b=2k+2s$ where $s\neq0$,
    \begin{equation}
        \begin{split}
           & \mathbb E\left[ \left( C(a,b) - T_{ab}(\mathsf X,\mathsf Y)\right)^2 \right] = \mathbb E\left[ \left( T_{ab}(\mathsf X,\mathsf Y) \right)^2 \right]\\
           =&\mathbb E\left[\left(\frac{\rho}{n} \sqrt{\frac{a+1}{b}} \sum_{j=1}^n h_{(a+1)/2-1}(\mathsf X_j)h_{b/2-1}(\mathsf Y_j)\right)^2\right] =\frac{\rho^2}{n^2} \frac{a+1}{b}\sum_{j=1}^n \mathbb E\left[  h^2_{(a+1)/2-1}(\mathsf X_j)h^2_{b/2-1}(\mathsf Y_j)\right]\\
           =&\frac{\rho^2}{n} \frac{a+1}{b} \mathbb E\left[  h^2_{(a+1)/2-1}(\mathsf X_1)h^2_{b/2-1}(\mathsf Y_1)\right]\leq \frac{a+1}{bn}\left(\mathbb E\left[h^4_{(a+1)/2-1}(\mathsf X_1)\right]\right)^{\frac{1}{2}}\left(\mathbb E\left[h^4_{b/2-1}(\mathsf X_1)\right]\right)^{\frac{1}{2}}\,.
        \end{split}
    \end{equation}
    Similarly if $a=2k,b=2k+2s-1$ where $s\neq 0$,
    \begin{equation}
        \mathbb E\left[ \left( C(a,b) - T_{ab}(\mathsf X,\mathsf Y)\right)^2 \right]\leq \frac{a}{(b+1)n}\left(\mathbb E\left[h^4_{a/2-1}(\mathsf X_1)\right]\right)^{\frac{1}{2}}\left(\mathbb E\left[h^4_{(b+1)/2-1}(\mathsf X_1)\right]\right)^{\frac{1}{2}}\,.
    \end{equation}
    Note that $\|C\|_{\operatorname{op}}=\max\{1+\rho,1-\rho\},\|C^{-1}\|_{\operatorname{op}} = \max\{\frac{1}{1+\rho},\frac{1}{1-\rho}\}$, by Lemmas~\ref{lem:CLT-Stein-lemma} and Lemma~\ref{lem:fourth-moment-hermite-estimate}, we have
    \begin{equation}
        \mathcal{W}^2_1\left(W,Z\right) \leq \frac{1+|\rho|}{(1-|\rho|)^2}\frac{1}{n}\exp(cM)\,,
    \end{equation}
    where $c$ is an absolute number independent of $\rho$ and $n$.
\end{proof}

\section*{Acknowledgment}
The bulk of this work was developed in 2024--2025.
All ideas and proofs in this paper were fully human-generated, and the manuscript was written by the authors. AI tools were used only for copy editing, sanity check, and assisting with the simulation in \prettyref{sec:numerical-results}.

Part of the work was supported by the National Science Foundation under
Grant No. DMS-1928930, while the authors were in residence at the Simons Laufer Mathematical Sciences
Institute in Berkeley, California, during the Spring 2025 semester. J. Xu is supported in part by
an NSF CAREER award CCF-2144593.

S.~Gong would like to thank Zhangsong Li for discussions on the low-degree hardness aspect of this model.

 \bibliographystyle{alpha}
\bibliography{references,references1} 

\newcommand{\etalchar}[1]{$^{#1}$}
\begin{thebibliography}{WWXY22}

\bibitem[APZ17]{abid2017linear}
Abubakar Abid, Ada Poon, and James Zou.
\newblock Linear regression with shuffled labels.
\newblock {\em arXiv preprint arXiv:1705.01342}, 2017.

\bibitem[BH05]{bai2005broken}
Zhidong Bai and Tailen Hsing.
\newblock The broken sample problem.
\newblock {\em Probability theory and related fields}, 131:528--552, 2005.

\bibitem[BK13]{berend2013sharp}
Daniel Berend and Aryeh Kontorovich.
\newblock A sharp estimate of the binomial mean absolute deviation with
  applications.
\newblock {\em Statistics \& Probability Letters}, 83(4):1254--1259, 2013.

\bibitem[CMK18]{cullina2018fundamental}
Daniel Cullina, Prateek Mittal, and Negar Kiyavash.
\newblock Fundamental limits of database alignment.
\newblock In {\em 2018 IEEE International Symposium on Information Theory
  (ISIT)}, pages 651--655. IEEE, 2018.

\bibitem[DDCP88]{delong1988comparing}
Elizabeth~R. DeLong, David~M. DeLong, and Daniel~L. Clarke-Pearson.
\newblock Comparing the areas under two or more correlated receiver operating
  characteristic curves: A nonparametric approach.
\newblock {\em Biometrics}, 44(3):837--845, 1988.

\bibitem[DFG71]{degroot1971matchmaking}
Morris~H DeGroot, Paul~I Feder, and Prem~K Goel.
\newblock Matchmaking.
\newblock {\em The Annals of Mathematical Statistics}, pages 578--593, 1971.

\bibitem[DG80]{degroot1980estimation}
Morris~H DeGroot and Prem~K Goel.
\newblock Estimation of the correlation coefficient from a broken random
  sample.
\newblock {\em The Annals of Statistics}, pages 264--278, 1980.

\bibitem[DMWX21]{ding2021efficient}
Jian Ding, Zongming Ma, Yihong Wu, and Jiaming Xu.
\newblock Efficient random graph matching via degree profiles.
\newblock {\em Probability Theory and Related Fields}, 179:29--115, 2021.

\bibitem[EH23]{elimelech2023phase}
Dor Elimelech and Wasim Huleihel.
\newblock Phase transitions in the detection of correlated databases.
\newblock In {\em International Conference on Machine Learning}, pages
  9246--9266. PMLR, 2023.

\bibitem[FS69]{fellegi1969theory}
Ivan~P Fellegi and Alan~B Sunter.
\newblock A theory for record linkage.
\newblock {\em Journal of the American statistical association},
  64(328):1183--1210, 1969.

\bibitem[FSWW26]{fan2026empirical}
Zhou Fan, Yandi Shen, Haoyu Wang, and Yihong Wu.
\newblock Empirical {Bayes} linear regression in high dimensions: Method of
  moments and sub-linear sample complexity, 2026.
\newblock \href{https://arxiv.org/abs/2608.16771}{arXiv:2608.16771}.

\bibitem[Gon10]{gonzalez2010measurement}
{\'A}lvaro Gonz{\'a}lez.
\newblock Measurement of areas on a sphere using fibonacci and
  latitude--longitude lattices.
\newblock {\em Mathematical Geosciences}, 42(1):49--64, 2010.

\bibitem[GZ22]{gamarnik2022sparse}
David Gamarnik and Ilias Zadik.
\newblock Sparse high-dimensional linear regression. estimating squared error
  and a phase transition.
\newblock {\em The Annals of Statistics}, 50(2):880--903, 2022.

\bibitem[Hop18]{hopkins2018thesis}
Samuel Hopkins.
\newblock Statistical inference and the sum of squares method.
\newblock {\em Cornell University}, 2018.

\bibitem[JWX25]{jiao2025DetectionBroken}
Simiao Jiao, Yihong Wu, and Jiaming Xu.
\newblock The broken sample problem revisited: Proof of a conjecture by
  bai-hsing and high-dimensional extensions.
\newblock {\em arXiv preprint arXiv:2503.14619. To appear in Bernoulli}, 2025.

\bibitem[KB80]{kaas1980mean}
Rob Kaas and Jan~M Buhrman.
\newblock Mean, median and mode in binomial distributions.
\newblock {\em Statistica Neerlandica}, 34(1):13--18, 1980.

\bibitem[KBSB08]{kaplan2008assignment}
Lance~M. Kaplan, Yaakov Bar-Shalom, and William~D. Blair.
\newblock Assignment costs for multiple sensor track-to-track association.
\newblock {\em IEEE Transactions on Aerospace and Electronic Systems},
  44(2):655--677, 2008.

\bibitem[Li26]{Zhangsong25+}
Zhangsong Li.
\newblock A computational transition for detecting multivariate shuffled linear
  regression by low-degree polynomials.
\newblock {\em IEEE Transactions on Information Theory}, 2026.

\bibitem[LS59]{leonov1959method}
V.~P. Leonov and A.~N. Shiryaev.
\newblock On a method of calculation of semi-invariants.
\newblock {\em Theory of Probability and Its Applications}, 4(3):319--329,
  1959.

\bibitem[MBR{\etalchar{+}}13]{marques2013spherical}
Ricardo Marques, Christian Bouville, Micka{\"e}l Ribardi{\`e}re, Lu{\'i}s~Paulo
  Santos, and Kadi Bouatouch.
\newblock Spherical fibonacci point sets for illumination integrals.
\newblock {\em Computer Graphics Forum}, 32(8):134--143, 2013.

\bibitem[MS06]{muralidhar2006data}
Krishnamurty Muralidhar and Rathindra Sarathy.
\newblock Data shuffling---a new masking approach for numerical data.
\newblock {\em Management Science}, 52(5):658--670, 2006.

\bibitem[MW23]{mazumder2023linear}
Rahul Mazumder and Haoyue Wang.
\newblock Linear regression with partially mismatched data: local search with
  theoretical guarantees.
\newblock {\em Mathematical Programming}, 197(2):1265--1303, 2023.

\bibitem[NPR10]{nourdin2010multivariate}
Ivan Nourdin, Giovanni Peccati, and Anthony R{\'e}veillac.
\newblock Multivariate normal approximation using stein's method and malliavin
  calculus.
\newblock In {\em Annales de l'IHP Probabilit{\'e}s et statistiques},
  volume~46, pages 45--58, 2010.

\bibitem[OWN97]{oppenheim1997signals}
Alan~V. Oppenheim, Alan~S. Willsky, and S.~Hamid Nawab.
\newblock {\em Signals and Systems}.
\newblock Prentice Hall, 2nd edition, 1997.

\bibitem[PH24]{paslev2023testing}
Vered Paslev and Wasim Huleihel.
\newblock Testing dependency of unlabeled databases.
\newblock {\em IEEE Transactions on Information Theory}, 2024.

\bibitem[PT11]{peccati2011wiener}
Giovanni Peccati and Murad~S. Taqqu.
\newblock {\em Wiener Chaos: Moments, Cumulants and Diagrams}, volume~1 of {\em
  Bocconi \& Springer Series}.
\newblock Springer, Milan, 2011.

\bibitem[PW25]{polyanskiy2025information}
Yury Polyanskiy and Yihong Wu.
\newblock {\em Information theory: From coding to learning}.
\newblock Cambridge university press, 2025.

\bibitem[PWC17]{pananjady2017linear}
Ashwin Pananjady, Martin~J Wainwright, and Thomas~A Courtade.
\newblock Linear regression with shuffled data: Statistical and computational
  limits of permutation recovery.
\newblock {\em IEEE Transactions on Information Theory}, 64(5):3286--3300,
  2017.

\bibitem[RXZ19]{reeves2019all}
Galen Reeves, Jiaming Xu, and Ilias Zadik.
\newblock The all-or-nothing phenomenon in sparse linear regression.
\newblock In {\em Conference on Learning Theory}, pages 2652--2663. PMLR, 2019.

\bibitem[SBD19]{slawski2019linear}
Martin Slawski and Emanuel Ben-David.
\newblock Linear regression with sparsely permuted data.
\newblock {\em Electronic Journal of Statistics}, 13(1):1--36, 2019.

\bibitem[SBDL20]{slawski2020two}
Martin Slawski, Emanuel Ben-David, and Ping Li.
\newblock Two-stage approach to multivariate linear regression with sparsely
  mismatched data.
\newblock {\em The Journal of Machine Learning Research}, 21(1):8422--8463,
  2020.

\bibitem[Sta13]{stanley2013algebraic}
Richard~P Stanley.
\newblock Algebraic combinatorics.
\newblock {\em Springer}, 20(22):4, 2013.

\bibitem[UHV18]{unnikrishnan2018unlabeled}
Jayakrishnan Unnikrishnan, Saeid Haghighatshoar, and Martin Vetterli.
\newblock Unlabeled sensing with random linear measurements.
\newblock {\em IEEE Transactions on Information Theory}, 64(5):3237--3253,
  2018.

\bibitem[Ver18]{vershynin2018high}
Roman Vershynin.
\newblock {\em High-dimensional probability: An introduction with applications
  in data science}, volume~47.
\newblock Cambridge university press, 2018.

\bibitem[WWXY22]{wang2022random}
Haoyu Wang, Yihong Wu, Jiaming Xu, and Israel Yolou.
\newblock Random graph matching in geometric models: the case of complete
  graphs.
\newblock In {\em Conference on Learning Theory}, pages 3441--3488. PMLR, 2022.

\bibitem[WZW{\etalchar{+}}21]{wang2021multisensor}
Jun Wang, Yajun Zeng, Shaoming Wei, Zixiang Wei, Qinchen Wu, and Yvon Savaria.
\newblock Multi-sensor track-to-track association and spatial registration
  algorithm under incomplete measurements.
\newblock {\em IEEE Transactions on Signal Processing}, 69:3337--3350, 2021.

\bibitem[XT06]{xiao2006anatomy}
Xiaokui Xiao and Yufei Tao.
\newblock Anatomy: Simple and effective privacy preservation.
\newblock In {\em Proceedings of the 32nd International Conference on Very
  Large Data Bases}, pages 139--150, 2006.

\bibitem[ZCB{\etalchar{+}}23]{zhu2023robust}
Bokai Zhu, Shuxiao Chen, Yunhao Bai, Han Chen, Guanrui Liao, Nilanjan
  Mukherjee, Gustavo Vazquez, David~R. McIlwain, Alexandar Tzankov, Ivan~T.
  Lee, Matthias~S. Matter, Yury Goltsev, Zongming Ma, Garry~P. Nolan, and Sizun
  Jiang.
\newblock Robust single-cell matching and multimodal analysis using shared and
  distinct features.
\newblock {\em Nature Methods}, 20(2):304--315, 2023.

\bibitem[ZL20]{zhang2020optimal}
Hang Zhang and Ping Li.
\newblock Optimal estimator for unlabeled linear regression.
\newblock In {\em International Conference on Machine Learning}, pages
  11153--11162. PMLR, 2020.

\bibitem[ZSL19]{zhang2019permutation}
Hang Zhang, Martin Slawski, and Ping Li.
\newblock Permutation recovery from multiple measurement vectors in unlabeled
  sensing.
\newblock In {\em 2019 IEEE International Symposium on Information Theory
  (ISIT)}, pages 1857--1861. IEEE, 2019.

\end{thebibliography}

\end{document}